\documentclass[a4paper,reqno,12pt]{amsart}
\usepackage[utf8]{inputenc}
\usepackage{amssymb}
\usepackage{amsxtra}
\usepackage{ stmaryrd }
\usepackage{amssymb,amsmath,amsthm,stmaryrd}
\usepackage{color}
\usepackage[all]{xy}
\usepackage{hyperref}
\usepackage{thmtools}
\usepackage{appendix}
\usepackage{cancel}
\usepackage[tracking,spacing]{microtype}
\theoremstyle{plain}

\declaretheorem{theorem}
\declaretheorem{lemma}
\declaretheorem{proposition}
\declaretheorem{remark}
\declaretheorem{example}
\declaretheorem{corollary}
\declaretheorem{definition}

\newcommand{\source}   {\mathsf{s}}

\newcommand{\tg}       {\mathsf t} 
\newcommand{\s}        {\mathsf s} 

\newcommand{\pt}{*}

\newcommand{\R}{\mathbb{R}}

 \DeclareMathOperator{\pr}{pr}
\DeclareMathOperator{\Hom}{Hom}

\DeclareMathOperator{\End}{End}

\DeclareMathOperator{\id}{id}

\DeclareMathOperator{\holim}{holim}

\newcommand*{\emptycomment}[1]{}

\DeclareMathOperator{\coker}{coker}
\DeclareMathOperator{\Aut}{Aut}

\newcommand*{\horiprod}{\mathbin{\cdot_\textup h}}
\newcommand*{\vertprod}{\mathbin{\cdot_\textup v}}
\newcommand*{\horicirc}{\mathbin{\circ_\textup h}}

\newcommand{\huaB}{\mathcal{B}}

\newcommand{\huaA}{\mathcal{A}}

\newcommand{\huaR}{\mathcal{R}}
\newcommand{\huaE}{\mathcal{E}}
\newcommand{\huaF}{\mathcal{F}}
\newcommand{\huaG}{\mathcal{G}}

\newcommand{\huaV}{\mathcal{V}}

\newcommand{\huaX}{\mathcal{X}}

\newcommand{\huaP}{\mathcal{P}}
\newcommand{\huaC}{{\mathcal{C}}}

\newcommand{\huaH}{\mathcal{H}}

\newcommand{\huaT}{\mathcal{T}}

\newcommand{\Chi}{\mathbb{X}}

\newcommand{\Rep}{\textsf{Rep}}
\newcommand{\Vect}{\textsf{Vect}}
\newcommand{\VectBd}{\textsf{VectBd}}
\newcommand{\GBd}{\huaG\textsf{Bd}}
\newcommand{\ngrd}{n\textsf{Grpd}}
\newcommand{\Auto}{\textsf{Auto}} 
\newcommand{\gpd}{\textsf{Gpd}} 

\DeclareMathOperator{\Fun}{Fun} 
\DeclareMathOperator{\pFun}{PsFun} 
\DeclareMathOperator{\smap}{\mathrm{Map}}

\newcommand{\Z}{\ensuremath{\mathbb Z}}

\newcommand{\Mfd}{\mathsf{Mfd}}
\newcommand{\bgpd}{\textsf{Bibun}} 

\newcommand{\String}{\textup{String}}
\newcommand{\faith}{\mathsf{faith}}

\newcommand{\sTran}{\textup{sTran}}

\newcommand{\thuaV}{\Tilde{\huaV}}
\newcommand{\tgamma}{\Tilde{\gamma}}
\newcommand{\tphi}{\Tilde{\phi}}
\newcommand{\talpha}{\Tilde{\alpha}}
\newcommand{\tbeta}{\Tilde{\beta}}
\newcommand{\tW}{\Tilde{W}}
\newcommand{\tU}{\Tilde{U}}
\newcommand{\tV}{\Tilde{V}}
\newcommand{\thuaE}{\Tilde{\huaE}}
\newcommand{\tfrakX}{\tilde{\mathfrak{X}}}

\begin{document}

\title{2-Representations of Lie 2-groups and 2-Vector Bundles}
\author[Z. Huan]{Zhen Huan}
\address{Zhen Huan, Center for Mathematical Sciences, Huazhong University of Science and Technology, Hubei 430074, China} \curraddr{}
\email{2019010151@hust.edu.cn}
\date{\today}


\begin{abstract}
We develop a systematic 2-representation theory for coherent Lie 2-groups and its geometric incarnation via 2-vector bundles. For any coherent Lie 2-group $\huaG_{\bullet}$, we construct a bicategory $2\Rep(\huaG_{\bullet})$ of 2-representations valued in $k$-prelinear categories, and prove that biequivalent 2-groups induce biequivalent 2-representation bicategories. When $\huaG_{\bullet}$ arises from an action groupoid $H \ltimes G \rightrightarrows G$, we construct a functor 
from the category of faithful $H$-representations  to $2\Rep(\huaG_{\bullet})$ 
that is injective on objects.  Applied to the string 2-group $\String_k(G)$, this yields, for each positive energy representation of the loop group, an associated 2-representation of the string 2-group. Geometrically, we construct 2-stacks of 2-vector bundles and principal 2-bundles, and establish an associated bundle construction as a strong transformation of 2-stacks. We further extend this to an equivariant framework over Lie groupoids, providing a foundation for higher gauge theories with symmetries.
\end{abstract}

\maketitle

\tableofcontents

\section{Introduction}

\subsection{Background and motivation}

Symmetries in classical gauge theory are elegantly encoded by principal bundles and their associated vector bundles. When we attempt to categorify this framework---replacing Lie groups by Lie 2-groups and vector spaces by 2-vector spaces---we enter a world where parallel transport is implemented by bifunctors, and gauge symmetries are governed by 2-group actions. It is a necessary step for understanding the geometry of string theory, higher gauge fields, and extended topological quantum field theories. A schematic comparison is given in Figure  \ref{cpr_1vect_2vect}.

\begin{figure}[t]
\begin{center}
\renewcommand{\arraystretch}{1.3}
\begin{tabular}{|c || c | c |}
\hline
\textbf{Feature} & \textbf{1-vector bundle} & \textbf{2-vector bundle} \\
\hline \hline
Fibers & Vector spaces & 2-vector spaces \\
\hline
Parallel Transport & Linear maps & Linear functors \& \\
 && natural transformations \\
\hline
Symmetries & Group actions & 2-group actions \\
\hline
Twists \& obstructions & $H^1(X; G)$ & $H^2(X; \mathcal{G}_\bullet)$ and gerbes \\
\hline
Physical realizations & Gauge theory, QM & Higher gauge theory, TQFT \\
\hline 
\end{tabular}
\caption{Comparing 1-vector bundles and 2-vector bundles}
\label{cpr_1vect_2vect}
\end{center}
\end{figure}

Despite growing interest in higher representation theory \cite{barrett2006categorical, BDSV2015, baez2012infinite, Bunk:2016rta, DR2018, KLW2Rep22, murray-robert-wockel}, an obstruction arises in the study of one of the most central objects: the string 2-group. As demonstrated by Murray, Roberts, and Wockel \cite{murray-robert-wockel}, the natural loop group model for the string 2-group---the one that  encodes the circle action on loops---\emph{cannot} be realized as a strict Lie 2-group. It is forced to be a \emph{coherent} Lie 2-group, where the group laws hold only up to higher coherent isomorphisms. This is not a defect; it is a fundamental feature of the geometry. Any   2-representation theory for the string 2-group must therefore handle coherence conditions from the start.

\subsection{Main results}

This paper develops a systematic 2-representation theory for  coherent Lie 2-groups, with the string 2-group as its main application. Our contributions are listed as follows.

\begin{enumerate}
    \item \textbf{A bicategory of 2-representations}. For a coherent Lie 2-group $\huaG_{\bullet}$, we construct a bicategory $2\Rep(\huaG_{\bullet})$ whose objects are pseudofunctors $B\huaG_{\bullet} \to 2\Vect$.  Here $2\Vect$ is the 2-category of $k$-prelinear categories. While several models of 2-vector spaces exist in the literature---Kapranov--Voevodsky model \cite{KapranovVoevodsky1994}, the Baez--Crans model \cite{baezcrans}, 2-Hilbert spaces \cite{Baez_hilbert_1996} and others. We work with $k$-prelinear categories, following  the frameworks of \cite{douglas-reutter:18, BDSV2015, Lurie08_TFT} adapted to the smooth setting.  All coherence conditions  from the 2-group structure are handled explicitly.

    \item \textbf{Relation between 1-representations and 2-representations}. For a Lie 2-group $\huaG_{\bullet}$ whose underlying Lie groupoid is an action groupoid $H \ltimes G \rightrightarrows G$, we construct a functor from  the  category $\Rep^{\faith} H$ of faithful representations of $H$ into the 2-representation bicategory $2\Rep(\huaG_{\bullet})$. Moreover, every representation $(V, \zeta)$ of $H$ gives a 2-representation of $\huaG_{\bullet}$ via the functor category $\Fun(\huaG_{\bullet}, \{V\}/\zeta(H))$, and intertwiners lift to 1-morphisms. This bridges the classical and categorified theories directly.

    \item \textbf{Biequivalence invariance}. We prove that if two coherent Lie 2-groups are biequivalent---the correct notion of equivalence for 2-groups---then their bicategories of 2-representations are biequivalent:
    \[ \huaG_{\bullet} \simeq_{\text{bieq}} \huaH_{\bullet} \quad\Longrightarrow\quad 2\Rep (\huaG_{\bullet} )\simeq_{\text{bieq}} 2\Rep(\huaH_{\bullet}). \]
    This result guarantees that our 2-representation theory is a genuine invariant of the 2-group's homotopy type, liberating us to work with whichever model is most convenient---strict or coherent.

    \item \textbf{Geometric incarnation: associated 2-vector bundles and equivariant models}. On the geometric side, we build a 2-stack $2\VectBd^+$ of 2-vector bundles over smooth manifolds and a 2-stack $\GBd^+$ of principal $\huaG_{\bullet}$-2-bundles. For each $\huaG_{\bullet}$-2-representation $\rho$, we establish an \emph{associated bundle construction} as a strong transformation of 2-stacks $\rho_* : \GBd^+ \Rightarrow 2\VectBd^+$, sending a principal 2-bundle to an associated 2-vector bundle. We then extend the framework of 2-vector bundles to the \emph{equivariant setting} over Lie groupoids, defining the bicategory $2\VectBd^+_{\huaG_{\bullet}}(X_\bullet)$ of $\huaG_{\bullet}$-equivariant 2-vector bundles over $X_\bullet$. This provides a natural geometric language for higher gauge theories with symmetries.

    \item \textbf{Explicit examples}. We construct several concrete, non-trivial examples of 2-representations, verifying all coherence conditions in the bicategorical setting. These examples demonstrate that the theory is far from empty and illustrate the rich additional structure that emerges in the coherent case compared to the strict one.
\end{enumerate}
Taken together, these results establish 2-representation theory as a bicategorical invariant of Lie 2-groups, and provide a systematic geometric framework—from algebraic representations to associated bundles—for higher gauge theory and string geometry.

\subsection{Comparison with prior work}

The string 2-group has been realized in various geometric models \cite{baez:str-gp, schommer:string-finite-dim, murray-robert-wockel}; we work with the coherent model of Murray--Roberts--Wockel, which is based on the free loop group. For 2-representations, our framework is complementary to the analytic approach of Kristel--Ludewig--Waldorf \cite{KLW2Rep22}, focusing instead on an algebraic and geometric model valued in $k$-prelinear categories. A more detailed discussion of related work is deferred to Section \ref{app_add}.

\subsection*{Methodological remark}
The framework developed in this paper encodes smoothness and continuity not via 
analytic conditions on representation maps, but via the higher categorical 
structure of stacks and pseudofunctors.  The coherence conditions of pseudofunctors (Section 
\ref{def:2rep_bicat}) and the plus construction (Section 
\ref{2vbconstr}) together provide the gluing and compatibility conditions 
that, in an analytic approach such as \cite{KLW2Rep22}, would be imposed on 
Hilbert space structures. 
This categorical approach is not an approximation to the analytic one, but a distinct perspective that captures the algebraic skeleton of higher representation theory. The two perspectives are complementary.

\subsection{Outline of the paper}

\begin{itemize}
    \item \textbf{Section \ref{2grp_basic}} lays the groundwork: we recall the language of coherent Lie 2-groups and define the bicategory $2\Rep(\huaG_{\bullet})$.
    \item \textbf{Sections \ref{ex2rep} and \ref{ex_2rep_cr}} bring the theory to life with a series of explicit, fully verified examples of 2-representations.
    \item \textbf{Section \ref{rel_12rep}} builds the bridge between the old and the new, establishing the functor from the category of faithful representations into our 2-representation bicategory.
    \item \textbf{Section \ref{str_2grp}} turns to our main character: the string 2-group. After reviewing its various models, for each positive energy representation of the loop group, we associate a 2-representations of the string 2-group.
    \item \textbf{Section \ref{biequiv_2rep}} contains the proof of the crucial biequivalence invariance theorem, a technical centerpiece of the paper.
    \item \textbf{Section \ref{2Vect}} shifts gears into geometry. We construct the 2-stacks of 2-vector bundles and principal 2-bundles, culminating in the associated 2-vector bundle construction.
    \item \textbf{Section \ref{equiv2VB}} introduces symmetries to the geometric picture, developing the theory of 2-equivariant 2-vector bundles over Lie groupoids.
    \item Finally,   \textbf{Section \ref{app_add}} situates our work in the broader literature and charts a course for future investigation, from higher Dixmier--Douady theory to applications in elliptic cohomology. This circle of ideas traces back to the work of Stolz and Teichner \cite{st, st2011}, who proposed a geometric interpretation of elliptic cohomology via representations of loop groups and 2-dimensional field theories. The framework developed in this paper can be seen as a step toward a higher categorical realization of their program for the string 2-group.
\end{itemize}

\subsection{Conventions}  \label{sec:conventions}


Throughout the paper, $k$ denotes a field of characteristic zero (typically $\mathbb{R}$ or $\mathbb{C}$). All categories and bicategories are assumed to be small with respect to a suitable Grothendieck universe unless otherwise stated. Smooth manifolds are assumed to be paracompact and finite-dimensional. Lie 2-groups are assumed to be coherent (weak) unless explicitly stated as strict.
Throughout the paper, we use the  symbols  as listed in Figure \ref{symbol_list}.

\begin{figure} \begin{center} 
\begin{tabular}{|c | c | } 
 \hline 
 \textbf{Symbol}  & \textbf{Meaning} 
\\ \hline
$(B_0, B_1, B_2)$ & a bicategory with objects $B_0$,\\
 & $1$-morphisms $B_1$, and $2$-morphisms $B_2$ \\ \hline 
 $-\horiprod -$  & the horizontal product of two $2$-morphisms  \\
  & in the given bicategory \\
 \hline
 $-\vertprod - $ & the vertical product of two $2$-morphisms \\
  & in the given bicategory  \\ \hline
  $- \horicirc -$ &  the multiplication of two elements  \\
   & in a coherent Lie 2-group \\ \hline
  $-\horicirc -$  & the horizontal composition of $1$-morphisms \\
    &in a bicategory \\ \hline
   $1_{f}$  & the identity $2$-morphism \\
    & at the $1$-morphism $f$ \\
   \hline $\id_x$ & the identity $1$-morphism \\
   & at the object $x$ \\ \hline
   $ \Hom_{\huaC}(A, B)$ & the set of morphisms from object $A$ \\
    & to $B$ in a category $\huaC$\\ \hline
   $\Fun(\huaA, \huaB)$   & the category of functors from $\huaA$ to $\huaB$\\
   \hline
  $ \pFun(-, -) $& the bicategory of pseudofunctors, \\ & strong transformations, and modifications \\ \hline
   $\smap(-, -)$ & simplicial mapping space \\ \hline
   $\Auto(-)$ & the category of automorphisms of an object\\ & in a bicategory \\ \hline
   $\Aut(-)$ & the group of automorphisms of an object \\ & in a category \\ \hline
\end{tabular} \caption{  The Symbol List }  \label{symbol_list}
\end{center} \end{figure}

\section{Lie 2-groups and 2-representations} \label{2grp_basic}

\subsection{Basic concepts}

The concept of a Lie 2-group is well-established in the literature (see, for example, \cite{WZ, schommer:string-finite-dim, baez:2gp}). Several equivalent definitions exist, and here we adopt the definition of a group object within a suitable bicategory.

\begin{definition}\label{defi:2group} 
Let $C$ be a bicategory with finite products. A group object in $C$ (or $C-$group, for brevity) consists of the following data:
\begin{itemize}
    \item an object $G$ in $C$
    \item  a list of 1-morphisms \begin{itemize} \item $m : G \times G\longrightarrow G$ (the multiplication)
\item $u : \ast\longrightarrow G$ (the unit) \end{itemize}
such that
\begin{equation}(pr_1, m) : G \times G \longrightarrow  G \times G 
\end{equation}is an equivalence in $C$.
\item a list of invertible 2-morphisms
\begin{itemize} \item $a : m \circ (m \times 1) \Rightarrow m \circ (1\times m)$ (the associator)
\item $l : m \circ (u \times 1) \Rightarrow 1$ (the left unit constraint) \item
$r : m \circ (1\times u) \Rightarrow 1$ (the right unit constraint)
\end{itemize} where  $1$ is the identity morphism,
subject to the following coherence conditions:  
\begin{itemize}
    \item 
the pentagon identity for associator, 
\item the triangle identity for the left and right unit laws. 
\end{itemize}

\end{itemize} \label{2group}
If $C=\Mfd$ is the category of smooth manifolds and smooth maps, then such a group object in $\Mfd$ is called a Lie 2-group. 
 \end{definition}
In this subsection we will also use $g f$ to denote their multiplication $m(g, f)$ and use $I$ to denote the image of the unit $u$ in $G$.
 
\begin{remark}
The coherence conditions in Definition \ref{2group}, as detailed in \cite[Definition 2.8]{baez:2gp}, require the associator $a$ to satisfy a pentagon identity. Furthermore, there are three diagrams related to the unit laws that should commute, as follows:
\begin{equation} \xymatrix{m \circ (m\times 1)\circ (1\times u \times 1) 
  \ar[rr]^{a\circ (1\times u \times 1) } \ar[dr]_{m\circ (r\times 1)}
& &  m \circ (1\times m) \circ (1\times u \times 1) \ar[dl]^{m\circ (1\times 1) } \\& m} \label{unitdef}
\end{equation}
                
\begin{equation} \xymatrix{m \circ (m\times 1)\circ (1\times 1 \times u)  \ar[rr]^{a\circ (1\times 1\times u) } \ar[dr]_r & & m \circ (1 \times m) \circ (1\times 1 \times u)  \ar[dl]^{m\circ (1\times r) } \\ & m} \label{unitright}
\end{equation}

\begin{equation} \xymatrix{m \circ (m\times 1)\circ (u\times 1 \times 1) \ar[rr]^{a\circ (u\times 1 \times 1)}\ar[dr]_{m\circ (l \times 1)  }
                &  & m \circ (1\times m)\circ (u\times 1 \times 1)  \ar[dl]^{m\circ l} \\
                & m}\label{unitleft} \end{equation}
However, as specified in the definition of a coherent 2-group  \cite[Definition 2.8]{baez:2gp}, the diagram \eqref{unitdef} is sufficient. This is because, according to  the coherence conditions, the maps $r: u \times u\longrightarrow u$ and $l: u\times u\longrightarrow u$ are identical. Then, following a similar approach to the proof in \cite[Exercies 1, VII.1]{maclane:cat-math}, we can demonstrate that the pentagon diagram for the associator $a$ and the diagram \eqref{unitdef} imply that \eqref{unitright} and \eqref{unitleft} are both commutative (as shown in the following lemma).
\end{remark}
                
\begin{lemma}
The pentagon diagram of the associator $a$ and  \eqref{unitdef} implies \eqref{unitright} and \eqref{unitleft}.\label{maclaneex7.1}
\end{lemma}
\begin{proof}
Let $g$ and $f$ be two objects in $G$. 

From the pentagon diagram for $((gf)I)I$ and the diagram \eqref{unitdef}, we have the following commutative diagrams:
\begin{equation}\xymatrix{(gf)I &((gf)I)I\ar[l]^{r_{(gf)I}}\ar[ld]^{a_{gf, I, I}}\ar[r]^{a_{g, f, I}\horiprod 1_I} &(g(fI))I\ar[d]^{a_{g, fI, I}}\\
(gf)(II)\ar[u]^{1_{gf}\horiprod  l_I} \ar[d]^{a_{g, f, II}} &&g((fI)I)\ar[lld]_{1_g\horiprod  a_{f, I, I}}\ar[d]^{1_g\horiprod  r_{fI}}\\ g(f(II))\ar[rr]^{1_g\horiprod (1_f\horiprod  l_I)} &&g(fI) }\label{pentagon:11:2}\end{equation}

Moreover, by naturality, we have
\begin{equation}\xymatrix{(gf)(II) \ar[r]^a \ar[d]^{1_{gf}\horiprod  r_I} &g(f(II)) \ar[d]^{1_g\horiprod (1_f\horiprod  r_I)} \\
(gf)I\ar[r]^a & g(fI)}  \label{natii}\end{equation}

Thus, \eqref{pentagon:11:2} and \eqref{natii} give us the commutative pentagon diagram below.

\begin{equation}\xymatrix{((gf)I)I \ar[rr]^{a_{g, f, I}\horiprod 1_I} \ar[dd]^{r_{(gf)I}} &&
(g(fI))I\ar[rd]^{a_{g, fI, I}} \ar[dd]^{r_{g(fI)}} & \\
&&& g((fI)I)\ar[ld]^{1_g\horiprod r_{fI}} \\
(gf)I\ar[rr]^{a_{g, f, I}} &&g(fI) &} 
\end{equation}
By naturality, the left square diagram commutes. Thus, the triangle diagram on the right commutes.

Finally, we apply the naturality of $r$ we have the commutative diagram

\begin{equation}\xymatrix{(gf)I\ar[rd]^{r_{gf}} \ar[rr]^{a_{g, f, I}}&& g(fI)\ar[ld]^{1_g\horiprod r_f} \\
&gf &&.}  \end{equation}

The commutativity of \eqref{unitright}  can be proved similarly. 
\end{proof}

\subsection{2-Representations} \label{def:2rep_bicat}

Recall that a representation of a group $G$ is a group homomorphism  from $G$ to $ \Aut(V)$ for a vector space $V$. This is a functor from $BG$,  the group $G$ viewed as a category with a single object $\pt$, to the category $\Vect$ of vector spaces. The single object $\pt$ maps to   $V$, and each element $g$ in $G$ maps to an automorphism of $V$. We can define $2$-representation based on a higher version of this idea. 

\begin{definition}
    Let $\huaG_{\bullet}$ be a coherent Lie 2-group. 
A {\em 2-representation} of a 2-group is a pseudofunctor $\rho: B\huaG_{\bullet}\longrightarrow 2\Vect$ 
from $B\huaG_{\bullet}$, the delooping of a coherent 2-group $\huaG_{\bullet}$ viewed as a bicategory with a single object $\pt$,  to the bicategory  $2\Vect$ of $k$-prelinear  categories. The \emph{bicategory of 2-representations} of $\huaG_{\bullet}$ 
with values in $2\Vect$ is then
\begin{equation}
2\Rep(\huaG_{\bullet} ):= \pFun (B\huaG_{\bullet}, 2\Vect), \label{2repdef}
\end{equation} where $\pFun (B\huaG_{\bullet}, 2\Vect)$
denotes the bicategory of pseudofunctors $B\huaG_{\bullet}\to 2\Vect$, strong transformations, and modifications. 
\end{definition}

\begin{remark}
    
In this article we choose  the strict 2-category 2\Vect \hspace{0.1cm} to be the model of 2-vector spaces. The objects of 2\Vect 
 \hspace{0.07cm} are  {\em $k$-prelinear  categories}, 
which are categories enriched in $k$-vector spaces (where $k$ is a field). More precisely, a $k$-prelinear category is a category whose morphism sets are $k$-vector spaces such that compositions are $k$-bilinear. 
The 1-morphisms of $2\Vect$ 
are $k$-linear functors and 2-morphisms are natural transformations between $k$-linear functors. Note that since direct sums of vector spaces are biproducts, a category enriched in $k$-vector spaces is necessarily an additive category. 

\end{remark}

\begin{remark}
\label{rem:pseudofunctor-essential}
The definition of a 2-representation as a pseudofunctor, rather than a strict 2-functor, 
is essential for the theory to capture the correct higher geometric structure. 

Consider the simplest class of Lie 2-groups: those arising from abelian Lie groups. 
 a \emph{strict} $2$-functor $BA \to 2\Vect$, with $A$ an abelian Lie group, 
would indeed yield arbitrary representations of $A$ without any continuity requirement, 
as one would obtain for discrete groups. In contrast, a \emph{pseudofunctor} 
$BA \to 2\Vect$ carries additional coherence data, the lax functoriality constraint $F^2$ and the lax unity constraint $F^0$, which satisfy pentagon and triangle identities. These coherence 
conditions encode the smoothness of the representation at the higher categorical 
level, in the same way that the descent condition for a stack encodes the smooth 
gluing of local data (see Section   \ref{2vbconstr}). 

This illustrates a general principle of the framework: smoothness and continuity 
are not imposed as external analytic conditions, but are built into the categorical 
structure itself via the choice of pseudofunctors over strict functors.
\end{remark}

\begin{remark}
    We work with pseudofunctors rather than general lax functors because the invertibility of their coherence data is essential for our main structural result: biequivalent 2-groups induce biequivalent bicategories of 2-representations, i.e. Theorem \ref{biequiv2rep_grp}. For lax functors, the lax functoriality constraint and the lax unity constraint need not be invertible, which obstructs the construction of the required local equivalences on Hom-categories.
\end{remark}

\begin{example}[Category of groupoid representations]
For a discrete groupoid $\Gamma_{\bullet} = (\Gamma_1 \rightrightarrows \Gamma_0)$, let $\Rep^{\Gamma_{\bullet}}$ denote the category of functors $\Gamma_{\bullet} \to \Vect$ and natural transformations between them. This is the prototype of a 2-vector space proposed in \cite{murray-robert-wockel}\footnote{MRW use the notation $\Vect^{\Gamma_{\bullet}}$, but we reserve this for $\Gamma_{\bullet}$-graded vector spaces when $\Gamma_{\bullet}$ is a group.}. 

If $\Gamma_{\bullet}$ is a Lie groupoid, the same notation $\Rep^{\Gamma_{\bullet}}$ denotes the category of representations of $\Gamma_{\bullet}$, where  objects are vector bundles $V \to \Gamma_0$ equipped with a Lie groupoid morphism $\Gamma_{\bullet} \to \Auto(V)$ \cite{Macdouble2}, and morphisms are $\Gamma_{\bullet}$-equivariant bundle maps. Unlike the discrete case, the presence of differential geometry makes this category non-abelian.
\end{example}

\begin{example}[Morphism in $2\Rep \huaG_{\bullet}$]  \label{nt2repdef}

Let $\huaG_{\bullet}= (\huaG_1 \rightrightarrows \huaG_0)$ denote a coherent 2-group with multiplication $m$, unit $u$, associator $a$, left and right unit constraint $l$ and $r$. 
We explicitly describe the data below for a 1-morphism $T$ in  $2\Rep \huaG_{\bullet}$, based on the definition of a strong transformation between pseudofunctors between bicategories \cite[Section 4.2]{JY:Bicat}.

Recall $\huaG_{\bullet}$-2-representations are pseudofunctors from $B\huaG_{\bullet}$ to $2\Vect$. Let $\huaF_1 = (F_1, F_1^2, F_1^0)$ and $\huaF_2 = (F_2, F^2_2, F^0_2)$ denote two $\huaG_{\bullet}$-2-representations,  where, for $i=1, 2$,
\begin{itemize} 
\item Each $F_i(\ast)$ is an object in $2\Vect$    \hspace{0.05cm}, i.e. a $k$-prelinear category.
\item Each $F_i: \huaG_{\bullet}\rightarrow \Auto(F_i(\ast))$ is a local functor which sends an object  $g_0$ in $\huaG_{\bullet}$ to an automorphism $F_{i, g_0}: F_i(\ast)\to F_i(\ast)$, and a morphism $ ( g_0\xrightarrow{g_1} g_0') \in \huaG_1$ to a natural transformation $F_{i, g_1}: F_{i, g_0}\Rightarrow F_{i, g_0'}$.

\item \[F^2_i: m'_i\circ (F_i \times F_i ) \Rightarrow F_i \circ m\] and \[F^0_i: u'_i \Rightarrow  F_i \circ u \] are  natural isomorphisms,
\begin{equation}
\xymatrix{ \ar @{} [dr] |{\nearrow_{F^2_i}} \huaG_{\bullet} \times \huaG_{\bullet} \ar[d]_{F_i \times F_i  } \ar[r]^{m} & \huaG_{\bullet} \ar[d]^{F_i } \\
\Auto(F_i(\ast)) \times \Auto(F_i(\ast))  \ar[r]^>>>>>{m'_i}& \Auto(F_i(\ast))} \quad \xymatrix{\ar @{}[dr] |{\nearrow_{F^0_i}} \ast\ar[r]^{ u } \ar@{=}[d] & \huaG_{\bullet} \ar[d]^{F_i} \\
\ast\ar[r]_{u'_i} &\Auto(F_i(\ast))} 
\end{equation} which satisfy lax associativity, lax left and right unity \cite[(4.1.3), (4.1.4)]{JY:Bicat}.
\end{itemize}

A strong transformation $T: \huaF_1 \longrightarrow \huaF_2 $ is defined by the following data and axioms. 
\begin{itemize}
    \item Data
\begin{enumerate}
\item A morphism $T(\ast): F_1(\ast)\longrightarrow F_2(\ast)$ in $2\Vect$, i.e. a functor between the $k-$prelinear categories.

\item 2-isomorphisms
\begin{equation}\xymatrix{\ar @{} [dr] |{\nearrow_{T}} \huaG_0\ar[r]^{F_{1}} \ar[d]_{F_2} & \Auto(F_1(\ast))\ar[d]^{T(\ast)_*} \\
    \Auto(F_2(\ast))\ar[r]^{T(\ast)^*} &\Hom_{2\Vect} (F_1(\ast), F_2(\ast))}\label{nt2repdata2}\end{equation}
where  $T(\ast)_*$ denotes the functor
\[\Auto(F_1(\ast)) \rightarrow \Hom_{2\Vect} (F_1(\ast), F_2(\ast)), \quad f\mapsto T(\ast) \circ f \]
and   $T(\ast)^*$ denotes the functor
\[\Auto(F_2(\ast)) \rightarrow \Hom_{2\Vect} (F_1(\ast), F_2(\ast)), \quad f\mapsto f\circ T(\ast) . \]
Thus, for any $g_0\in \huaG_0$,  there is a  $2$-isomorphism \[T_{g_0}: F_{2, g_0}\circ T(\ast) \Rightarrow T(\ast)\circ F_{1, g_0}.\]
\end{enumerate}

\item Axioms

The diagrams below commute.
\begin{enumerate} 
    \item \textbf{Lax Naturality}
    \begin{equation}\xymatrix{(F_{2, g_0}\circ F_{2, h_0} ) \circ T(\ast) \ar@{=}[d] \ar[rrr]^{F^2_{2, (g_0, h_0)} \horiprod 1_{T(\ast)}} 
    &&& F_{2, g_0\horicirc h_0} \circ T(\ast)  \ar[d]^{T_{g_0\horicirc h_0}}\\
    F_{2, g_0}\circ (F_{2, h_0}\circ T(\ast)) \ar[d]_{1_{F_{2, g_0}} \horiprod T_{h_0}} &&& T(\ast)\circ F_{1, g_0\horicirc h_0} \\
    F_{2, g_0}\circ (T(\ast)\circ F_{1, h_0})
    \ar@{=}[d] &&&  T(\ast)\circ (F_{1, g_0}\circ F_{1, h_0}) \ar[u]_{1_{T(\ast)} \horiprod F^2_{1, (g_0, h_0)}} \\
    (F_{2, g_0} \circ T(\ast)) \circ F_{1, h_0} \ar[rrr]_{T_{g_0}\horiprod 1_{F_{1, h_0}}} &&&  (T(\ast) \circ F_{1, g_0})\circ F_{1, h_0}   \ar@{=}[u] \\ }\label{nt2repaxiom1}\end{equation}
    
    \item \textbf{Lax Unity}
    \begin{equation}\xymatrix{I'_2\circ T(\ast) \ar[r]^{l_{T(\ast)}} \ar[d]_{F^0_2\horicirc 1_{T(\ast)}} &T(\ast) \ar[r]^{r^{-1}_{T(\ast)}} 
    & T(\ast)\circ I'_1\ar[d]^{1_{T(\ast)}\horicirc F^0_1} \\
    F_2I\circ T(\ast) \ar[rr]_{T_{I}} &&T(\ast)\circ F_1I}\label{nt2repaxiom2}\end{equation}
\end{enumerate}
\end{itemize}

\end{example}

\begin{remark} \label{rep_cat}
    Since $2\Vect$ is a strict 2-category, the composition of $1$-morphisms in $2\Rep(\huaG_{\bullet})$, i.e. the horizontal product of strong transformations, is strictly associative and strictly unital. Thus, the pseudofunctors $B\huaG_{\bullet}\rightarrow 2\Vect$ and strong transformations between them define a category. When there is no confusion, we will still use the same symbol $2\Rep (\huaG_{\bullet} )$ to denote the category. 
\end{remark}

\begin{lemma}\label{klinear}
For a Lie groupoid $\Gamma_{\bullet}$, $\Rep^{\Gamma_{\bullet}}$ is a  category enriched in $k$-vector spaces, but it is not necessarily an abelian category.
\end{lemma}
\begin{proof}
A morphism $f: \rho \to \rho'$ in $\Rep^{\Gamma_{\bullet}}$ consists of a $\Gamma_{\bullet}$-equivariant, fiberwise linear map. For each $x \in \Gamma_0$, the component $f_x$ lies in $\Hom_{\Vect}(V_x, V'_x)$. The $k$-linear structure on the morphism sets is defined pointwise over $\Gamma_0$, and composition is bilinear: for any $\alpha \in k$,
\[
f \circ (\alpha f') = \alpha (f \circ f'), \qquad (\alpha f) \circ f' = \alpha (f \circ f').
\]
Hence $\Rep^{\Gamma_{\bullet}}$ is enriched in $k$-vector spaces.

However, $\Rep^{\Gamma_{\bullet}}$ is not abelian. If a morphism $f$ does not have constant rank, its kernel may fail to be a vector bundle, and therefore does not define an object of $\Rep^{\Gamma_{\bullet}}$. Consequently, $\Rep^{\Gamma_{\bullet}}$ is not an abelian category.

\end{proof}

\subsection{Examples of 2-representations}\label{ex2rep}
Let $\huaG_{\bullet}= (\huaG_1 \rightrightarrows \huaG_0)$ denote a coherent 2-group with multiplication $m$, unit $u$, associator $a$, left and right unit constraint $l$ and $r$. In Lemma \ref{lem:right:gen-case} we give an example of $\huaG_{\bullet}$-2-representation, $\huaF= (F, F^2, F^0)$, explicitly.
\begin{lemma} \label{lem:right:gen-case}
A 2-group $\huaG_{\bullet}= (\huaG_1 \rightrightarrows \huaG_0)$ can represent on $2\Vect$ in the following way.
\begin{itemize}
    \item Object level: the only object $\ast $ in $B\huaG_{\bullet}$, maps to $\Rep^{\huaG_{\bullet}}$. Here the symbol $\huaG_{\bullet}$ in $\Rep^{\huaG_{\bullet}}$ also denotes the underlying Lie groupoid of the 2-group $\huaG_{\bullet}$ by an abuse of notation.
    \item 1-Morphism level:
    \[ F: 
\huaG_0 \to \Fun(\Rep^{\huaG_{\bullet}}, \Rep^{\huaG_{\bullet}} )_0\footnote{Here $\Fun(-,-)_0$ denotes the set of functors; $\Fun(-,-)_1$ denotes the set of natural transformations.}, \quad g \mapsto F_g,
\]
is defined by: for all $x, y\in \huaG_0$, any $(x\xrightarrow{b}y)\in \huaG_1$, any objects $\rho$, $\rho'$ in $\Rep^{\huaG_{\bullet}}$, and any $1$-morphism $f$ in $\Rep^{\huaG_{\bullet}}$ (a natural transformation $\rho \buildrel{f}\over\Rightarrow \rho'$), 
\[
\begin{split}
    &F_g(\rho): x \mapsto \rho(x\horicirc  g)(=V_{x\horicirc  g}), \quad (x\xrightarrow{b} y) \mapsto (V_{x\horicirc  g} \xrightarrow{\rho(b\horiprod 1_g)} V_{y\horicirc  g})\\
    &F_g(\rho \buildrel{f}\over{\Rightarrow}\rho')=(F_g(\rho) \buildrel{F_g(f)}\over{\Rightarrow} F_g(\rho')),  \text{ 
    a natural transformation with      $F_g(f)_x = f_{x\horicirc  g}$.}
\end{split}
\]
    \item 2-Morphism level:  \[ F: 
\huaG_1 \to \Fun(\Rep^{\huaG_{\bullet}}, \Rep^{\huaG_{\bullet}}  )_1, \quad (g_0 \xrightarrow{g_1} g'_0)  \mapsto F_{g_1},
\] is defined by 
\begin{equation}\label{eq:eta-gamma}
  F_{g_1}(\rho)_x: F_{g_0}(\rho)(x) \to F_{g'_0}(\rho)(x),\quad v\mapsto \rho(1_x \horiprod g_1)(v), \quad \forall v \in F_{g_0}(\rho)(x)
.\end{equation}
\end{itemize}
\end{lemma}
\begin{proof} 
On the level of 1-morphisms, we need to prove that $F_g(\rho)$ is indeed a functor, $F_g(f)$ is indeed a natural transformation, and $F_g$ is $k$-linear, that is, $F_g(f+f')=F_g(f)+F_g(f')$ and $F_g(\lambda f)=\lambda F_g(f)$ for any $
\lambda \in k$. The fact that $F_g(\rho)$ is a functor is obvious, and the $k$-linearity of $F_g$ follows from that of $\Hom$-space of $\Rep^{\huaG_{\bullet}}$.  

The fact that $F_g(f)$ is indeed a natural transformation follows from the fact that $f$ is a natural transformation. Indeed, given $ (x\xrightarrow{b} y )\in \huaG_1$, the diagram
\begin{equation}\label{eq:Fgf}
\xymatrix{F_g(\rho)(x)\ar[r]^{F_g(f)_x}\ar[d]_{F_g(\rho)(b)} & F_g(\rho')(x) \ar[d]^{F_g(\rho')(b)}\\
F_g(\rho)(y)\ar[r]^{F_g(f)_y} &F_g(\rho')(y) } \end{equation} i.e., the diagram
\begin{equation}\xymatrix{\rho(x\horicirc  g)\ar[r]^{f_{x\horicirc  g}}\ar[d]_{\rho(b\horiprod 1_g)} &\rho'(x\horicirc  g)\ar[d]^{\rho'(b\horiprod 1_g)}\\
\rho(y\horicirc  g) \ar[r]^{f_{y\horicirc  g}} & \rho'(y\horicirc g)}\label{eq:Fgf-explicit} \end{equation} commutes because $f$ is a natural transformation.

On the level of 2-morphisms, we need to check that, given $(g_0\xrightarrow{ g_1} g'_0) \in \huaG_1$, what we defined $F_{g_1}(\rho): F_{g_0}(\rho) \Rightarrow F_{g'_0} (\rho)$ is indeed a natural transformation. We need to check the commutativity of the  following diagram for all $( x\xrightarrow{b} y ) \in \huaG_1$,
\begin{equation}\label{eq:eta-g1}\xymatrix{F_{g_0}(\rho)(x)\ar[r]^{F_{g_1} (\rho)_x} \ar[d]_{F_{g_0}(\rho)(b)} &F_{g_0'}(\rho)(x)
\ar[d]^{ F_{g_0'}(\rho)(b)}\\F_{g_0}(\rho)(y)\ar[r]^{F_{g_1} (\rho)_y} &F_{g_0'}(\rho)(y) },
\end{equation}
which is exactly
\begin{equation}
\xymatrix{V_{x\horicirc g_0} \ar[d]_{\rho(b\horiprod 1_{g_0})} \ar[rr]^{\rho(1_x \horiprod g_1)} &  &V_{x\horicirc  g'_0} \ar[d]^{\rho(b\horiprod 1_{g'_0})} \\
V_{y\horicirc g_0} \ar[rr]^{\rho(1_y \horiprod g_1)} & & V_{y\horicirc g'_0}
}.\label{xxxxxxx}
\end{equation} The diagram \eqref{xxxxxxx} follows from the fact that $\rho$ is a functor and the following equation, 
\begin{align*} &(b\horiprod 1_{g'_0})\vertprod (1_x \horiprod g_1) =(b\vertprod 1_x)\horiprod (1_{g'_0}\vertprod g_1)
= b\horiprod g_1   \\   = & (1_y\vertprod b)\horiprod (g_1\horiprod 1_{g_0})   = (1_y \horiprod g_1) \vertprod (b \horiprod 1_{g_0}).\end{align*}
\emptycomment{
Since $\rho$ is a functor, we have 
\begin{align}\rho((b\horiprod 1_{g'_0})\vertprod (1_x \horiprod g_1)) &= \rho ((1_y \horiprod g_1) \vertprod (b \horiprod 1_{g_0})).\\
\rho((b\horiprod 1_{g'_0}))\vertprod\rho( (1_x \horiprod g_1)) &= \rho ((1_y \horiprod g_1) )\vertprod \rho ((b \horiprod 1_{g_0}))\end{align} i.e. the diagram \eqref{xxxxxxx} commutes.}

To show that the data in the above lemma gives us a morphism from $B\huaG_{\bullet}$ to $2\Vect$, we still need to verify the following (see \cite[1.1]{Leinster}).
\begin{enumerate}
\item[1)] \label{itm:1} $F: \huaG_{\bullet} \to \Auto(\Rep^{\huaG_{\bullet}})$ from the category $\huaG_{\bullet}$ (whose structure is the underlying Lie groupoid of the 2-group $\huaG_{\bullet}$) to the category $\Auto(\Rep^{\huaG_{\bullet}})$ of autofunctors from $\Rep^{\huaG_{\bullet}}$ to itself and the natural transformations between them, is a functor. 

\item[2)] \label{itm:2} There exist natural isomorphisms 
\begin{equation} \label{eq:phi}
    F^2: m'\circ (F\times F) \Rightarrow F\circ m, 
\end{equation}
and 
\begin{equation} \label{eq:psi}
    F^0: u' \Rightarrow F\circ u, 
\end{equation}

\begin{equation}
\xymatrix{ \ar @{} [dr] |{\nearrow_{F^2}} \huaG_{\bullet} \times \huaG_{\bullet} \ar[d]_{F\times F} \ar[r]^{m} & \huaG_{\bullet} \ar[d]^{F} \\
\Auto(\Rep^{\huaG_{\bullet}}) \times \Auto(\Rep^{\huaG_{\bullet}})  \ar[r]^>>>>>{m'}& \Auto(\Rep^{\huaG_{\bullet}})}, \quad \xymatrix{\ar @{} [dr] |{\nearrow_{F^0}} \ast\ar[r]^{ u } \ar@{=}[d] & \huaG_{\bullet} \ar[d]^{F} \\
\ast\ar[r]_{u'} &\Auto(\Rep^{\huaG_{\bullet}})}\label{1morcoh}.
\end{equation}

Here $ m'(-, -)$ denotes the composition of functors and horizontal contacnation of natural transformations in $\Auto(\Rep^{\huaG_{\bullet}})$ and $u': \ast\longrightarrow \Auto(\Rep^{\huaG_{\bullet}})$ is the unit in $\Auto(\Rep^{\huaG_{\bullet}})$. Thus we write $m'(F, F')$ also simply as $ F\circ F'$. 

\item[3)]\label{itm:3}
Moreover, $F^2$ and $F^0$  satisfy lax associativity, lax left and right unity, i.e. they make the following diagrams commute. Here we use $I'$ to denote $u'(\ast)$. 
\begin{equation} \label{eq:6-arrow}
    \xymatrix{& F_{g_0\horicirc  h_0} \circ F_{f_0} \ar[rd]^-{F^2_{(g_0\horicirc  h_0,f_0)}} &\\ (F_{g_0} \circ F_{h_0} ) \circ F_{f_0}  \ar[ru]^-{F^2_{(g_0, h_0)} \horiprod 1_{F_{f_0}}} \ar@{=}[d] & & F_{(g_0\horicirc  h_0) \horicirc  f_0} \ar[d]^{F_{a_{g_0, h_0, f_0}}} \\
    F_{g_0} \circ (F_{h_0} \circ F_{f_0}) \ar[rd]_-{1_{F_{g_0}} \horiprod F^2_{(h_0, f_0)}} & & F_{g_0\horicirc  (h_0 \horicirc  f_0)} \\
    & F_{g_0} \circ F_{h_0\horicirc  f_0} \ar[ru]_-{F^2_{(g_0, h_0\horicirc  f_0)}}  &}
\end{equation}

\begin{equation}
\xymatrix{ &F_{I}\circ F_{g_0} \ar[rd]^{F^2_{(I, g_0)}}  &\\ I'\circ F_{g_0}\ar[ru]^{F^0\horiprod 1_{F_{g_0}}} \ar[d]_{\id}  & &F_{I\horicirc  g_0}\ar[d]^{F_l}\\
F_{g_0}\ar@{=}[rr] &&F_{g_0} } \label{morax32}
\end{equation}

\begin{equation}\xymatrix{ & F_{g_0}\circ F_{I} \ar[rd]^{F^2_{(g_0, I)}}  & \\F_{g_0}\circ I' \ar[ru]^{1_{F_{g_0}}\horiprod F^0} \ar[d]_{\id} & &F_{g_0\horicirc  I} \ar[d]^{F_{r}}\\
F_{g_0} \ar@{=}[rr] && F_{g_0}}  \label{morax22}\end{equation}

\end{enumerate}
We now verify the three items one by one. \\
\\
\underline{Verification of item 1):}
To verify 1), we need to see that $F$ preserves identities and compositions. That is, we need to verify that $F_{1_{g_0}}=1_{F_{g_0}}$ and $F_{g'_1\vertprod g_1} = F_{g'_1}\circ F_{g_1}$. These follow from some direct calculation. The morphism $g_0\buildrel{1_{g_0}}\over\longrightarrow g_0$ is sent to $F_{1_{g_0}}$, which is defined by: for any $x, g_0\in \huaG_0$, any composable $g_1, g'_1 \in \huaG_1$, and   any   object  $\rho$   of $\Rep^{\huaG_{\bullet}}$, 
    \[F_{1_{g_0}}(\rho)_x: F_{g_0}(\rho)(x)\longrightarrow F_{g_0}(\rho)(x), \quad v\mapsto \rho(1_x\horiprod 1_{g_0}) (v)=\rho(1_{x\horicirc  g_0})(v)=v;\] and, 
\begin{align*}
F_{g'_1 \vertprod g_1}(\rho)_x &= \rho(1_x \horiprod (g'_1 \vertprod g_1))
= \rho((1_x\vertprod 1_x) \horiprod (g'_1 \vertprod g_1))
\\
& =
\rho((1_x \horiprod g'_1)\vertprod(1_x \horiprod g_1)) = \rho(1_{x} \horiprod g'_1) \circ \rho(1_x \horiprod g_1) \\&= (F_{g'_1} \circ F_ { g_1})(\rho)_x.
\end{align*}
\\
\underline{Verification of item 2):} it is easy to see that
\begin{align}
F_{g_0\horicirc   \hat{g}_0}(\rho)(x)&= \rho(x\horicirc   (g_0\horicirc  \hat{g}_0)) \\
&\xleftarrow{ \rho(a_{x, g_0, \hat{g}_0}) } \rho ((x\horicirc   g_0)\horicirc   \hat{g}_0)  \\
&= F_{g_0} (F_{\hat{g}_0}(\rho)) (x)= (F_{g_0} \circ  F_{\hat{g}_0})(\rho)(x).
 \end{align}

This shows that there is a natural   isomorphism  $F_{g_0\horicirc  \hat{g}_0} \xleftarrow{F^2_{(g_0, \hat{g}_0)}} F_{g_0} \circ F_{\hat{g}_0}$ between these two autofunctors, given by \begin{equation}(F^2_{(g_0, \hat{g}_0)})(\rho) = \rho(a_{-, g_0, \hat{g}_0}): \hspace{0.5cm}  
 (F_{g_0} \circ F_{\hat{g}_0})(\rho) \Rightarrow   F_{g_0\horicirc \hat{g}_0}(\rho)  .  \label{rep.mul}\end{equation} at each object $\rho$ in $\Rep^{\huaG_{\bullet}}$. 
Next we verify that $F^2_{(g_0, \hat{g}_0)}$ defined above is indeed a natural transformation. Let $$\alpha: \rho\Rightarrow \rho'$$ be a 1-morphism in $\Rep^{\huaG_{\bullet}}$. Then we need to show the diagram commutes.
\begin{equation} \xymatrix{F_{g_0\horicirc  \hat{g}_0}(\rho)\ar[d]_{F_{g_0\horicirc  \hat{g}_0}(\alpha)}   &&(F_{g_0} \circ F_{\hat{g}_0})(\rho) \ar[d]^{(F_{g_0} \circ F_{\hat{g}_0})(\alpha)} 
\ar[ll]_{{\rho(a_{-, g_0, \hat{g}_0})}}   \\ F_{g_0\horicirc  \hat{g}_0}(\rho') &&(F_{g_0} \circ F_{\hat{g}_0})(\rho') \ar[ll]_{\rho'(a_{-, g_0, \hat{g}_0})} .}\end{equation}
It suffices to show the diagram below commutes for any $x$ \begin{equation} \xymatrix{\rho(x\horicirc (g_0\horicirc  \hat{g}_0))\ar[d]_{\alpha_{x\horicirc (g_0\horicirc  \hat{g}_0)}}   &&\rho((x\horicirc  g_0)\horicirc  \hat{g}_0) 
\ar[ll]_{{\rho(a_{x, g_0, \hat{g}_0})}}    \ar[d]^{\alpha_{(x\horicirc  g_0)\horicirc  \hat{g}_0}} \\ \rho'(x\horicirc (g_0\horicirc  \hat{g}_0))  &&\rho'((x\horicirc  g_0)\horicirc  \hat{g}_0)
\ar[ll]_{\rho'(a_{x, g_0, \hat{g}_0})}   .}\end{equation}
which commutes because of the naturality of the associator.
Then $$F^2: m'\circ (F\times F) \Rightarrow F\circ m,$$ given by $F^2_{(g_0, \hat{g}_0)}$ at each   $(g_0, \hat{g}_0) \in \huaG_0\times \huaG_0$ is the natural transformation that we desire. The naturality of $F^2$, 
\begin{equation}\label{diag:nat-phi}
    \xymatrix{F_{g_0 \horicirc  \hat{g}_0} \ar[d]_{F_{g_1\horiprod \hat{g}_1}}   && F_{g_0} \circ F_{\hat{g}_0} 
     \ar[ll]_{F^2_{(g_0, \hat{g}_0)}}   \ar[d]^{\rho((1_{-}\horiprod g_1) \horiprod \hat{g}_1)} \\
    F_{g'_0 \horicirc  \hat{g}'_0} &&F_{g'_0} \circ F_{\hat{g}'_0}  \ar[ll]_{F^2_{(g'_0, \hat{g}'_0)}} , }
\end{equation}
follows from the naturality of $a$ and the functoriality of $\rho$. 
\emptycomment{here calculation to be omitted: we can put it into emptycomment when the article is finished otherwise, it could be too long. but for our record now of calculation, let's write still: diagram \eqref{diag:nat-phi} says
\[
(\rho((1_{-}\horiprod g_1) \horiprod \hat{g}_1)\circ F^2_{(g_0, \hat{g}_0)} = F^2_{(g'_0, \hat{g}'_0)} \circ F_{g_1\horiprod \hat{g}_1} ,
\]
which is equivalent to
\[
(\rho((1_{-}\horiprod g_1) \horiprod \hat{g}_1)\circ \phi_{(g_0, \hat{g}_0)}
 = \rho(a^{-1}_{-, g'_0, \hat{g}'_0})\circ \rho(1_{-}\horiprod(g_1\horiprod \hat{g}_1)).
\] This is implied by naturality of $a$ and functoriality of $\rho$. 
}

Now we construct the natural transformation $F^0$ in \eqref{1morcoh} by \begin{equation}
F^0_\ast: I'=u'(\ast)\to F(\ast)=F_I, \quad F^0_\ast(\rho)(x)=\rho(r_x)^{-1}.\label{rep.unit}\end{equation} 
The naturality of $r$ and the functoriality of $\rho$ implies that $F^0_\pt$ is indeed a natural transformation between the two autofunctors $I'$ and $F_I$. Then
the naturality of $F^0$ is trivial to be shown because there is only identity morphism in $\ast$. \emptycomment{ Calculation to be omitted: For a morphism $\rho\xrightarrow{f} \rho'$, to show that $F^0_\pt$ is indeed a natural transformation between the two autofunctors $I'$ and $F_I$, we need to show that 
\[
\xymatrix{I'(\rho) \ar[r]^{\phi_\pt (\rho) } \ar[d]^{I'(f)} & F_I(\rho) \ar[d]_{F_I(f)} \\
I'(\rho') \ar[r]^{\phi_\pt(\rho')} & F_I(\rho')
}
\] Since $I'$ is the identity in $\Aut(\Vect^{\huaG_{\bullet}})$, this is equivalent to show that for each point $x\in \huaG_0$, 
\[
\xymatrix{\rho(x) \ar[r]^{\rho(r^{-1}_x)} \ar[d]^f_x & \rho(x\horicirc  I) \ar[d]^{f_{x\horicirc  I}} \\
\rho'(x)\ar[r]^{\rho'(r_x^{-1})} & \rho'(x\horicirc  I) 
}
\] Since $\rho$ is a functor from the underlying groupoid of $\huaG_{\bullet}$ to $GL(V)$, $\rho(r_x^{-1})=\rho(r_x)^{-1}$. Thus the above diagram is implied by the naturality of $r$.  
}
\\
\noindent \underline{Verification of item 3):}
The 6-arrow diagram \eqref{eq:6-arrow} when evaluated on an object $\rho$ in $\Rep^{\huaG_{\bullet}}$ becomes
\begin{equation} \label{eq:6-arrow-rho}
    \xymatrix{ & F_{g_0\horicirc  h_0} \circ F_{f_0}(\rho)  \ar[rd]^-{\rho(a_{-, g_0\horicirc  h_0,f_0})}  &\\(F_{g_0} \circ F_{h_0} ) \circ F_{f_0} (\rho)  \ar[ru]^-{\rho(a_{-, g_0, h_0}) \horiprod 1_{F_{f_0}}(\rho)} \ar@{=}[d]  & & F_{(g_0\horicirc  h_0) \horicirc  f_0} (\rho) \ar[d]^{F_{a_{g_0, h_0, f_0}}(\rho)} \\
    F_{g_0} \circ (F_{h_0} \circ F_{f_0}) (\rho) \ar[rd]_-{1_{F_{g_0}}(\rho) \horiprod \rho(a_{-, h_0, f_0})} & & F_{g_0\horicirc  (h_0 \horicirc  f_0)}(\rho) \\
     & F_{g_0} \circ F_{h_0\horicirc  f_0} (\rho) \ar[ru]_-{\rho(a_{-, g_0, h_0\horicirc  f_0})}   & }
\end{equation}
Notice that $1_{F_{f_0}}(\rho) = \rho(1_{f_0})$, and when we further evaluate on $x\in \huaG_0$, the natural transformation $\rho(a_{-, g_0, h_0}) \horiprod \rho(1_{f_0})$ at $x$ as a morphism $$\rho(((x \horicirc  g_0)\horicirc  h_0) \horicirc  f_0)\to \rho((x\horicirc  (g_0\horicirc  h_0)) \horicirc  f_0), $$   is exactly $\rho(a_{x, g_0, h_0} \horiprod 1_{f_0})$. 
So the diagram follows from the pentagon condition of $a$ and the functoriality of $\rho$.

In the diagram \eqref{morax32} the composition 
\begin{align*}
    (F_{l}\vertprod F^2_{(I, g_0)} )(\rho)(x)&=\rho (1_{x } \horiprod l_{g_0})\vertprod \rho(a_{x, I, g_0} )\\
&= \rho( (1_{x } \horiprod l_{g_0})\vertprod a_{x, I, g_0})= \rho (r_x\horiprod 1_{g_0})\\ &=\rho(r_x)\horiprod \rho(1_{g_0})
\end{align*}
Thus, $(F_{I}\vertprod F^2_{(I, g_0)} )(\rho)(x)\vertprod (F^0\horiprod 1_{F_{g_0}})$ is identity.

By Lemma \ref{maclaneex7.1}, we also have the commutative diagram \begin{equation}\xymatrix{(x\horicirc g_0)\horicirc  I\ar[r]^{a_{x, g_0, I}} \ar[rd]^{  r_{x\horicirc g_0}} & x \horicirc (g_0\horicirc I ) \ar[d]^{1_x\horiprod r_{g_0}}\\
&x\horicirc g_0}\end{equation}

Thus, 
in the diagram \eqref{morax22} the composition \begin{align*}(F_r\vertprod F^2_{(g_0, I)})(\rho)(x)&=\rho (1_x\horiprod r_{g_0})\vertprod \rho(a_{x, g_0, I} )\\
&= \rho ((1_x\horiprod r_{g_0})\vertprod a_{x, g_0, I} )\\
&= \rho (r_{x\horicirc g_0}).\end{align*}

Thus, the composition $(F_r\vertprod F^2_{(g_0, I)})(\rho)(x)\vertprod ( 1_{F_{g_0}}\horiprod  F^0)(\rho)(x)$ is the identity.

\bigskip

In conclusion, $\huaF:= (F, F^2, F^0)$ defines a $\huaG_{\bullet}$-2-representation. 
\end{proof}

Let $$\Rep_{tr}^{\huaG_{\bullet}}$$ denote the full subcategory of $\Rep^{\huaG_{\bullet}}$ which consists of trivial vector bundles. Each object $\rho$ in $\Rep_{tr}^{\huaG_{\bullet}}$ comes with a trivial bundle of the form $V\times \huaG_0$ with $V$ a fixed vector space.

\begin{corollary}\label{lem:trivial-rep}
For a Lie 2-group $\huaG_{\bullet}$, there is another 2-representation similar to what is given in Lemma \ref{lem:right:gen-case} but restricted to the subcategory $\Rep^{\huaG_{\bullet}}_{tr}$, that is, we replace $\Rep^{\huaG_{\bullet}}$ by $\Rep^{\huaG_{\bullet}}_{tr}$ everywhere.  More precisely, there is a 2-representation of $\huaG_{\bullet}$ with 
\begin{itemize}
    \item Object level: the only object $\pt$ in the bicategory corresponding to a 2-group, maps to $\Rep^{\huaG_{\bullet}}_{tr}$. 
    \item 1-Morphism level: 
    \[
F: \huaG_0 \to \Fun(\Rep^{\huaG_{\bullet}}_{tr} , \Rep^{\huaG_{\bullet}}_{tr})_0, \quad g \mapsto F_g,
\]
 is defined by: for all $x, y\in \huaG_0$, any $(x\xrightarrow{b}y)\in \huaG_1$, any object $\rho$ in $\Rep^{\huaG_{\bullet}}_{tr}$, the only element of which in the image of $\huaG_0$ is denoted by $V$, and for any $1$-morphism $\rho \buildrel{f}\over\Rightarrow \rho'$ in $\Rep^{\huaG_{\bullet}}_{tr}$, 
 \[
\begin{split}
    F_g(\rho): x\mapsto V, 
    \quad (x\xrightarrow{b} y) \mapsto (V \xrightarrow{\rho(b\horiprod 1_g)} V )\\
    F_g(\rho \buildrel{f}\over\Rightarrow \rho')= \bigl( F_g( \rho) \xrightarrow{F_g(f)=f} F_g(\rho') \bigr).
\end{split}
\]
    \item 2-Morphism level:    \[
F: \huaG_1 \to \Fun(\Rep^{\huaG_{\bullet}}_{tr} , \Rep^{\huaG_{\bullet}}_{tr})_1, \quad (g_0 \xrightarrow{g_1} g'_0) \mapsto F_{g_1},
\] where
\begin{equation}\label{eq:eta-gamma_3}
  F_{g_1}(\rho)_x: F_{g_0}(\rho)(x) \to F_{g'_0}(\rho)(x),\quad v\mapsto \rho(1_x \horiprod g_1)(v), \quad \forall v \in F_{g_0} (\rho) (x).
\end{equation}
\end{itemize}
\end{corollary}
\begin{proof}The 2-representation we describe above is what we have in Lemma \ref{lem:right:gen-case} but restricting to $\Rep_{tr}^{\huaG_{\bullet}}$. Because the inclusion $\Rep_{tr}^{\huaG_{\bullet}} \subset \Rep^{\huaG_{\bullet}}$ is fully faithful, this gives us another 2-representation of $\huaG_{\bullet}$.   
\end{proof}

\subsection{Relation between 1-representations and 2-representations} \label{rel_12rep}

Suppose there is a Lie 2-group $\huaG_{\bullet}=(\huaG_1 \rightrightarrows \huaG_0)$ whose underlying Lie groupoid is an action groupoid, that is, $\huaG_0=G$, $\huaG_1=H\ltimes G$ with the structure morphisms defined by  
\begin{align*} & \source(h, g): = g,  \quad   \tg(h, g) : = h\cdot g, \quad u(g) := (e, g), \quad 
 (h_1, g_1)\vertprod  (h_2, g_2) : = (h_1h_2, g_2), 
\end{align*} where $H$ acts on $G$ from the left with $\cdot$.    
Then an object $\rho$ in $\Rep^{\huaG_{\bullet}}_{tr}$ reads for us a trivial vector bundle $V\times G$ and each $\rho(h, g)\in \End(V)$. Therefore, it is not hard to see that an $H$-representation $\rho_0$ gives us an element $\rho \in \Rep^{\huaG_{\bullet}}_{tr}$, with each
\[
\rho(h, g):= \rho_0(h). 
\] In fact, we may further see that
$\Rep(H) \xrightarrow{i} \Rep^{\huaG_{\bullet}}_{tr}$ is a full subcategory.

Let $(V, H\xrightarrow{\zeta} \Aut(V))$ be a representation of a Lie group $H$. Let $\{V\}/\zeta(H)$ denote the groupoid with one object $V$ and $\Hom(V, V) =\zeta(H)$. Define a $k$-prelinear category \[C_{V}:=\Fun(\huaG_{\bullet}, \{V\}/\zeta(H)),\]  where $\huaG_{\bullet}$ on the right denotes the underlying Lie groupoid of the action groupoid $ \bigl( H\ltimes G \rightrightarrows G \bigr)$. 
There is a fully faithful embedding $C_{V} \hookrightarrow \Rep^{\huaG_{\bullet}}_{tr} $, sending a functor $\rho$ to the trivial bundle $G\times V$ with the action $\rho(h, g) = \zeta(h)$.  

Moreover, we have another corollary of Lemma \ref{lem:right:gen-case}, whose proof is analogous to that of Corollary \ref{lem:trivial-rep}. 

\begin{corollary}\label{lem:CV}
Let $\huaG_{\bullet}$ be a Lie 2-group with underlying Lie groupoid given by an action groupoid $H\ltimes G \Rightarrow G$.  There is a $\huaG_{\bullet}$-2-representation  similar to the construction of $\huaF$ in Lemma \ref{lem:right:gen-case} restricted to the subcategory $C_V$, that is, we replace $\Rep^{\huaG_{\bullet}}$ in the lemma by $C_V$ everywhere. This $\huaG_{\bullet}$-2-representation is denoted  by $$\huaF_V = ( f_V, f^2_V, f^0_V).$$ 
\end{corollary}

In the rest part of the subsection, we work with the category $\Rep^{\faith}(H)$ of \emph{faithful} representations of a Lie group $H$. The objects in the category are pairs $(V, \zeta)$ where $\zeta: H \to \Aut(V)$ is an injective group homomorphism, and the morphisms are $H$-intertwiners. This is a full subcategory of the usual representation category $\Rep(H)$, and the restriction to faithful representations ensures that distinct group elements act by distinct automorphisms, which is imposed to ensure injectivity in Proposition \ref{fi_i_prop}.
The construction of the pseudofunctor $i$ itself in Theorem \ref{rep_2rep_f} does not require it.


\begin{theorem} \label{rep_2rep_f}

Let $\huaG_{\bullet}$ be a Lie 2-group whose underlying Lie groupoid is an action groupoid $\bigl(  H\ltimes G \rightrightarrows G \bigr)$. We view $\Rep^{\faith} H$ as a discrete bicategory. There is a pseudofunctor \[i: \Rep^{\faith} H\longrightarrow 2\Rep (\huaG_{\bullet} ) \]defined as follows.

\begin{itemize}
    \item On the level of objects: for $(V, \zeta) \in \Rep^{\faith} H$, $i(V, \zeta):= \huaF_V$. 
    \item On the level of $1$-morphisms: 
    for each morphism \[(V_0, \zeta_{0})\buildrel{\omega}\over\rightarrow (V_1,\zeta_{1})\] in $\Rep^{\faith} H$, we define a morphism
    $i(\omega)$ \[i( V_0, \zeta_{0})\xrightarrow{(T(\ast), T_{\omega})} i( V_1, \zeta_{1})\] in $2\Rep (\huaG_{\bullet} )$ as below.
   For each $h\in H$,  we have the commutative diagram 
\[\xymatrix{(V_0, \zeta_{0})\ar[r]^{\omega}\ar[d]_{\zeta_{0}(h)} &(V_1, \zeta_{1})\ar[d]^{\zeta_{1}(h)} \\ (V_0, \zeta_{0})\ar[r]^{\omega} &(V_1, \zeta_{1}) .}\] 
For each $\omega$, we can define  a functor 
\[
F_\omega: \{V_0\}/\zeta_0(H)  \to \{V_1\}/\zeta_1(H)\] by sending the object 
 $V_0$ to $ V_1$ and sending each morphism $\zeta_0(h)$ to $\zeta_1(h)$. 

    \begin{enumerate}
\item  $T(\ast)$ is defined to be $F_{\omega *}: C_{V_0}\rightarrow C_{V_1}$. Explicitly, 
an object  $\big(\huaG_{\bullet}\xrightarrow{\rho} \{V_0\}/\zeta_0( H)\big)$ in $C_{V_0}$ is mapped to $$\big(\huaG_{\bullet}\xrightarrow{\rho} \{V_0\}/\zeta_0(H) \xrightarrow{F_\omega} \{V_1\}/\zeta_1(H)\big) .$$ 
 
\item The $2$-isomorphism $T_{\omega}$
\[\xymatrix{\ar @{} [dr] |{\nearrow_{T_{\omega}}} G \ar[r]^{f_{V_0}} \ar[d]_{f_{V_1}} & \Auto(C_{V_0})\ar[d]^{T(\ast)_*} \\
    \Auto(C_{V_1})\ar[r]^{T(\ast)^*} &\Fun (C_{V_0}, C_{V_1})}\]
 is defined to be identity.
\end{enumerate}

\item  On the level of $2$-morphisms: each identity $2$-morphism in $\Rep^{\faith} H$ is mapped to identity modifications in $2\Rep ( \huaG_{\bullet}  )$. 

\end{itemize}

\label{thm:repHinc2repG}
\end{theorem}

\begin{proof}
This is a direct specialization of the construction in Lemma \ref{ex2rep} and Corollary \ref{lem:trivial-rep} to the full subcategory $C_{V}$
 of $\Rep^{\huaG_{\bullet}}_{tr}$. For each intertwiner $\omega: (V_0, \zeta_0) \to (V_1, \zeta_1)$, the data $(T(*), T_\omega)$ defines a strong transformation between the 2-representations $\huaF_{V_0}$ and $\huaF_{V_1}$. Since $T_\omega$ is taken to be the identity, For each intertwiner $\omega: (V_0, \zeta_0) \to (V_1, \zeta_1)$, the data $(T(*), T_\omega)$ defines a strong transformation between the 2-representations $\huaF_{V_0}$ and $\huaF_{V_1}$. Since $T_\omega$ is taken to be the identity, the lax naturality and lax unity conditions \eqref{nt2repaxiom1} and \eqref{nt2repaxiom2} reduce to the functoriality of $F_{\omega *}$, which follows immediately from the construction. The assignment is strictly functorial and respects identity 2-morphisms, so $i$ is a strict functor, hence a pseudofunctor, from the discrete bicategory $\Rep^{\faith} H $ to $ 2\Rep (\huaG_{\bullet} ) $.
\end{proof}

As discussed in Remark \ref{rep_cat}, $2\Rep (\huaG_{\bullet} )$ has the structure of a category. 
When restricted to the 1-categorical level, the pseudofunctor $i$ constructed in Theorem \ref{rep_2rep_f} yields a functor between the underlying 1-categories.
\begin{proposition} \label{fi_i_prop}
    The  functor $i: \Rep^{\faith} H\longrightarrow 2\Rep (\huaG_{\bullet} )$ is injective on isomorphism classes of objects. 
\end{proposition}

\begin{proof}

If $(V_1, \zeta_1) $ and $ (V_2, \zeta_2)$ are non-isomorphic $H$-representations, then the corresponding groupoids $C_{V_1}$ and $C_{V_2}$ are non-isomorphic because $C_V$ determines the group $\zeta(H)$ up to isomorphism, and $\zeta$ is faithful. Hence, the 2-representations $i(V_1, \zeta_1)$ and $i(V_2, \zeta_2)$ are distinct.

\end{proof}

\begin{lemma} \label{reptr_klinear}
Let $\huaG_{\bullet} = (\huaG_1\rightrightarrows \huaG_0)$ be a coherent Lie 2‑group whose underlying Lie groupoid is an action groupoid $\bigl(  H\ltimes G \rightrightarrows G \bigr)$.
Then $\Rep_{tr}^{\huaG_{\bullet}}$ is a cocomplete $k$-linear category, hence a 2‑vector space.

\end{lemma}
\begin{proof}
As a subcategory of $\Rep^{\huaG_{\bullet}}$, $\Rep_{tr}^{\huaG_{\bullet}}$ is a category enriched in $k$-vector spaces.
Notice that direct sum of vector spaces is a biproduct, thus, a category enriched in $k$-vector spaces is an additive category. 
Now we show that $\Rep_{tr}^{\huaG_{\bullet}}$ is abelian. We define biproduct $\oplus$ in $\Rep_{tr}^{\huaG_{\bullet}}$ in the following way:  for any $ \rho_1, \rho_2 \in \Rep_{tr}^{\huaG_{\bullet}}, $  for any $  g_0 $  in   $ \huaG_0$  and any $g_0\xrightarrow{g_1} g_0' $ in $\huaG_1$, 
\[\rho_1\oplus\rho_2:  \huaG_{\bullet}\longrightarrow \Aut(V_1\oplus V_2), \quad g_0\mapsto \rho_1(g_0)\oplus \rho_2(g_0), \quad  g_1\mapsto \rho_1(g_1)\oplus \rho_2(g_1). \]
where $V_i$ is the trivial vector bundle over $\huaG_0$ on which $\huaG_{\bullet}$ represents by $\rho_i$. 
Let $\underline{0}$ denote the 0-dimensional vector bundle over $\huaG_0$. The zero object in $\Rep^{\huaG_{\bullet}}_{tr}$ is the morphism sending $\huaG_{\bullet}$ to $\Aut(\underline{0})$.
Let $\rho_i$ be an object in $\Rep^{\huaG_{\bullet}}_{tr}$ whose image is a vector space $V_i$, $i=1, 2$. Let $f: \rho_1\longrightarrow \rho_2$ denote a morphism in $\Rep^{\huaG_{\bullet}}_{tr}$. Then it is an $H$-equivariant linear map from $V_1$ to $V_2$.  (And given any $H$-equivariant linear map from $V_1$ to $V_2$, we can define a functor from $\rho_1$ to $\rho_2$.)
The kernel of $f$ is the functor $\ker f$ from $\huaG_{\bullet}$ to the vector space $\ker f$, which is an $H-$invariant subspace of $V_1$. We have the evident  monomorphism from $\ker f$ to $\rho_1$. 
And the cokernel of $f$ is the functor from $\huaG_{\bullet}$ to the vector space $\coker f$, which is an $H-$invariant subspace of $V_2$. We have the evident epimorphism from $\coker f$ to $\rho_2$. 
Since the category  of $H-$representations is an abelian category,  for any $H-$equivariant linear map $h: V_1\longrightarrow V_2$, when it is a monomorphism, it is a kernel, and when it is a epimorphism, it is a cokernel. 
Thus, $\Rep^{\huaG_{\bullet}}_{tr}$ is also an abelian category.
The cocompleteness of $\Rep^{\huaG_{\bullet}}_{tr}$ comes from the cocompleteness of the category of $H-$representations.
\end{proof}

\subsection{2-Representations of crossed modules} \label{ex_2rep_cr}

The 2-representations constructed in Section  \ref{ex2rep} apply to general coherent Lie 2-groups. When the 2-group is \emph{strict}---that is, when it arises from a crossed module of Lie groups---there is a  simpler construction that takes advantage of the strictness. We present this construction here for completeness, following the formulation in \cite{murray-robert-wockel}, which uses a right action and is adapted to the string 2-group model.

\begin{definition}A {\em crossed
module of Lie groups} $(H\xrightarrow{\partial} G, \alpha)$ consists of a Lie group morphism $\partial: H\to G$ and a right action
$\alpha$ of $G$ on $H$, namely a Lie group morphism $\alpha: G^{op}\to \Aut(H)$, such that
\begin{itemize}
\item[(Eq)] $\partial(\alpha(g) h) = g^{-1} \partial(h) g$ (Equivariance)
\item[(Pf)] $\alpha(\partial(h))  h' = h^{-1} h' h$ (Pfeiffer identity).
\end{itemize}
\end{definition}
Such a crossed module of Lie groups gives us a strict Lie 2-group whose underlying Lie groupoid is  $H \ltimes G \rightrightarrows G$, where $H$ acts from left on $G$ through
$\partial$. More precisely, we have
\[
\tg(g_1, h)=g_2=\partial(h) g_1 , \quad \s (g_1, h) = g_1, \quad
\xymatrix@1@C+2em{
  \bullet &
  \ar@/_1pc/[l]_{g_2}_{}="0"
  \ar@/^1pc/[l]^{g_1}^{}="1"
  \ar@{<=}"0";"1"^{h}
  \bullet
}.
\]
The groupoid multiplication is visualized as vertical product,
\[
(g_2, h_2) \vertprod (g_1, h_1)  =(g_1, h_2 h_1),
\quad \xymatrix@1@C+1.5em{
  \bullet &
  \ar@/_1.5pc/[l]_{g_3}_{}="0"
  \ar[l]|{\phantom{|}g_2\phantom{|}}^{}="1"_{}="1b"
  \ar@/^1.5pc/[l]^{g_1}^{}="2"
  \ar@{<=} "0";"1"^{h_2}
  \ar@{<=} "1b";"2"^{h_1}
  \bullet},
\]
where $g_3=\partial(h_2) g_2 = \partial(h_2) \partial(h_1)g_1  = \partial(h_2 h_1)g_1 $.

The group structure gives further a horizontal product,
\[(g_1, h) \horiprod (g'_1, h')= (g_1g'_1, h\horiprod h'), \quad
\xymatrix@1@C+2em{
  \bullet&
  \ar@/_1pc/[l]_{g_2}_{}="0"
  \ar@/^1pc/[l]^{g_1}^{}="1"
  \ar@{<=}"0";"1"^{h}
  \bullet &
  \ar@/_1pc/[l]_{g'_2}_{}="2"
  \ar@/^1pc/[l]^{g'_1}^{}="3"
  \ar@{<=}"2";"3"^{h'}
  \bullet}
\;\mapsto\;
\xymatrix@1@C+2em{
  \bullet & \ar@/_1pc/[l]_{g_2g'_2}_{}="4"
  \ar@/^1pc/[l]^{g_1 g'_1}^{}="5"
  \ar@{<=}"4";"5"^{h\horiprod h'}
  \bullet,
}
\]
where $h \horiprod h' = h \alpha(g^{-1}_1) (h') $.  Indeed,
\[
g_2 g'_2 = \partial(h) g_1 \partial(h') g'_1  = \partial(h) g_1 \partial(h') g^{-1}_1 g_1 g'_1
= \partial(h \alpha'(g_1^{-1}) (h') )g_1 g'_1 .
\] Thus, $\bigl( H\ltimes G\rightrightarrows G\bigr) $ is a group object in the (strict) 2-category of Lie groupoid, Lie groupoid morphisms, and 2-morphisms. All the structure  2-morphisms of it 
are identity. In other words, it is a strict Lie 2-group.

From now on we do not distinguish strict Lie 2-groups and crossed modules. 
\begin{lemma} \label{lem:right:strict-case}
A crossed module $(H\xrightarrow{\partial} G, G^{op}\xrightarrow{\alpha} Aut(H))$ represents on $2\Vect$ in the following way
\begin{itemize}
    \item Object level: the only object $\pt$ in the 2-category corresponding to a crossed module, which is a 2-group, maps to the category of $H$-representations, that is $\pt\mapsto   \Rep H$.
    \item 1-Morphism level:
    \[
G\to \Fun(\Rep H , \Rep H)_0, \quad g \mapsto F_g,
\]
where
\[
F_g\big( (V, \rho) \big) = (V, \rho \circ \alpha(g) ), \quad F_g\big((V, \rho) \xrightarrow{f} (V', \rho')\big)=(V, \rho \circ \alpha(g))
\xrightarrow{f} (V', \rho'\circ \alpha(g) ).
\]
for any  $H$-representation $( V, H\xrightarrow{\rho} \Aut(V))$.

    \item 2-Morphism level: \[
G\times H \to \Fun(\Rep H , \Rep H)_1, \quad (g_1, h) \mapsto F_{(g_1,h)},
\]
where
\[
F_{(g_1, h)}(V, \rho): \big(V, \rho\circ \alpha(g_1) \big) \to \big(V, \rho \circ \alpha(g_2) \big), \quad v \mapsto \rho(\alpha(g_1)(h^{-1})) (v).
\]
\end{itemize}
\end{lemma}
\emptycomment{For the definition of $F$, I also tried $ v \mapsto \rho(h^{-1}) (v)$ and $v \mapsto \rho(h) (v)$ and $ v \mapsto \rho(\alpha(g_1)(h)) (v)$,
 All are not $H-$equivariant.}
\begin{proof}
We first check whether $F_g$ is indeed a functor. For this, we need to see whether $f: (V, \rho \circ \alpha(g) ) \to (V', \rho' \circ \alpha(g)
)$ is still $H$-equivariant, thus, still a morphism between $H$-representations. This follows from a direct calculation,
\[
f(h\cdot v) = f(\rho(\alpha(g)(h))v) = \rho'(\alpha(g)(h)) f(v) = h \cdot_{'} f(v).
\]

We now check whether $F_{(g_1, h)}(V, \rho)$ is equivariant: \begin{align*} &
F_{(g_1, h)}(V, \rho)(\hat{h} \cdot v) = \rho(\alpha(g_1)(h^{-1}))( \hat{h} \cdot v) = \rho(\alpha(g_1)(h^{-1}))(\rho( \alpha(g_1) (\hat{h}))(v))
\\= &\rho(\alpha(g_1)(h^{-1}) \alpha(g_1)(\hat{h}) )(v),
\end{align*}
and
\begin{align*}   &\hat{h}\cdot_{'}
F_{(g_1, h)}(V, \rho)(v) =\rho(\alpha(g_2) \hat{h}) (\rho(\alpha(g_1)(h^{-1}))v) \\ = & \rho(\alpha(\partial(h)g_1) (\hat{h})) (\rho(\alpha(g_1)(h^{-1}))v)
= \rho(\alpha(g_1)\alpha(\partial(h)) (\hat{h})\alpha(g_1)(h^{-1}))(v)\\
= & \rho(\alpha(g_1)(h^{-1}\hat{h} h) \alpha(g_1)(h^{-1})) (v) = \rho(\alpha(g_1)(h^{-1}\hat{h} )) (v)
\\  = & \rho(\alpha(g_1)(h^{-1}) \alpha(g_1)(\hat{h}) )(v) = F(\hat{h} \cdot v)\end{align*}
where $g_2 =\partial (h) g_1$. Thus, $F$ is equivariant and, thus, it defines a morphism between
representations.

The following diagram 
\begin{equation}\label{diag:etar}
    \xymatrix{
F_{g_1}\big( (V, \rho)\big) \ar[rr]^{F_{(g_1, h)}(V, \rho)} \ar[d]_{F_{g_1}(f)}& & F_{g_2}\big( (V, \rho) \big) \ar[d]^{F_{g_2}(f)} \\
F_{g_1}\big( (V', \rho') \big) \ar[rr]^{F_{(g_1, h)} (V', \rho')} & &F_{g_2}\big( (V', \rho') \big)
}
\end{equation} commutes. In fact, this boils down to the fact that $f(\rho(\alpha(g_1)(h^{-1}))v) = \rho'(\alpha(g_1)(h^{-1})) f(v)$, which is implied by the fact that $f$ is a morphism of
representations. Thus, $F$, that we construct, is a natural transformation.

Now to verify that these data gives us a morphism from $B\huaG_{\bullet}$ to $2\Vect$, where $\huaG_{\bullet}$ is the strict Lie 2-group corresponds to our crossed module, we still need to verify three items as in Lemma \ref{lem:right:gen-case}. But the strictness simplifies a lot our proof. Now we check these items one by one:\\
\underline{Verification of 1)}: We need to check whether 
\begin{equation} F_{(g_2, h_2)}\circ     F_{(g_1, h_1)} =F_{(g_2, h_2) \vertprod (g_1, h_1)} = F_{(g_1, h_2h_1)}.\label{etaf}\end{equation}
 Consider the composition
\[ F_{g_1}(V, \rho) \xrightarrow{F_{(g_1, h_1)}(V, \rho)}  F_{g_2}(V, \rho) \xrightarrow{F_{(g_2, h_2)} (V, \rho)} F_{g_3}(V, \rho) ,\]
\begin{align}&F_{(g_2, h_2)} (V, \rho)\circ F_{(g_1, h_1)} (V, \rho) (v)\\ =&F_{(g_2, h_2)} (V, \rho)(\rho(\alpha(g_1)(h_1^{-1}))(v))
= \rho(\alpha(g_2)(h_2^{-1}))\rho(\alpha(g_1)(h_1^{-1}))(v)\\
= &\rho(\alpha (\partial (h_1) g_1)(h_2^{-1}) \alpha(g_1)(h_1^{-1}) )(v) \\
= & \rho(\alpha(g_1)\alpha(\partial(h_1))(h_2^{-1}) \alpha(g_1)(h_1^{-1}) )(v)
\\ =& \rho (\alpha(g_1) (h_1^{-1} h_2^{-1} h_1 h_1^{-1}))(v) = \rho(\alpha(g_1)((h_2h_1)^{-1}) ) (v)
\\ =&  F_{(g_1, h_2h_1)}(V, \rho) (v) =F_{(g_2, h_2) \vertprod (g_1, h_1)} (V, \rho) (v) \label{etafre}
\end{align}
\underline{Verification of 2) and 3)}: We see that $F_{g \hat{g}} = F_g \circ F_{\hat{g}}$,
\[
F_{g\hat{g}}(V, \rho)= (V, \rho \circ \alpha(g\hat{g}) ) = (V, \rho  \circ \alpha(\hat{g}) \circ \alpha(g)) =F_{g}(V, \rho
\circ\alpha(\hat{g}))= F_g F_{\hat{g}}(V, \rho).
\] Therefore, the natural isomorphism $F^2$ in \eqref{eq:phi} is the identity. It is then similar to show that $F^0$ in \eqref{eq:psi} is also the identity. Thus the coherence conditions in 3) automatically hold. 
\end{proof}

In addition, applying Lemma \ref{lem:right:gen-case} to the crossed module $(H \xrightarrow{\partial} G, \alpha)$, we obtain another 2-representation of the corresponding strict 2-group.

\begin{corollary} \label{crm_2rep_sqg}
Let $\huaG_{\bullet}$ denote the strict 2-group associated to the crossed module $(\partial, \alpha)$.  We have the two $\huaG_{\bullet}$-representations below.
  \begin{enumerate}
      \item We have the  $\huaG_{\bullet}$-representation $\huaF= (F, F^2, F^0)$ defined below.

\begin{itemize}
    \item Object level: the only object $\ast $ in $B\huaG_{\bullet}$, maps to $\Rep^{\huaG_{\bullet}}$. 
    \item 1-Morphism level:
    \[ F: 
G \to \Fun(\Rep^{\huaG_{\bullet}}, \Rep^{\huaG_{\bullet}} )_0, \quad g \mapsto F_g,
\]
is defined by: for all $x, y\in G$, any $(x\xrightarrow{(x, h)}y)\in H\ltimes G$, any objects $\rho$, $\rho'$ in $\Rep^{\huaG_{\bullet}}$, and any $1$-morphism $\rho \buildrel{f}\over\Rightarrow \rho'$   in $\Rep^{\huaG_{\bullet}}$ , 
\[
\begin{split}
    &F_g(\rho): x \mapsto \rho(x g), \quad (x\xrightarrow{(x, h)} y) \mapsto (V_{x  g} \xrightarrow{\rho(xg, h)} V_{y  g})\\
    &F_g(\rho \buildrel{f}\over{\Rightarrow}\rho')=(F_g(\rho) \buildrel{F_g(f)}\over{\Rightarrow} F_g(\rho')),  \text{  with      $F_g(f)_x = f_{xg}$.}
\end{split}
\]
    \item 2-Morphism level:  \[ F: 
H\ltimes G\to \Fun(\Rep^{\huaG_{\bullet}}, \Rep^{\huaG_{\bullet}}  )_1, \quad (g, h)  \mapsto F_{(g, h)},
\] is defined by 
\begin{equation}\label{eq:eta-gamma2}
  F_{(g, h)}(\rho)_x: \rho(x g)\to \rho(x (\partial{h})g), \quad v\mapsto \rho( xg, \alpha(x)^{-1}(h))(v), \quad \forall v \in F_{g}(\rho)(x).
\end{equation}
\item the structure $2$-morphism $F^2$ and $F^0$ are both trivial. 
\end{itemize}

\item  We have the  $\huaG_{\bullet}$-representation $\huaF_{tr}= (F_{tr}, F^2_{tr}, F^0_{tr})$ defined below.

\begin{itemize}
    \item Object level: the only object $\ast $ in $B\huaG_{\bullet}$, maps to $\Rep^{\huaG_{\bullet}}_{tr}$. 
    \item 1-Morphism level:
    \[ F_{tr}: 
G \to \Fun(\Rep^{\huaG_{\bullet}}_{tr}, \Rep^{\huaG_{\bullet}}_{tr} )_0, \quad g \mapsto F_{tr, g},
\]
is defined by: for all $x, y\in G$, any $(x\xrightarrow{(x, h)}y)\in H\ltimes G$, any object $\rho$ in $\Rep^{\huaG_{\bullet}}_{tr}$, the only element of which in the image of $\huaG_0$ is denoted by $V$, and any $1$-morphism $\rho \buildrel{f}\over\Rightarrow \rho'$   in $\Rep^{\huaG_{\bullet}}$ , 
\[
\begin{split}
    &F_{tr, g}(\rho): x \mapsto V , \quad (x\xrightarrow{(x, h)} y) \mapsto (V \xrightarrow{\rho(xg, h)} V)\\
    &F_{tr, g}(\rho \buildrel{f}\over{\Rightarrow}\rho')=(F_{tr, g}(\rho) \buildrel{F_g(f)}\over{\Rightarrow} F_{tr, g}(\rho')),  \text{  with      $F_{tr, g}(f)_x = f_{xg}$.}
\end{split}
\]
    \item 2-Morphism level:  \[ F: 
H\ltimes G\to \Fun(\Rep^{\huaG_{\bullet}}_{tr}, \Rep^{\huaG_{\bullet}}_{tr}  )_1, \quad (g, h)  \mapsto F_{tr, (g, h)},
\] is defined by 
\begin{equation}\label{eq:eta-gamma3}
  F_{tr, (g, h)}(\rho)_x: \rho(x g)\to \rho(x (\partial{h})g) ,\quad v\mapsto \rho( xg, \alpha(x)^{-1}(h))(v), \quad \forall v \in F_{g_0}(\rho)(x).
\end{equation}

\item the structure $2$-morphism $F^2$ and $F^0$ are both trivial. 
\end{itemize}
  \end{enumerate}
\end{corollary}

\section{Loop group and string 2-group} \label{str_2grp}

In this section we apply the general 2-representation theory developed in Section \ref{2grp_basic} to the string 2-group. The string 2-group $\String(G)$ is the universal higher extension of a compact simple Lie group $G$ that kills the generator of $\pi_3(G)$; it plays a central role in higher gauge theory, string geometry, and loop group representations.

We begin in Section \ref{MRWmodel} by recalling the two models of the string 2-group that will be used: the coherent model $S(QG)$ of Murray–Roberts–Wockel, which is built from the free loop group and is the natural setting for loop group representations, and its strict model  $S(Q_{\flat, *}G)$, which is biequivalent to it. The biequivalence invariance theorem (Theorem \ref{biequiv2rep_grp} below) guarantees that our 2-representation theory does not depend on which model we choose.

In Section \ref{biequiv_2rep} we prove the biequivalence invariance theorem for 2-representations: biequivalent 2-groups have biequivalent 2-representation bicategories. This is the technical core of the section and justifies working with the coherent model.

In Section \ref{ex_str2grp_rep}, we apply the construction of Section \ref{rel_12rep} to the coherent model $S(QG)$: for each positive energy representation of the loop group $\widehat{LG}$, we construct an associated 2-representation of $\String(G)$. This gives a categorical shadow of the loop group representation theory.

\subsection{Models of string Lie 2-group}\label{MRWmodel}

Let  $G$  be a compact simply-connected simple Lie group. The free loop group  $LG$  consists of smooth maps  $\gamma: \R \to G $ with  $\gamma(t+1) = \gamma(t)$, $\forall t\in \R$; it is a Fr\'{e}chet Lie group under pointwise multiplication. The set of quasi-periodic paths
\[ QG:= \{ \gamma\in C^{\infty}(\R, G) \mid \gamma(t+1) \cdot \gamma(t)^{-1} \text{   is constant } \} \]
contains  $LG$  as a subgroup. Let
\begin{equation} \xymatrix{
 1\ar[r] &  S^1 \ar[r] &  \widehat{LG} \ar[r]^\pi & LG \ar[r] & 1 .} \label{lg_cext} \end{equation}
be the standard central extension. The action groupoid
\[ S(QG):= \bigl( (\widehat{LG} \ltimes QG )\rightrightarrows QG\bigr),\] 
is   not a strict Lie 2-group \cite[Section 2]{murray-robert-wockel}.

To obtain a strict model, consider the subset
 $$ Q_{\flat, *} G :=\{ \gamma \in  QG \mid \gamma(0)=e, \; \gamma(t+1)\gamma(t)^{-1}  \text{ is constant  } \}. $$ 

Define  \[ \Omega_{\flat}G :=  LG \cap  Q_{\flat, *} G,  \] and let  $ \widehat{\Omega_{\flat} G} $  be the restriction of the central extension to  $\Omega_{\flat}G $. The action groupoid
\[  S(Q_{\flat, *}G): =  \bigl( (\widehat{\Omega_{\flat} G} \ltimes Q_{\flat, *} G  ) \rightrightarrows Q_{\flat, *} G \bigr) \]
carries a strict Lie 2-group structure \cite[Proposition 3.5]{murray-robert-wockel}. 
The inclusion  $\iota: S(Q_{\flat, *}G) \hookrightarrow S(QG) $ is a weak equivalence of Lie groupoids\cite[Theorem 3.6]{murray-robert-wockel}. Its quasi-inverse  $\xi$  transfers the Lie 2-group structure from the strict model  $S(Q_{\flat, *}G)$  to  $S(QG)$, yielding a  non-strict Lie 2-group structure on  $S(QG)$. Both  $S(QG)$  and  $S(Q_{\flat, *}G)$  are models for the string 2-group  $\String(G)$; we will work with the coherent model  $S(QG)$  in what follows.

For the explicit construction of the quasi-inverse  $\xi$  using a smooth staircase cut-off function, we refer the reader to    \cite[Section 3]{murray-robert-wockel}. The details are lengthy and not needed for the 2-representation constructions in later sections; only the existence of the coherent model and its biequivalence to a strict model will be used.

\subsection{2-Representations of biequivalent 2-groups} \label{biequiv_2rep}

In this section we show that  the 2-representation theory is invariant under biequivalence of coherent 2-groups, ensuring that the theory depends only on the essential 2-group structure and not on the specific model chosen.

For ordinary groups, an isomorphism  $G \cong G'$  induces an isomorphism of representation categories $\Rep(G) \cong \Rep(G') $. More generally, if  $BG$  and  $BG'$  are equivalent as categories, then  $\Rep(G) \simeq \Rep(G') $ are equivalent as categories. The following theorem is the 2-categorical analogue of this classical fact:  if two coherent Lie 2-groups are biequivalent, then their bicategories of 2-representations are biequivalent. The proof is a direct categorification of the classical proof, with the main technical point being the lifting of the equivalence data (pseudofunctors, strong transformations, and modifications) from the 2-groups to their 2-representation bicategories.

\begin{theorem} \label{biequiv2rep_grp}

    Let $\huaG_{\bullet}$ and $\huaG_{\bullet}'$ denote two biequivalent coherent 2-groups, i.e. there exists a biequivalence between the bicategories $B\huaG_{\bullet}$ and $B\huaG_{\bullet}'$. Then, the bicategories $2\Rep(\huaG_{\bullet})$ and $2\Rep(\huaG_{\bullet}')$ are biequivalent.
    
\end{theorem}

The proof is a direct categorification of the classical argument above. Given a biequivalence between  $B\huaG_\bullet$  and  $B\huaG'_\bullet$, we precompose with the pseudofunctors to obtain functors between the 2-representation bicategories, and post-whisker the coherence data to obtain the required strong transformations and modifications. The verification is lengthy but entirely formal; it is included in Appendix \ref{proof3.2}.

\begin{remark}
The proof of Theorem \ref{biequiv2rep_grp}  uses only the abstract bicategory 
structure of $2\Vect$. In particular, the statement holds for any other 
choice of bicategory used to model 2-vector spaces.
\end{remark}

\begin{remark} \label{biequiv2rep_ex}

Theorem \ref{biequiv2rep_grp} also provides a way to construct examples of $2$-representations  of coherent 2-group. For instance, if the given coherent 2-group $\huaG_{\bullet}$ is biequivalent to another coherent 2-group $\huaG_{\bullet}'$ via the internal equivalence $$\Bigl(  B\huaG_{\bullet} \xrightarrow{\rho}B\huaG_{\bullet}', \quad B\huaG_{\bullet}'\xrightarrow{\iota} B\huaG_{\bullet}, \quad \id_{B\huaG_{\bullet}} \buildrel{\eta}\over\Rightarrow \iota \circ \rho,  \quad 
     \rho\circ \iota \buildrel{\epsilon}\over\Rightarrow \id_{B\huaG_{\bullet}'}\Bigr), $$ then, any $\huaG_{\bullet}'$-representation $ F: B\huaG_{\bullet}'\longrightarrow 2\Vect$ gives a $\huaG_{\bullet}$-representation $F\circ \rho$.


\end{remark}

\subsection{From positive energy representations of $LG$ to 2-representation of $\String(G)$} \label{ex_str2grp_rep}

We now apply the general construction of Section \ref{2grp_basic}  to the coherent string 2-group $\String(G)$. Recall from Section \ref{MRWmodel} that $\String(G)$ is modeled by the coherent Lie 2-group $S(QG):= \bigl( (\widehat{LG} \ltimes QG )\rightrightarrows QG\bigr)$. This fits exactly into the setting of Corollary \ref{lem:CV}. 

We recall the standard notion of positive energy representation of the loop group; see \cite[Chapter 8]{Segal:loop} for details. A positive energy representation of $\widetilde{LG}$ is a unitary representation on a Hilbert space $V$ together with a compatible action of the rotation group $\mathbb{T}$ such that the generator of the rotation action has positive spectrum.

Let $\Rep^+(\widetilde{LG})$ denote the full subcategory of $\Rep(\widetilde{LG})$ consisting of positive energy representations.

For each positive energy representation $(V, \zeta)$ of $\widetilde{LG}$ at level $k$, 
we obtain a $k$-linear category
\[
C_V := \Fun(S(QG), \{V\}/\zeta(\widetilde{LG})),
\]
and Corollary \ref{lem:CV} yields a 2-representation
\[
\mathcal{F}_V = (f_V, f^2_V, f^0_V)
\]
of $\String_k(G)$ on $2\Vect$.  No faithfulness condition on \((V, \zeta)\) 
is required.

\begin{remark}
The construction above extracts the algebraic skeleton of a positive energy 
representation and builds a 2-representation from it. It does not require the 
positive energy representation to be faithful, nor does it preserve the full analytic 
structure (e.g., the topology on $V$ or the unitarity of the $\widetilde{LG}$-action). 
This is analogous to the standard fact that every continuous representation of a Lie 
group defines an algebraic representation of its underlying abstract group: the 
analytic information is forgotten, but the algebraic core remains intact. The 
resulting 2-representation is therefore not a replacement for the original loop group 
representation, but a categorified shadow that encodes its algebraic structure.
\end{remark}


\section{2-Vector bundles and Associated 2-vector bundles} \label{2Vect}

\subsection{
2-Vector bundles over smooth manifolds}\label{2vbconstr}

Motivated by the construction of principal 2-bundles in \cite{Sheng-Xu-Zhu}, we develop a model of 2-vector bundles over smooth manifolds. The construction follows the pattern of stackification: we first define a $(3,1)$-presheaf $2\VectBd$ of pre-2-vector bundles, and then apply the plus construction of Nikolaus and Schweigert \cite{Nikolaus-Schweigert} to obtain a 2-stack $2\VectBd^+$ of genuine 2-vector bundles.



Let $Bicat$ denote the tricategory of small bicategories, pseudofunctors, strong transformations and modifications. For any two bicategories $\huaA$ and $\huaB$, there is a bicategory $\pFun(\huaA, \huaB)$ of pseudofunctors from $\huaA$ to $\huaB$, strong transformations and modifications. We first recall 
    the bicategory of $(3,1)$-presheaves over $\Mfd$. The category $\Mfd$ can be viewed as a sub-tricategory of $Bicat$ where each manifold is regarded as a bicategory with only identity 1-morphisms and identity 2-morphisms, and each smooth map is regarded as a pseudofunctor with trivial laxity constraints. The bicategory of $(3,1)$-presheaves over $\Mfd$ can be interpreted as a sub-bicategory of  $\pFun(\Mfd^{op}, Bicat)$ of strictly unitary pseudofunctors, strong transformations and invertible modifications. A  $(3, 1)$-presheaf $\huaX$ assigns to each manifold $M$ a bicategory $\huaX(M)$, and to each smooth map $U\xrightarrow{\psi} V$ a    
    pseudofunctor \[ \huaX(\psi) :    \huaX(V) \longrightarrow \huaX(U). \]  
    Moreover, there is an extra technical condition on $\huaX$ that it preserves products, i.e.  disjoint union of smooth manifolds.

Next, we formulate the $(3,1)$-presheaf $2\VectBd$ of 2-vector bundles: $\Mfd^{op} \to 2\gpd$, which  is 
constructed in the following way.
\begin{itemize}
\item \textbf{On the level of objects.} For each manifold $U$, the bicategory $2\VectBd(U)$ is the nerve of the following truncated simplicial object:
\begin{itemize}
    \item Objects: trivial bundles $U \times \huaV$ with $\huaV \in 2\Vect$.
    \item 1-morphisms: functor-valued maps $U \to \Auto(\huaV)_0$.
    \item 2-morphisms: natural-isomorphism-valued maps $U \to \Auto(\huaV)_1$.
\end{itemize}

 \item \textbf{On the level of morphisms},  for any morphism $U\xrightarrow{\psi} V$, the associated pseudofunctor $2\VectBd(V) \to 2\VectBd(U)$ is defined by pre-composition with $\psi$. 

\end{itemize}

\begin{remark}

viewed as a pseudofunctor $\Mfd^{op} \rightarrow Bicat$, The $(3, 1)$-presheaf $2\VectBd$  sends a manifold $U$ to the bicategory $$(\{\huaV\}, \Auto(\huaV)_0, \Auto(\huaV)_1)$$ where  $\huaV$ is a 2-vector space,  and sends each smooth map $U\xrightarrow{\psi} V$ to the identity strict functor. Moreover, the laxity constraints are defined as: $2\VectBd^2_{\phi, \psi} = l_{\id}$ for any composable smooth maps $\phi$ and $\psi$; $2\VectBd^0_{\phi} $ is the identity. 
    
\end{remark}

\begin{remark} \label{2vb_ff_rmk}
Note that for any surjective submersion 
$U\xrightarrow{\psi} V$, the pullback $2\VectBd(V) \xrightarrow{\psi^*} 2\VectBd(U)$ is fully faithful.

For any $1-$morphism $f$ and $f'$ in $2\VectBd(V)$, let $\alpha: f\Rightarrow f'$ denote a natural transformation between them, which consists of  a natural transformation  valued map $V \to \Auto(\huaV)_1$, i.e. for each $v\in V$, we have an element  \[\alpha_v: f(v)\rightarrow f'(v)\] in $\Auto(\huaV)_1$. The pullbacks of the natural transformation $\alpha$ is $\alpha\horicirc\id_{\psi}: \psi^*f \Rightarrow \psi^* f'$.
If $\alpha, \alpha': f\Rightarrow f'$ are two different 2-morphisms, so are $\alpha\horicirc\id_{\psi}$ and $\alpha'\horicirc\id_{\psi}$. Thus,
$\psi^*$ on the Hom-categories is faithful. In addition,  
each natural transformation $f\circ \psi\Rightarrow  f'\circ \psi$ can be expressed in the form $\alpha' \horicirc\id_{\psi} $. Thus, $\psi^*$ is fully.

Thus, by \cite[Def 2.12.]{Nikolaus-Schweigert},  
 the $(3,1)$-presheaf $2\VectBd$ is a 2-prestack.

\end{remark}



Now we do a plus construction as in \cite{Nikolaus-Schweigert} to obtain the  $(3,1)$-sheaf $2\VectBd^+$ which describes the entire data of 2-vector bundles, their $1$-morphisms and 2-morphisms.  
For any smooth manifold $M$, let $\{ U_i\}_i$ denote a cover of it and $U(M)_\bullet$ denote the \v{C}ech nerve of $\{ U_i\}$.

\begin{itemize}
    \item An object in $\holim 2\VectBd(U(M)_\bullet)$  is a strictly unitary pseudofunctor from $U(M)_\bullet$ to the bicategory $2\Vect$. 
    This means it assigns data to each level of the \v{C}ech nerve, satisfying certain compatibility conditions. Let's detail these:
\begin{itemize}
    \item On each $U_i$, we have a $2$-vector bundle, which is locally trivial: $U_i\times \huaV$ for some $2$-vector space $\huaV$.
    \item On each intersection $U_{ij} = U_i\cap U_j$, we have a gluing functor $\gamma_{ij}: U_{ij} \to \Auto(\huaV)_0$, which is a functor valued map assigning to each point $x\in U_{ij}$ an automorphism $\gamma_{ij}(x)$ of $\huaV$. Think of this as a functor gluing the fiber of the 2-vector bundles over $U_i$ and $U_j$ together along the overlap.
    \item On each triple intersection $U_{ijk} = U_i\cap U_j\cap U_k$, we have  a natural isomorphism  $
    \gamma_{ij}\circ \gamma_{jk} \xleftarrow[]{\phi_{ijk}} \gamma_{ik}
    $ for each $2$-simplex, which fit together and give 
    a natural-isomorphism valued map $\phi_{ijk}: U_{ijk} \to \Auto(\huaV)_1$. 
 
    \item On each quadruple intersection $U_{ijkl}$, a pentagon condition for $\phi_{ijk}$'s ensures the triple overlaps fit together in a coherent way:
    \begin{equation}
            \xymatrix{(\gamma_{ij}\circ \gamma_{jk})\circ \gamma_{kl} \ar@{=}[dd]& \gamma_{ik}\circ \gamma_{kl} \ar[l]_>>>>>{\phi_{ijk}\circ \id } &\\
            & & \gamma_{il}\ar[lu]_{\phi_{ikl}} \ar[dl]^{\phi_{ijl}} \\
            \gamma_{ij}\circ (\gamma_{jk}\circ \gamma_{kl} ) & \gamma_{ij} \circ \gamma_{jl} \ar[l]^>>>>>{\id\circ \phi_{jkl}},  } 
    \end{equation}
\end{itemize}

\item 
A 1-morphism in $\holim 2\VectBd(U(M)_\bullet)$ from $(U_i\times \Tilde{\huaV}, \Tilde{\gamma}_{ij}, \Tilde{\phi}_{ijk})$ to $(U_i\times \huaV, \gamma_{ij}, \phi_{ijk})$ is a lax transformation. 

\begin{itemize}
    \item For each $U_i$, it consists of a map $\gamma_i: U_i \to \Hom_{2\Vect}(\Tilde{\huaV}, \huaV)$, which is a $1$-morphism between the $2$-vector bundles
    \[
    U_i \times \huaV \xleftarrow{\gamma_i} U_i\times \Tilde{\huaV}.
    \]  
    \item For each double overlap $U_{ij}$, an element $\phi_{ij} \in 2\VectBd(\sqcup_{ij}U_{ij})_2$ which provides a 2-morphism making the following diagram 2-commute.
    \[
    \xymatrix{U_j \times \huaV \ar[dd]_{\gamma_{ij}} & & U_j\times \Tilde{\huaV} \ar[ll]_{\gamma_j} \ar[dd]^{\Tilde{\gamma}_{ij}}\\
    &\buildrel{\phi_{ij}}\over{\Leftarrow}&\\ 
    U_i \times \huaV  & & U_i\times \Tilde{\huaV} \ar[ll]_{\gamma_i} 
    }
    \]
    \item There is a higher coherence condition between these 2-morphisms, expressing the lax naturality of the transformation.   

\begin{equation} \label{2vectbd_1mor_lax_nat}
\xymatrix{ 
 U_k\times \thuaV \ar@/^0.8pc/[rr]^{ \tgamma_{ik} } \ar[dd]_{\gamma_k} && U_i \times \thuaV   \ar[dd]^{ \gamma_i \hspace{0.6cm} = } \\ && \\
  U_k\times \huaV  \ar@{}[rruu]|{\buildrel{\phi_{ik}}\over\Leftarrow }
        \ar@/^0.8pc/[rr]|{\gamma_{ik} }   \ar@/_0.8pc/[rd]_{ \gamma_{jk} }
      &&  U_i \times \huaV  \\
        & U_j \times \huaV    \ar@/_0.8pc/[ru]_{ \gamma_{ij}}   \ar@{}[u]|{\Uparrow { \phi_{ijk}} } & 
}
\xymatrix{ 
  U_k\times \thuaV \ar@/^0.8pc/[rr]^{ \tgamma_{ik}  }  \ar@/_0.8pc/[rd]|{ \tgamma_{jk} }  \ar[dd]_{\gamma_k} & & U_i \times \thuaV  \ar[dd]^{ \gamma_i } \\ & U_j\times \thuaV \ar@/_0.8pc/[ru]|{ \tgamma_{ij}}  \ar[dd]|{ \gamma_j }  \ar@{}[u]|{ \Uparrow  \tphi_{ijk} } & \\ 
   U_k\times \huaV \ar@{}[ru]|{\buildrel{ \phi_{jk} }\over\Leftarrow }
      \ar@/_0.8pc/[rd]_{\gamma_{jk}}
     &&  U_i \times \huaV \\   &   U_j \times \huaV    \ar@{}[ruuu]|{\buildrel{\phi_{ij}}\over\Leftarrow } \ar@/_0.8pc/[ru]_{ \gamma_{ij} } &
}
\end{equation}

\end{itemize}

\item A 2-morphism in $\holim 2\VectBd(U(M)_\bullet)$ from $(\tgamma_i, \tphi_{ij})$ to $(\gamma_i, \phi_{ij})$ is a modification. 
\begin{itemize}
    \item It consists of, for each $U_i$, a $2$-morphism $\phi_i\in 2\VectBd(\sqcup_{ij}U_{ij})_2$ from $\tgamma_i$ to $\gamma_i$;
    \item A coherence condition is required to hold on $U_{ij}$, as shown in \eqref{2mor_modif}, ensuring that the $2$-morphisms are compatible with the gluing functors.

    \begin{equation} \label{2mor_modif}
           \xymatrix{U_j \times \huaV \ar[dd]_{\gamma_{ij}} && U_j\times \Tilde{\huaV} \ar[ll]_{\gamma_j} \ar[dd]^{\Tilde{\gamma}_{ij} \hspace{0.5cm} = }\\
           &\buildrel{\phi_{ij}}\over{\Leftarrow}& \\ 
    U_i \times \huaV  
     \ar@{}@<-2.0ex>[rr]|{\Uparrow \phi_i }
      & &U_i\times \Tilde{\huaV} \ar[ll]_{\gamma_i}    \ar@/^2pc/[ll]|{\tgamma_{i} }  
    }   
           \xymatrix{U_j \times \huaV \ar[dd]_{\gamma_{ij}}      \ar@{}@<+2.0ex>[rr]|{\Uparrow \phi_j } && U_j\times \Tilde{\huaV} \ar[ll]^{\tgamma_j} \ar[dd]^{\Tilde{\gamma}_{ij}  }
           \ar@/_2pc/[ll]|{\gamma_{j} }  
      \\
           &\buildrel{\tphi_{ij}}\over{\Leftarrow}& \\ 
    U_i \times \huaV  
      & &U_i\times \Tilde{\huaV} \ar[ll]^{\tgamma_i}     
    }   
    \end{equation}
\end{itemize}

\end{itemize}

Applying the plus construction as outlined by Nikolaus and Schweigert \cite{Nikolaus-Schweigert}, we obtain a $(3,1)$-sheaf, denoted as $2\VectBd^+:\Mfd^{op}
 \to 2\gpd$, which encodes the complete structure of $2$-vector bundles over smooth manifolds, including their $1$-morphisms and $2$-morphisms. This sheaf is constructed by taking homotopy limits of the bicategories $2\VectBd$ on the \v{C}ech nerve associated to open covers. Specifically, for each smooth manifold $M$, the $2$-groupoid $2\VectBd^+ (M)$ is built as follows:
\begin{itemize}
    \item $2\VectBd^+(M)_0$: an object is a pair $(\{U_i\}, \huaE)$, where $\{U_i\}$ is an open cover of $M$ and $\huaE$ is an element in $\holim 2\VectBd(U(M)_\bullet)_0$;
    \item $2\VectBd^+(M)_1$: a 1-morphism between $(\{U_i\}, \huaE)$ and $(\{\tU_i\}, \thuaE)$ is a pair consisting of a common refinement $\{V_i\}$ of $\{U_i\}$ and $\{\tU_i\}$ and an element $\alpha \in \holim 2\VectBd(V(M)_\bullet)_1 $; 
    \item $2\VectBd^+(M)_2$: a 2-morphism between $(\{V_i \}, \alpha)$ and $(\{\tV_i\}, \talpha)$ consists of a common refinement $\{W_i\}$ of $\{V_i\}$ and $\{\tV_i\}$ and an element $\beta \in \holim 2\VectBd(W(M)_\bullet)_2$. Moreover, $(\{W_i\}, \beta)$ and $(\{\tW_i\}, \tbeta)$ are identified if $\beta$ and  $\tbeta$ agree on a further common refinement of $\{W_i\}$ and $\{\tW_i\}$. 
\end{itemize}

Thus, applying the plus construction of \cite{Nikolaus-Schweigert}, we obtain a 
$(3, 1)$-sheaf $2\VectBd^+ $ on $\Mfd$. The sheaf encodes the complete data of 2-vector bundles over smooth manifolds, including their 1-morphisms and 2-morphisms.

\begin{remark}
\label{2vb_stack_smooth}
The definition of a 2-vector bundle given above should be contrasted with a naive 
assignment of 2-vector spaces to points of a manifold. The plus construction 
enforces the descent condition, which ensures that local data glue smoothly. 
This descent condition is the higher categorical encoding of continuity and 
smoothness for the 2-vector bundle. Without it, the definition would reduce to 
an assignment of algebraic data with no compatibility conditions; with it, we 
obtain a genuine geometric object in the sense of higher stack theory.
\end{remark}

\subsection{ Principal $\huaG_{\bullet}$-bundle and the associated 2-vector bundles} \label{2principalbs}

In this section we construct a model of principal 2-bundles and the 2-vector bundles associated to it.

Let $\huaG_{\bullet}= (\huaG_1 \rightrightarrows \huaG_0)$ be a coherent Lie 2-group. We define a $(3,1)$-presheaf, $\GBd: \Mfd^{op}\to 2\gpd$  as follows.
\begin{itemize}
    \item  \textbf{On the level of objects.} 
For each smooth manifold $U$, the bicategory $\GBd(U)$ consists of: 
\begin{itemize}
    \item $\GBd(U)_0$ consists of a single object, denoted by $U\times \huaG_{\bullet}$;
    \item $\GBd(U)_1$ consists of   1-morphisms  
    $U\times \huaG_{\bullet} \rightarrow U\times \huaG_{\bullet} $ given by maps $U\to \huaG_0$; 
    \item $\GBd(U)_2$ consists of   2-simplices with edges $g_{ij}$ for $0\le i<j\le 2$, which are given by  2-morphisms $f: U\to \huaG_1$ from $g_{02}$ to $g_{01} \horicirc g_{12}$.  
\end{itemize}
\item \textbf{On the level of morphisms.}  Given a morphism $U\xrightarrow{\psi} V$, the associated functor $$\GBd(V) \xrightarrow{\GBd(\psi)} \GBd(U)$$ is given by precomposition with $\psi$. 
\end{itemize} 

\begin{remark}
    Viewed as a pseudofunctor $\Mfd^{op} \rightarrow Bicat$,  the $(3, 1)$-presheaf $\GBd$  sends each smooth manifold $U$ to the delooping bicategory $B\huaG_{\bullet}$,  and each smooth map $U\xrightarrow{\psi} V$ to the identity strict $2$-functor on $B\huaG_{\bullet}$.  Moreover, the laxity constraints are defined as: $2\GBd^2_{\phi, \psi} = l_{\id}$ for any composable smooth maps $\phi$ and $\psi$; $2\GBd^0_{\phi} $ is the identity for any $\phi$. 
    
\end{remark}

Moreover, with plus construction, we arrive at a 2-stack $$\GBd^+: \Mfd^{op}\to 2\gpd,$$ which organizes the information of $\huaG_{\bullet}$-principal bundles, their $1$-morphisms and 2-morphisms together.

\begin{example}[A non-trivial principal $\String_k(G)$-2-bundle]
Let $G$ be a compact, simply connected simple Lie group and fix a level $k \in \mathbb{Z}$. 
Let $\String_k(G)$ be the coherent Lie 2-group model of the string 2-group 
constructed in \cite{murray-robert-wockel}.

Take  a closed smooth manifold $X$ of dimension at least $4$, and  a principal $G$-bundle 
$P_0 \to X$. The obstruction to lifting its structure group to $\String_k(G)$ is given by the     fractional Pontryagin class $\frac{k}{2}p_1(P_0) \in H^4(X,\Z)$. If the obstruction vanishes, $P_0$ carries a level $k$  string structure; a detailed discussion of these obstructions can be found in \cite{BunkeNikolausVolkl2016, kottke_melrose_2021, Waldorf2023Encyclopedia}.

Such a string structure is precisely a lift of the structure group from $G$ to the 
2-group $\String_k(G)$. Equivalently, it gives a principal 
$\String_k(G)$-2-bundle
\[
\huaP \longrightarrow X.
\] in the sense of \cite{Sheng-Xu-Zhu, Bartels:Higher_gauge}. 
This 2-bundle may be  non-trivial even when the underlying  $G$-bundle  $P_0$  is trivial, as the set of string structures forms  a torsor under  the group $H^3(X, U(1))$  and carries additional higher geometric data  (see, e.g., \cite[Corollary 1.3.2]{Waldorf:StringConnections}).
\end{example}

Analogous to how a representation of a Lie group $G$ allows us to associate a vector bundle to a $G$-principal bundle, a 2-representation of the Lie 2-group $\huaG_{\bullet}$ enables us to associate a 2-vector bundle to each  principal $\huaG_{\bullet}$-bundle.
 We show in Remark \ref{rep_cat} that the objects and $1$-morphisms in $2\Rep(\huaG_{\bullet})$ define a category. We first have the conclusion below.

\begin{proposition}
    There is a functor $\huaR_{\huaG_{\bullet}}$ from the category  $2\Rep( \huaG_{\bullet} )$ to  the category $$\sTran(\GBd, 2\VectBd)$$ of strong transformations from $\GBd$ to  $2\VectBd$ and modifications.
\end{proposition}

\begin{proof}
    For each $\huaG_{\bullet}$-2-representation $\rho: B\huaG_{\bullet}\rightarrow 2\Vect$, we can define a strong transformation $\huaR_{\huaG_{\bullet}}(\rho): \GBd \Rightarrow 2\VectBd$ by: \begin{itemize}
        \item for each smooth manifold $U$, $\huaR_{\huaG_{\bullet}}(\rho)_U: \GBd(U) \longrightarrow 2\VectBd(U)$ is defined by $\rho$; 
        \item for each smooth map $U\xrightarrow{\psi} V$, $\huaR_{\huaG_{\bullet}}(\rho)_{\psi}: \rho\circ \id_{B\huaG_{\bullet}}\rightarrow \id_{\rho(\ast)} \circ \rho$ is defined to be $l^{-1}_{\rho}\circ r_{\rho}$. 
    \end{itemize} It is straightforwards to check that $\huaR_{\huaG_{\bullet}}(\rho)$ satisfies the naturality, lax naturality and lax unity.  Thus, It is a well-defined strong transformation. 

    In addition, for any strong transformation $  \xymatrix@1@C+2em{ 
2\Vect &
  \ar@/_1pc/[l]_{\rho}_{}="2"
  \ar@/^1pc/[l]^{\rho'}^{}="3"
  \ar@{<=}"3";"2"_{\alpha}
  B\huaG_{\bullet}}$, we define a modification $\huaR_{\huaG_{\bullet}}(\alpha): \huaR_{\huaG_{\bullet}}(\rho) \longrightarrow \huaR_{\huaG_{\bullet}}(\rho')$. For any smooth manifold $U$, 
  we define a $2$-isomorphism $\huaR_{\huaG_{\bullet}}(\alpha)_U$ to be $\alpha_{\ast}$ exactly, where $\ast$ is the only object in $B\huaG_{\bullet}$. The modification axiom follows by the naturality of $l$ and $r$.  

  In addition, we can check directly that, for any composable strong transformations $\rho\buildrel{\alpha}\over\Rightarrow \rho' \buildrel{\alpha'}\over\Rightarrow  \rho'' $ in $2\Rep (\huaG_{\bullet})$, $ \huaR_{\huaG_{\bullet}} (\alpha'\horicirc \alpha)  = \huaR_{\huaG_{\bullet}}(\alpha')\circ \huaR_{\huaG_{\bullet}}(\alpha)$; for any $\huaG_{\bullet}$-2-representation $\rho$, $\huaR_{\huaG_{\bullet}}(1_{\rho}) = 1_{\huaR_{\huaG_{\bullet}}(\rho)}$.

\end{proof}

Then, by the functoriality of the plus construction,  there is a functor $$\text{Plus}(\GBd, 2\VectBd ): \quad \sTran(\GBd, 2\VectBd) \rightarrow \sTran(\GBd^+, 2\VectBd^+).$$ 

By composing $\text{Plus}(\GBd, 2\VectBd )$ and   $\huaR_{\huaG_{\bullet}}$, we obtain the following result.
\begin{theorem} \label{asso2vectbdl}
    The composition of the  functors $\text{Plus}(\GBd, 2\VectBd )$ and $\huaR_{\huaG_{\bullet}}$  gives a  functor $$\huaR^{+}_{\huaG_{\bullet}}: = \text{Plus}(\GBd, 2\VectBd ) \circ \huaR_{\huaG_{\bullet}}:  \quad  2\Rep(\huaG_{\bullet})\longrightarrow     \sTran(\GBd^+, 2\VectBd^+).$$ 
\end{theorem} 

In particular, for a principal $\huaG_{\bullet}$-2-bundle 
$\huaP \in \GBd^+(X)$ and a 
$\huaG_{\bullet}$-2-representation $\rho \in 2\Rep (\huaG_{\bullet})$, 
evaluating the 
strong transformation $\huaR^{+}_{\huaG_{\bullet}}(\rho)$ on $\huaP$ yields an 
\emph{associated 2-vector bundle}
\[
\huaP \times_{\huaG_{\bullet}} \huaV:= \mathcal{R}_{\huaG_{\bullet}}^+(\rho)(\huaP) \in 2\VectBd^+(X).
\] where $\huaV = \rho(\ast) $.

The construction of this section therefore establishes the dictionary between 
$\huaG_{\bullet}$-2-representations and associated 2-vector bundles.

\section{A model of Equivariant 2-vector bundles}\label{equiv2VB}

In this section we develop a framework for \emph{2-equivariant 2-vector bundles} over 
Lie groupoids. The construction proceeds in two steps. In Section \ref{2vectbd_grpd}, we first 
generalize the notion of 2-vector bundles from smooth manifolds to Lie groupoids, 
using hypercovers to handle descent and obtaining a 2-stack $2\VectBd^+$ 
on the category $\gpd$ of Lie groupoids. In Section \ref{2eq_2vect}, we specialize to the case 
of a Lie groupoid  $X_\bullet$ equipped with an action of a coherent Lie 2-group $\huaG_{\bullet}$, and we define the bicategory 
$2\VectBd^+_{\huaG}(X_\bullet)$ of $\huaG_{\bullet}$-equivariant 2-vector 
bundles over $X_\bullet$. This notion reduces to ordinary 2-vector bundles 
when the $\huaG_{\bullet}$-action is trivial, and provides the natural setting for 
studying categorified gauge theories with symmetries.

\subsection{2-Vector bundles over Lie groupoids } \label{2vectbd_grpd}

To define equivariant 2-vector bundles, we first need to establish the notion of 2-vector bundles over Lie groupoids. Recall that $\Mfd$ denotes the category of smooth manifolds, and $\gpd$ the category of Lie groupoids. To define sheaves and stacks on $\gpd$, we need a suitable Grothendieck topology. Here we choose the Grothendieck pretopology $\tau$ on $\Mfd$ to be the surjective submersions.

While open covers are sufficient for manifolds, we need a more general notion of covering for Lie groupoids: hypercovers. A hypercover, in essence, provides a simplicial resolution of a Lie $n$-groupoid, analogous to a simplicial set resolving a topological space. The successive conditions in Definition \ref{def:hypercover} ensure that the approximation becomes increasingly accurate. Specifically, the $\tau$-covers in the lower degrees provide local trivializations, while the isomorphism at degree $n$ guarantees that the higher homotopy information of the original groupoid is faithfully captured by the hypercover. This construction allows for replacing a potentially complicated Lie $n$-groupoid with a weakly equivalent simplicial object built from simpler components, i.e. manifolds.

\begin{definition}  \label{def:hypercover}For any Lie $n-$groupoid $\Gamma_{\bullet}$ and $X_{\bullet}$, a morphism $\Gamma_{\bullet}\twoheadrightarrow X_{\bullet}$ is a hypercover if, for any $0\leq k\leq n-1$, the natural map from $\Gamma_k$ to the pull-back of the diagram
\[\smap(\partial\Delta[k], \Gamma_{\bullet})\rightarrow \smap(\partial\Delta[k], X_{\bullet})\leftarrow X_k\] is a $\tau-$cover and it is an isomorphism for $k=n$, where $\tau$ is a singleton Grothendieck pretopology on the category $\Mfd$. 
\end{definition}

To establish that our construction of 2-equivariant 2-vector bundles satisfies descent, we can employ hypercovers on the category of Lie groupoids. Hypercovers provide a powerful tool for dealing with descent properties by ensuring equivalence of the corresponding homotopy colimits.

We now extend the construction of the $(3,1)$-presheaf $2\VectBd$ from manifolds to Lie groupoids.

For any Lie groupoid $\Gamma_{\bullet}$, we use the symbol $\Gamma_2$ to denote the fiber product $\Gamma_1\times_{\Gamma_0} \Gamma_1$, and $\Gamma_n : = \Gamma_1\times_{\Gamma_0} \cdots \times_{\Gamma_0} \Gamma_1$    analogously.  We define a bicategory $2\VectBd(\Gamma_{\bullet})$ as follows:

\begin{itemize}
    \item \textbf{Objects} are pairs $(\Gamma_0, \huaV)$, where $\huaV$ is an object of $2\Vect$. We denote such an object by $\Gamma_0 \times \huaV$.
    
    \item \textbf{1-morphisms} between $\Gamma_0 \times \huaV$ and $\Gamma_0 \times \huaV'$ are functor-valued maps $\Gamma_1 \to \Fun(\huaV, \huaV')_0$. (When $\huaV = \huaV'$, this is a map $\Gamma_1 \to \Auto(\huaV)_0$.)
    
    \item \textbf{2-morphisms} between two such 1-morphisms are natural-isomorphism-valued maps $\Gamma_2 \to \Fun(\huaV, \huaV')_1$. 
\end{itemize}

For a morphism of Lie groupoids $\psi_{\bullet}: \Gamma_{\bullet} \to \Lambda_{\bullet}$, the pullback pseudofunctor 
\[
\psi_{\bullet}^*: 2\VectBd(\Lambda_{\bullet}) \longrightarrow 2\VectBd(\Gamma_{\bullet})
\]
is defined by pre-composition with $\psi_{\bullet}$ on each level: objects pull back via $\psi_0$, 1-morphisms via $\psi_1$, and 2-morphisms via $\psi_2$.

\begin{definition}\label{X:stack}
Let $\mathfrak{Y}$ be a presheaf in bicategories on $(\ngrd, \tau)$. 
\begin{enumerate}
\item A presheaf $\mathfrak{Y}$ is called a $\tau-$prestack, if, for every hypercover $\Gamma_{\bullet}\twoheadrightarrow X_{\bullet}$ in $(\ngrd, \tau)$,
the functor of bicategories \[\tau_{\Gamma}: \mathfrak{Y}(X_{\bullet})\rightarrow  \mathfrak{Y}(\Gamma_{\bullet}) \] 
is fully faithful, 
i.e. all the functors on Hom-categories are equivalences of categories.

\item A presheaf $\mathfrak{Y}$ is called a $\tau-$stack, if, for every hypercover  $\Gamma_{\bullet}\twoheadrightarrow X_{\bullet}$ in $(\ngrd, \tau)$, the functor $\tau_{U}$ of bicategories is an equivalence of bicategories.
\end{enumerate}
\end{definition}

The definition of fully faithful bifunctors in the definition is given in \cite[Definition 7.0.1]{JY:Bicat}, \cite[Definition 2.4.9(ii)]{Gabber2004FoundationsFA}.

Moreover, via similar argument to  Remark \ref{2vb_ff_rmk}, we have the conclusion below.

\begin{proposition}
   For each hypercover  $\Gamma_{\bullet}\xrightarrow{\psi_{\bullet}} \Lambda_{\bullet}$ in $\gpd$, the pullback  $$2\VectBd(\Lambda_{\bullet} ) \xrightarrow{\psi^*_{\bullet}} 2\VectBd(\Gamma_{\bullet}) $$ is fully faithful.
   Thus, by Definition \ref{X:stack}, the $(3,1)$-presheaf $2\VectBd$ over $\gpd$ is a 2-prestack.
\end{proposition}

In addition, a hypercover leads to the equivalence of bicategories, as shown in the lemma below. 
 \begin{lemma}\label{keystep3.3}

Let $\tfrakX$ be a $2$-stack on $\gpd$.  Let $X_\bullet \to Y_\bullet$ be a hypercover of Lie groupoid. Then we have the following equivalence of bicategories: 
\begin{equation} \label{eq:y-xy}
    \tfrakX(Y_\bullet) \xrightarrow{\sim} \holim ( \tfrakX(X_\bullet) \Rightarrow \tfrakX(X^{[2]}_\bullet) \dots ) 
\end{equation}
\begin{equation}\label{eq:x-xy}
    \tfrakX(X_\bullet) \xrightarrow{\sim} \holim ( \tfrakX(X_\bullet) \Rightarrow \tfrakX(X^{[2]}_\bullet) \dots ) 
\end{equation}
thus \begin{equation} \label{eq:y-x}
    \tfrakX(Y_\bullet) \xrightarrow{\sim} \tfrakX(X_\bullet) \end{equation}     
 \end{lemma}
\begin{proof} 
\eqref{eq:y-xy} follows from \cite[Theorem 7.5]{Nikolaus-Schweigert} since hypercover of Lie groupoid especially implies a level-wise cover, which is a $\tau$-covering of simplicial manifolds.   \eqref{eq:x-xy} follows from the fact that the diagonal map 
\[
(X_\bullet \Leftarrow X_\bullet \dots )\to (X_\bullet \Leftarrow X^{[2]}_\bullet \dots )  
\]
is a strong equivalence proven in Lemma 8.1 of the same article. Thus 
\[
\holim(\tfrakX(X_\bullet) \Rightarrow \tfrakX(X^{[2]}_\bullet) \dots ) \xrightarrow{\sim} \holim (\tfrakX(X_\bullet) \Rightarrow \tfrakX(X_\bullet) \dots )\xrightarrow{\cong} \tfrakX(X_\bullet) 
\]
\end{proof}

The general procedure of higher sheafification \cite{Nikolaus-Schweigert} produces, 
from any 2-prestack $\tilde{\mathfrak{X}}$ on $\gpd$, a 2-stack $\tilde{\mathfrak{X}}^+$ 
whose descent condition can be verified locally on hypercovers. 

We now apply this procedure to the 2-prestack $2\VectBd$ constructed above. 
For any Lie groupoid $X_{\bullet}$, let $\Gamma_{\bullet}$ denote a hypercover of 
$X_{\bullet}$. The resulting 2-stack $2\VectBd^+$ encodes the complete data of 
2-vector bundles over Lie groupoids, together with their 1-morphisms and 2-morphisms.

\begin{itemize}
    \item An object in $\holim 2\VectBd(\Gamma_\bullet)$  is a pseudofunctor from $\Gamma_\bullet$ to the bicategory $2\Vect$ with trivial lax  unity constraint. 
    
    Explicitly, it is of the form $(\huaV, \alpha, \alpha^2)$, as shown below.
    \begin{itemize}
    \item on $\Gamma_0$, $\Gamma_0\times \huaV$;
    \item on $\Gamma_1$, a functor valued map $\alpha: \Gamma_1 \to \Auto(\huaV)_0$, i.e. for each $\gamma_1\in \Gamma_1$, $\alpha_{\gamma_1}: \huaV\rightarrow \huaV$ is a $1$-morphism;
    \item on  $\Gamma_2$, a natural-isomorphism valued map $\alpha^2: \Gamma_{2} \to \Auto(\huaV)_1$ given by 2-morphisms 
    \[\alpha^2_{\gamma_1, \gamma_2}: \alpha(\gamma_1)\circ \alpha(\gamma_2) \rightarrow \alpha(\gamma_1\circ \gamma_2)\]
    for any $(\gamma_1, \gamma_2)\in \Gamma_2$;
    \item on  $\Gamma_3$, a pentagon condition for $ \alpha_2$:  for any $(\gamma_{1},  \gamma_{2},  \gamma_{3}) \in \Gamma_3$, 
    \begin{equation}
            \xymatrix{\alpha ( (\gamma_{1}\circ \gamma_{2})\circ \gamma_{3} ) \ar@{=}[dd]    & \alpha(\gamma_1\circ \gamma_2 )\horicirc \alpha (\gamma_{3})  \ar[l]_>>>>>{\alpha^2_{\gamma_1\circ \gamma_2,  \gamma_3}}    &\\
            & & \alpha(\gamma_1) \horicirc \alpha(\gamma_2) \horicirc \alpha(\gamma_3)   \ar@/^1pc/[ld]^{\id_{\alpha(\gamma_1)} \horicirc \alpha^2_{\gamma_2,  \gamma_3}}  
            \ar@/_1pc/[lu]_{\alpha^2_{\gamma_1,  \gamma_2}\horicirc \id_{\alpha(\gamma_3)} }    \\
            \alpha( \gamma_{1}\circ (\gamma_{2}\circ \gamma_{3} ))  & \alpha(\gamma_1)\horicirc \alpha(\gamma_2\circ \gamma_3) \ar[l]^>>>>>{\alpha^2_{\gamma_1, \gamma_2\circ \gamma_3} }& }   
    \end{equation}
\end{itemize}

\item 
A 1-morphism in $\holim 2\VectBd(\Gamma_\bullet)$ from $(\Tilde{\huaV}, \Tilde{\alpha}, \Tilde{\alpha}^2)$ to $(\huaV, \alpha, \alpha^2)$ is a lax transformation. 

Explicitly, it consists of: 
\begin{itemize}
    \item $\beta : \Gamma_0\to \Hom_{2\Vect}(\Tilde{\huaV}, \huaV)$ which are given by 1-morphisms  
    \[
 \huaV \xleftarrow{\beta_{\gamma_0}}  \Tilde{\huaV}, 
    \] for each $\gamma_0\in \Gamma_0$; 
    \item an element $\beta  \in 2\VectBd(\Gamma_1)$ which provides a 2-morphism    $\beta_{\gamma_1}$ for each $(\gamma_0\xrightarrow{\gamma_1} \gamma_0')\in \Gamma_1$, making the following diagram 2-commute;
    \[
    \xymatrix{  \huaV \ar[dd]_{\alpha_{\gamma_1}} & & \Tilde{\huaV} \ar[ll]_{\beta_{\gamma_0}} \ar[dd]^{\Tilde{\alpha}_{\gamma_1} }\\
    &\buildrel{\beta_{\gamma_1}}\over{\Leftarrow}&\\ 
     \huaV  & &  \Tilde{\huaV} \ar[ll]_{\beta_{\gamma_0'}} 
    }
    \] 
    \item a higher coherence condition between 2-morphisms, i.e. the lax naturality below.  
For any $\gamma_0\xrightarrow{\gamma_1}   \gamma_0' \xrightarrow{\gamma_1'} \gamma_0''$ in $\Gamma_2$, 
    
\begin{equation} \label{2vectbd_1mor_lax_nat_grp}
\xymatrix{ 
 \thuaV \ar@/^0.8pc/[rr]^{ \Tilde{\alpha}_{\gamma_1'\circ \gamma_1} } \ar[dd]_{\beta_{\gamma_0}} &&   \thuaV   \ar[dd]^{ \beta_{\gamma''_0} \hspace{0.6cm} = } \\ && \\
   \huaV  \ar@{}[rruu]|{\buildrel{\beta_{\gamma_1'\circ \gamma_1}}\over\Leftarrow }
        \ar@/^0.8pc/[rr]|{\alpha_{\gamma_1'\circ \gamma_1} }   \ar@/_0.8pc/[rd]_{ \alpha_{\gamma_1} }
      &&   \huaV  \\
        &   \huaV    \ar@/_0.8pc/[ru]_{\alpha_{\gamma_1'}}   \ar@{}[u]|{\Uparrow { \alpha^{2}_{\gamma_1', \gamma_1} }} & 
}
\xymatrix{ 
    \thuaV \ar@/^0.8pc/[rr]^{ \Tilde{\alpha}_{\gamma_1'\circ \gamma_1}  }  \ar@/_0.8pc/[rd]|{ \Tilde{\alpha}_{\gamma_1} }  \ar[dd]_{\beta_{\gamma_0}} & &  \thuaV  \ar[dd]^{ \beta_{\gamma_0''} } \\ &   \thuaV \ar@/_0.8pc/[ru]|{ \Tilde{\alpha}_{\gamma_1'}}  \ar[dd]|{ \beta_{\gamma_0'} }  \ar@{}[u]|{ \Uparrow  \Tilde{\alpha}^2_{\gamma_1', \gamma_1} } & \\ 
   \huaV \ar@{}[ru]|{\buildrel{ \beta_{\gamma_1} }\over\Leftarrow }
      \ar@/_0.8pc/[rd]_{\alpha_{\gamma_1}}
     &&  \huaV .\\   &   \huaV    \ar@{}[ruuu]|{\buildrel{\beta_{\gamma_1'}}\over\Leftarrow } \ar@/_0.8pc/[ru]_{ \alpha_{\gamma_1'} } &
}
\end{equation}

\end{itemize}

\item A 2-morphism $\phi$ in $\holim 2\VectBd(\Gamma_\bullet)$ from  $\Tilde{\beta}$  to $\beta$   is a modification. Explicitly, it consists of 
\begin{itemize}
    \item For each $\gamma_0\in \Gamma_0$, a $2$-morphism $\phi_{\gamma_0}:  \Tilde{\beta}_{\gamma_0}\rightarrow   \beta_{\gamma_0}$;

    \item a coherence condition held on $\Gamma_1$, as shown in \eqref{2mor_modif}. For each $\gamma_0\xrightarrow{\gamma_1} \gamma_0'$,

    \begin{equation} \label{2mor_modif_grp}
           \xymatrix{  \huaV \ar[dd]_{\alpha_{\gamma_1}} && \Tilde{\huaV} \ar[ll]_{\beta_{\gamma_0}} \ar[dd]^{\Tilde{\alpha}_{\gamma_1} \hspace{0.5cm} = }\\
           &\buildrel{\beta_{\gamma_1}}\over{\Leftarrow}& \\ 
     \huaV  
     \ar@{}@<-2.0ex>[rr]|{\Uparrow \phi_{\gamma_0'} }
      & &  \Tilde{\huaV} \ar[ll]_{\beta_{\gamma'_0}}    \ar@/^2pc/[ll]|{\Tilde{\beta}_{\gamma'_0}}  
    }   
           \xymatrix{  \huaV \ar[dd]_{\alpha_{\gamma_1}}      \ar@{}@<+2.0ex>[rr]|{\Uparrow \phi_{\gamma_0} } &&   \Tilde{\huaV} \ar[ll]^{\Tilde{\beta}_{\gamma_0}} \ar[dd]^{\Tilde{\alpha}_{\gamma_1}  }
           \ar@/_2pc/[ll]|{\beta_{\gamma_0} }  
      \\
           &\buildrel{\Tilde{\beta}_{\gamma_1}}\over{\Leftarrow}& \\ 
     \huaV  
      & & \Tilde{\huaV} \ar[ll]^{\Tilde{\beta}_{\gamma'_0}}     
    }   .
    \end{equation}
\end{itemize}

\end{itemize}

Then, the $(3,1)$-sheaf $2\VectBd^+: \gpd^{op} \to 2\gpd$  consists of:
\begin{itemize}
    \item $2\VectBd^+(X_{\bullet})_0$: an object is a pair $(\Gamma_{\bullet}, \alpha)$, where $\Gamma_{\bullet}$ is a  hypercover of $X_{\bullet}$ and $\alpha$ is an element in $\holim 2\VectBd(\Gamma_{\bullet})_0$;
    \item $2\VectBd^+(X_{\bullet})_1$: a 1-morphism between $(\Gamma_{\bullet}, \alpha)$ and $(\Tilde{\Gamma}_{\bullet}, \Tilde{\alpha})$ is a pair consisting of a common refinement $\Lambda_{\bullet}$ of $\Gamma_{\bullet}$ and $\Tilde{\Gamma}_{\bullet}$ and an element $\beta \in \holim 2\VectBd(\Lambda_\bullet)_1 $; 
    \item $2\VectBd^+(X_{\bullet})_2$: a 2-morphism between $(\Lambda_{\bullet}, \beta)$ and $( \Tilde{\Lambda}_{\bullet}, \Tilde{\beta})$ consists of a common refinement $\Chi_{\bullet}$ of $\Lambda_{\bullet}$ and $\Tilde{\Lambda}_{\bullet}$ and an element $\psi \in \holim 2\VectBd(\Chi_\bullet)_2$. Moreover, $(\Chi_{\bullet}, \psi)$ and $(\Tilde{\Chi}_{\bullet}, \Tilde{\psi})$ are identified if $\psi$ and  $\Tilde{\psi}$ agree on a further common refinement of $\Chi_{\bullet}$ and $\Tilde{\Chi}_{\bullet}$. 
\end{itemize}

Thus, with the plus construction in \cite{Nikolaus-Schweigert}, we describe what are the 2-vector bundles on Lie groupoids, and their 1-morphisms and 2-morphisms. All is organized in the above (3,1)-sheaf $2\VectBd^+$. 


The construction of the 2-stack of principal $\huaG_{\bullet}$-2-bundles on 
Lie groupoids follows the same pattern as the manifold case treated in Section \ref{2principalbs}. 
One first defines a $(3,1)$-presheaf $\GBd_{\gpd}$ on $\gpd$ by 
assigning to each Lie groupoid $\Gamma_\bullet$ the bicategory of trivial principal   
$\huaG_{\bullet}$-2-bundles over $\Gamma_\bullet$, and then applies the plus 
construction to obtain the 2-stack $\GBd_{\gpd}^+$. 
Since the details are entirely analogous to those of Section \ref{2principalbs}, 
we omit them here.

Applying the associated bundle construction of Section \ref{2principalbs}---which 
is valid for principal $\huaG_{\bullet}$-2-bundles over arbitrary Lie groupoids 
by the same functoriality argument---to a principal 2-bundle $\huaP$ over $\Gamma_{\bullet}$ and a 
$\huaG_{\bullet}$-2-representation $\rho$, we obtain a 2-vector bundle
\[
E_\rho := \huaP \times_{\huaG_{\bullet}} \huaV \longrightarrow \Gamma_\bullet
\]
over the Lie groupoid $\Gamma_\bullet$, where $\huaV = \rho(\ast)$.

\begin{example}[Associated 2-Vector Bundle over the Atiyah Groupoid]
We now illustrate how the associated bundle construction of Section \ref{2principalbs} 
interacts with the Atiyah groupoid of a principal 2-bundle.

Let $X$ be a smooth manifold and let $X_{\bullet}$ denote its pair groupoid: the Lie 
groupoid with objects $X$ and a unique morphism $x\rightarrow y$ for any points $x$ and $y$. 
Let $\huaP \to X_{\bullet}$ be a principal $\huaG_{\bullet}$-2-bundle, which is 
equivalently classified by a morphism of Lie groupoids $X_{\bullet} \to B\huaG_{\bullet}$. 

Given a $\huaG_{\bullet}$-2-representation $\rho $, 
the associated bundle construction yields a 2-vector bundle
\[
E_\rho := \huaP \times_{\huaG_{\bullet}} \huaV \longrightarrow X_{\bullet},
\]
whose fiber over each $x \in X$ is equivalent to the 2-vector space 
$\huaV := \rho(\ast)$. 

Consider the Atiyah groupoid $\mathrm{At}(\huaP) \rightrightarrows X_{\bullet}$ of the 
principal 2-bundle $\huaP$. For strict Lie 2-groups, the Atiyah sequence and its 
higher categorical properties have been studied in \cite{CHATTERJEE_CK_2022, Chaudhuri_2024}. 
In our context, we only need the $1$-truncation of the symmetry 2-groupoid: $\mathrm{At}(\huaP)$ is the Lie groupoid whose objects are points 
of $X$, and whose morphisms from $x$ to $y$ are equivalence classes of 
$\huaG_{\bullet}$-equivariant isomorphisms $\huaP_x \to \huaP_y$ between the fibers, where two such 1-morphisms are considered equivalent if they are related by a 
2-isomorphism. This   1-truncated Atiyah groupoid does not require the acting 2-group to be strict; it depends only on the underlying Lie groupoid structure of $\huaG_{\bullet}$, and the fiber symmetries of 
$\huaP$.

There is a canonical forgetful morphism of Lie groupoids 
$\pi : \mathrm{At}(\huaP) \to X_{\bullet}$, which is the identity on objects and 
sends any $\huaG_{\bullet}$-equivariant  isomorphism to the unique morphism in $X_{\bullet}$. 
Pulling back the classifying map $X_{\bullet} \to B\huaG_{\bullet}$ along $\pi$ 
gives a morphism 
\[
\mathrm{At}(\huaP) \longrightarrow B\huaG_{\bullet}.
\]
Composing this with the 2-representation $\rho : B\huaG_{\bullet} \to 2\Vect$ 
yields a 2-vector bundle over the Atiyah groupoid $\mathrm{At}(\huaP)$, providing a natural example of the 
functoriality of our framework.
\end{example}

\subsection{2-equivariant 2-vector bundles over Lie groupoids} \label{2eq_2vect}

In this section we given the construction of a version of 2-equivariant 2-vector bundles over a Lie groupoid. 

In this section, we will use the symbol $F\circ g$, instead of $F\horicirc \id_g$, to denote the horizontal product  of the 2-morphisms 
\[ \xymatrix@1@C+2em{
  Z&
  \ar@/_1pc/[l]_{f}_{}="0"
  \ar@/^1pc/[l]^{f'}^{}="1"
  \ar@{<=}"0";"1"^{F}
  Y &
  \ar@/_1pc/[l]_{g}_{}="2"
  \ar@/^1pc/[l]^{g}^{}="3"
  \ar@{<=}"2";"3"^{\id_g}
  X}\]
when there is no confusion. Similarly, we will use $f\circ G$, instead of $\id_f\horicirc G$ to denote the horizontal product of the 2-morphisms
\[ \xymatrix@1@C+2em{
  Z&
  \ar@/_1pc/[l]_{f}_{}="0"
  \ar@/^1pc/[l]^{f}^{}="1"
  \ar@{<=}"0";"1"^{\id_f}
  Y &
  \ar@/_1pc/[l]_{g}_{}="2"
  \ar@/^1pc/[l]^{g'}^{}="3"
  \ar@{<=}"2";"3"^{G}
  X}\] when there is no confusion.

 Let $\bgpd$ denote the bicategory of Lie groupoids, bibundles, and equivariant bibundle maps. 
Let $$(\huaG_{\bullet}, m, e, a, l, r) $$ be a coherent Lie 2-group. Let  $ X_{\bullet}$ be Lie groupoid acted by  $\huaG_{\bullet}$ in the sense of \cite[Definition 44]{schommer:string-finite-dim}, which consists of:
\begin{itemize}
    \item an action map $\rho: \huaG_{\bullet}\times  X_{\bullet}\longrightarrow  X_{\bullet}$, which is a $1$-morphism in $\bgpd$.
    \item two invertible 2-morphisms 
\begin{align*}
    a_{\rho}: \rho \circ (m \times  \id) &\longrightarrow \rho \circ (\id \times \rho ), \\
    l_{\rho}: \rho \circ (e\times \id) &\longrightarrow \id,
\end{align*} satisfying the coherence condition \cite[Figure 7 and 8]{schommer:string-finite-dim}.

\end{itemize} 

In the section we give a construction of $\huaG_{\bullet}$-equivariant 2-vector bundles over $ X_{\bullet}$ with $\huaG_{\bullet}$ a strict Lie 2-group. We only consider those $\huaG_{\bullet}$-spaces whose action maps are  Lie functors. 
Let $pr_{X_{\bullet}}: \huaG_{\bullet}\times  X_{\bullet}\longrightarrow  X_{\bullet}$ denote the projection to the second factor.

\begin{definition} \label{def:2eq2vectgrpd}
    Let $$\huaE:= (\Gamma_{\bullet}, \quad \huaV, \quad \Gamma_1\xrightarrow{\alpha}\Auto (\huaV)_0,  \quad \alpha^2)$$ be a  2-vector bundle over $ X_{\bullet}$. It is a $\huaG_{\bullet}$-equivariant  $2$-vector bundle over $X_{\bullet}$ if there is 
    \begin{itemize}
        \item an action map $\huaT: pr_{X_{\bullet}}^*\huaE\longrightarrow \rho^*\huaE$, which is a 1-morphism in the bicategory $$2\VectBd^+
(\huaG_{\bullet}\times X_{\bullet});$$
\item two $2$-isomorphisms  of $2$-vector bundles \begin{align*} &\phi_{\huaT}: 
  & \bigl(pr_1\times (\rho\circ pr_{23}) \bigr)^*\huaT\circ pr^*_{23}\huaT  &\Rightarrow      (m\times pr_3 )^* \huaT 
&\mbox{                      over $\huaG_{\bullet}\times \huaG_{\bullet}\times X_{\bullet}$;  }\\
& l_{\huaT}: &\id      &\Rightarrow  (e\times \id)^*\huaT &\mbox{  over $ X_{\bullet}$,   }  \end{align*} which satisfy the coherence conditions below.
    \end{itemize}
\begin{enumerate}
    \item \textbf{Lax Associativity}
    
     Over $\huaG_{\bullet}\times\huaG_{\bullet}\times\huaG_{\bullet}\times     X_{\bullet}$, 

         \begin{equation}\label{equiv_2vect_cocycle_c}
        \mbox{\tiny\xymatrix{         \bigl( (m\circ pr_{12}) \times (\rho\circ pr_{34})  \bigr)^*\huaT \circ pr^*_{34}\huaT \ar[d]_-{ \bigl( (m\circ pr_{12}) \times pr_{34} \bigr)^*\phi_{\huaT} }   
       &\biggl( \Bigl(pr_1 \times \rho\circ \bigl(pr_2\times (\rho\circ pr_{34})   \bigr) \Bigr)^*\huaT  \circ \bigl(pr_2\times (\rho\circ pr_{34})\bigl)^*\huaT \biggr)  \circ pr^*_{34}\huaT 
        \ar[d]^{a_{\huaG^3X}} \ar@/_2pc/[l]_-{ \bigl(  pr_{12} \times (\rho \circ pr_{34}) \bigr)^*  \phi_{\huaT} \circ   pr^*_{34}\huaT }\\ 
        \biggl(\Bigl( m\circ \bigl(  (m\circ pr_{12}) \times pr_3  \bigr)  \Bigl) \times pr_4 \biggl)^*\huaT   &\Bigl(pr_1 \times \rho\circ \bigl(pr_2\times (\rho\circ pr_{34})   \bigr) \Bigr)^*\huaT  \circ\biggl( \bigl(pr_2\times (\rho\circ pr_{34})\bigl)^*\huaT  \circ pr^*_{34}\huaT  \biggr) 
        \ar[d]^-{ (pr_1\times a_{\rho})^*\huaT \circ  \bigl( pr^*_{234} \phi_{\huaT} \bigr) }\\
        \biggl( \Bigl(  m\circ \bigl( pr_1 \times (m\circ pr_{23}) \bigr)  \Bigr) \times pr_4 \biggr)^* \huaT  \ar[u]^-{(a_{\rho}\times pr_4 )^* \huaT} &   \Bigl( pr_1\times \rho\circ \bigl( (m\circ pr_{23}) \times pr_4  \bigr)    \Bigr)^*\huaT \circ    \bigl( (m\circ pr_{23}) \times pr_4  \bigr) ^*\huaT  \ar@/^2pc/[l]^-{\bigl(pr_1\times (m\circ pr_{23}) \times pr_4 \bigr)^* \phi_{\huaT}}  \\
      } }
    \end{equation} where $a_{\huaG^3X}$ is the associator in the bicategory $2\VectBd^+(\huaG_{\bullet}\times\huaG_{\bullet}\times\huaG_{\bullet}\times X_{\bullet})$.

 \item \textbf{Lax Unity}

 Over $\huaG_{\bullet}\times X_{\bullet}$, 
    \begin{equation} \label{equiv_2vect_lr_c}
        \xymatrix{ (pr_1\times e \times pr_2)^* \bigl( (m\times pr_3 )^*\huaT \bigr)   \ar[d]_-{\cong} 
        && (pr_1\times e \times pr_2)^* \Bigl( \bigl(pr_1\times (\rho\circ pr_{23})\bigr)^*\huaT\circ pr^*_{23}\huaT \Bigr) \ar[d]^-{\cong} 
        \ar[ll]_-{ (pr_1\times e \times pr_2)^* \phi_{\huaT} } 
        \\ 
        \Bigl( \bigl(m\circ (pr_1\times e) \bigr)  \times pr_2 \Bigr)^*\huaT   &&   \Bigl(pr_1 \times \bigl( \rho \circ (e\times pr_2)\bigr) \Bigr)^*\huaT \circ  (e\times pr_2)^*\huaT  \\
         &\huaT\ar[lu]|{(r\times pr_2)^*\huaT }   \ar[ru]|{ (pr_1\times l_{\rho})^* \huaT \circ  
         l_{\huaT} } &}
    \end{equation} \end{enumerate}
\end{definition}

\begin{remark} \label{ex_equiv_2vectbdl}
     If $\huaG_{\bullet}$ is a  Lie group $G$ and $X_{\bullet}$ is a smooth manifold $X$, the commutative diagram  \eqref{equiv_2vect_cocycle_c} reduces to the cocycle condition. More explicitly, over each point $(g, x)\in G\times X$, $\huaT_{(g,x)}$ is a map from $\huaV_x$ to $\huaV_{g\cdot x} = (g^*\huaV)_x$. 
     Over each $(g_1, g_2, x)\in G\times G\times X$, $(\phi_{\huaT})_{(g_1, g_2, x)}$ is a 2-isomorphism from 
$g_2^*\huaT_{g_1, x} \circ \huaT_{g_2, x}$ to $\huaT_{g_1g_2, x}$.  In addition, the coherence condition \eqref{equiv_2vect_cocycle_c} is exactly
\begin{equation}
    (\phi_{\huaT})_{g_1g_2, g_3, x} \circ (g_3^*(\phi_{\huaT})_{g_1, g_2, x} \times \id) = (\phi_{\huaT})_{g_1, g_2g_3, x} \circ (\id \times (\phi_{\huaT})_{g_2, g_3, x}). 
\end{equation}
And, the 
 coherence condition \eqref{equiv_2vect_lr_c} over each point $(g, x)\in G\times X$ is given by the equation \begin{equation}
(\phi_{\huaT})_{g, e, x} = \id_{\huaT_{g, x}}\horicirc (l_{\huaT})_{x}.
\end{equation}
\end{remark}

\bigskip

Next we define $1$-morphisms between $\huaG_{\bullet}$-equivariant 2-vector bundles over $X_{\bullet}$.    Let $$\huaE:= (\Gamma_{\bullet}, \quad \huaV, \quad  \Gamma_1\xrightarrow{\alpha}\Auto (\huaV)_0, \quad  \alpha^2)$$  denote a $\huaG_{\bullet}$-equivariant $2$-vector bundle over $X_{\bullet}$ with     \begin{itemize}
        \item action map $\huaT: pr_{X_{\bullet}}^*\huaE\longrightarrow \rho^*\huaE$;
\item two $2$-isomorphisms \begin{align*} &\phi_{\huaT}: 
  & \bigl(pr_1\times (\rho\circ pr_{23}) \bigr)^*\huaT\circ pr^*_{23}\huaT  &\Rightarrow      (m\times pr_3 )^* \huaT 
&\mbox{                      over $\huaG_{\bullet}\times \huaG_{\bullet}\times X_{\bullet}$;  }\\
& l_{\huaT}: &\id      &\Rightarrow  (e\times \id)^*\huaT &\mbox{  over $ X_{\bullet}$,   }  \end{align*} 
    \end{itemize}

 And let $$\huaE':= (\Gamma'_{\bullet}, \quad \huaV', \quad \Gamma'_1\xrightarrow{\alpha'}\Auto (\huaV')_0, \quad  \alpha^{'2})$$ denote a $\huaG_{\bullet}$-equivariant $2$-vector bundle over $X_{\bullet}$ with     \begin{itemize}
        \item action map $\huaT': pr_{X_{\bullet}}^*\huaE'\longrightarrow \rho^*\huaE'$;
\item two $2$-isomorphisms \begin{align*} &\phi_{\huaT'}: 
  & \bigl(pr_1\times (\rho\circ pr_{23}) \bigr)^*\huaT'\circ pr^*_{23}\huaT'  &\Rightarrow      (m\times pr_3 )^* \huaT' 
&\mbox{                      over $\huaG_{\bullet}\times \huaG_{\bullet}\times X_{\bullet}$;  }\\
& l_{\huaT'}: &\id      &\Rightarrow  (e\times \id)^*\huaT' &\mbox{  over $ X_{\bullet}$,   }  \end{align*} 
    \end{itemize}

\begin{definition}

A $\huaG_{\bullet}$-equivariant $1$-morphism $\vartheta$ from $\huaE$ to $\huaE'$ consists of:
\begin{itemize}
    \item a $1$-morphism $\beta$ from $\huaE$ to $\huaE'$ in   $2\VectBd^+(X_{\bullet})$; 
    \item a $2$-isomorphism $\zeta$  in  $2\VectBd^+(\huaG_{\bullet}\times X_{\bullet})$: 
\[ \xymatrix{ pr_{X_{\bullet}}^*\huaE \ar[d]_{pr^*_{X_{\bullet}} \beta} \ar[r]^{\huaT} & \rho^*\huaE\ar[d]^-{\rho^* \beta }  \\
pr_{X_{\bullet}}^*\huaE' \ar@{}@<-0.2ex>[ru]|{\buildrel{ \zeta  }\over\Rightarrow} \ar[r]_-{\huaT'} &\rho^*\huaE'   } \] which satisfies the coherence conditions below.
\end{itemize}
\begin{enumerate}
    \item \textbf{Lax Naturality}

 Over $\huaG_{\bullet}\times \huaG_{\bullet}\times X_{\bullet}$, 

 \begin{equation}
   \mbox{\tiny\xymatrix{  & \bigl(pr_1\times( \rho\circ pr_{23} )\bigl)^*\huaT' \circ  pr_{23}^*\bigl( \huaT'\circ pr_{X_{\bullet}}^*\beta \bigl) \ar@/^2pc/[rd]^-{\bigl(pr_1\times( \rho\circ pr_{23} )\bigl)^*\huaT'  \circ pr^*_{23}\zeta}  &
   \\
   & \Bigl( \bigl(pr_1\times( \rho\circ pr_{23} )\bigr)^*\huaT' \circ  pr_{23}^* \huaT' \Bigr)\circ pr_{3}^*\beta  \ar[u]^-{a_{\huaG^2X}}\ar[d]_-{\phi_{\huaT'}\circ pr_{3}^*\beta } &
   \Bigl(pr_1\times( \rho\circ pr_{23} )\Bigl)^*\huaT' \circ  \Bigl( (\rho\circ pr_{23})^*\beta \circ pr_{23}^*\huaT \Bigl)     
   \\ 
   &  (m\times pr_3)^*\huaT'\circ pr_{3}^*\beta \ar[d]_-{\cong}   &    \Bigl(\bigl(pr_1\times( \rho\circ pr_{23} )\bigl)^*\huaT' \circ   (\rho\circ pr_{23})^*\beta \Bigl) \circ pr_{23}^*\huaT  
   \ar[d]_-{\bigl(pr_1\times( \rho\circ pr_{23} )\bigl)^*\zeta \circ pr_{23}^*\huaT  }    \ar[u]^-{ a_{\huaG^2X}} 
 \\ 
 & (m\times pr_3)^*(\huaT'\circ pr_{X_{\bullet}}^*\beta) \ar[d]_-{ (m\times pr_3)^* \zeta} & \Bigl( \bigl(\rho\circ (pr_1\times( \rho\circ pr_{23} ))\bigl)^*\beta \circ \bigl(pr_1\times( \rho\circ pr_{23} )\bigl)^*\huaT \Bigl)    \circ pr_{23}^*\huaT   \ar[d]_-{a_{\huaG^2X} } \\ 
 & (m\times pr_3)^*(\rho^*\beta\circ \huaT)  \ar[d]_-{\cong} &  \bigl(\rho\circ (pr_1\times( \rho\circ pr_{23} ))\bigl)^*\beta \circ  \Bigl( \bigl(pr_1\times( \rho\circ pr_{23} )\bigl)^*\huaT  \circ pr_{23}^*\huaT  \Bigl)  \ar[d]_-{\bigl(\rho\circ (pr_1\times( \rho\circ pr_{23} ))\bigl)^*\beta  \circ \phi_{\huaT}} \\ 
 &  \bigl(\rho\circ (m\times pr_3)\bigl)^*\beta \circ (m\times \pr_3)^*\huaT \ar@/_2pc/[r]_-{a_{\rho}^*\beta \circ  (m\times \pr_3)^*\huaT }&  \bigl(\rho\circ (pr_1\times( \rho\circ pr_{23} ))\bigl)^*\beta \circ  (m\times \pr_3)^*\huaT  }   }
 \end{equation} where $a_{\huaG^2X}$ is the associator in the bicategory $2\VectBd^+ (\huaG_{\bullet}\times\huaG_{\bullet}\times X_{\bullet})$. 

\item \textbf{Lax Unity}

Over $\huaG_{\bullet}\times X_{\bullet}$, 

\begin{equation}
    \xymatrix{ (e\times pr_2)^*\huaT'\circ \beta \ar[dd]_-{\cong }& \id\circ \beta  \ar[l]_-{l'_{\huaT'}\circ \beta}  \\
     &   \beta \circ \id \ar[u]^-{ l_{\huaG X}^{-1}\circ r_{\huaG X}} \ar[d]_-{\beta \circ l_{\huaT} } \\  
    (e\times pr_2)^*(\huaT'\circ pr_{X_{\bullet}}^*\beta) \ar[d]_-{ (e\times pr_2)^*\zeta} & \beta \circ (e\times pr_2)^*\huaT   \\ 
     (e\times pr_2)^* (\rho^*\beta \circ  \huaT   )  \ar[r]_-{\cong} & (\rho\circ (e\times pr_2))^*\beta \circ  (e\times pr_2)^* \huaT   \ar[u]^-{l_{\rho}^*\beta \circ (e\times pr_2)^*\huaT}\\ 
    }
\end{equation}
where $l_{\huaG X}$ and  $r_{\huaG X}$ are the left unitor and  the right unitor in the bicategory $2\VectBd^+ (\huaG_{\bullet}\times X_{\bullet})$. \end{enumerate}
    
\end{definition}

\begin{definition}
    Let $$\vartheta = (\beta, \zeta) \mbox{ and   }\vartheta' = (\beta', \zeta')$$ denote two $\huaG_{\bullet}$-equivariant $1$-morphisms from $\huaE$ to $\huaE'$. A $\huaG_{\bullet}$-equivariant $2$-morphism $\varrho: \vartheta\longrightarrow \vartheta'$ is a $2$-isomorphism $\varrho: \beta\longrightarrow \beta'$ in $2\VectBd^+(X_{\bullet}) $ such that the diagram below commutes.
    \begin{equation}
        \xymatrix{ \huaT' \circ pr^*_{X_{\bullet}}\beta \ar[rr]^{\zeta}  \ar[d]_-{ \huaT'\circ pr^*_{X_{\bullet}}\varrho }   && \rho^*\beta \circ \huaT 
        \ar[d]^-{ \rho^* \varrho \circ \huaT}\\ 
        \huaT' \circ pr^*_{X_{\bullet}}\beta' \ar[rr]^{\zeta'}  && \rho^*\beta \circ \huaT}
    \end{equation}

\end{definition}

\begin{definition}
    The $\huaG_{\bullet}$-equivariant 2-vector bundles over $X_{\bullet}$, $\huaG_{\bullet}$-equivariant $1$-morphisms, and $\huaG_{\bullet}$-equivariant $2$-morphisms, form the bicategory 
    $2\VectBd^+_{\huaG_{\bullet}} (X_{\bullet}) $   of $\huaG_{\bullet}$-equivariant 2-vector bundles over $X_{\bullet}$.
\end{definition}

\begin{remark}[Coherent 2-groups and equivariance]
\label{rem:coherent-equiv}
The definition of 2-equivariant 2-vector bundles presented in this section assumes 
that the acting 2-group $\huaG_{\bullet}$ is \emph{strict}, i.e., a group 
object in the 1-category $\gpd$ of Lie groupoids. For a coherent Lie 
2-group (such as the string 2-group model of Section \ref{MRWmodel}), the 
construction of an intrinsic 2-equivariant 2-vector bundle theory requires a more 
elaborate framework, as coherent 2-groups are not group objects in $\gpd$. 
One expects that, by passing through a biequivalent strict model (which always 
exists for the string 2-group by \cite{murray-robert-wockel}), the resulting 
equivariant bundles should be suitably transportable, but a detailed development 
of this idea is beyond the scope of the present paper and is left for future work.
\end{remark}

\begin{example}
\begin{enumerate}
    \item 
    If $X_{\bullet}$ is the trivial groupoid, then, the bicategory  $$2\VectBd^+_{\huaG_{\bullet}} (X_{\bullet}) $$  is exactly the bicategory $2\Rep(\huaG_{\bullet})$  of $\huaG_{\bullet}$-representations.

\item 
    If the 2-group $\huaG_{\bullet}$ is a discrete bicategory and is given by a Lie group $G$, and the Lie groupoid $X_{\bullet}$ is given by a smooth manifold $X$, then the bicategory $2\VectBd^+_{\huaG_{\bullet}} (X_{\bullet}) $ is exactly the discrete bicategory of $G$-equivariant vector bundles over $X$. \end{enumerate}
\end{example}



\section{Related Work and Future Directions}\label{app_add}

The framework developed in this paper---2-representations of coherent Lie 2-groups, 
associated 2-vector bundles, and their equivariant generalizations---sits at the 
intersection of several active areas of research in higher geometry and mathematical 
physics. In this final section, we situate our results in the broader literature by 
comparing them with related constructions of string 2-group representations and flat 
string structures. We then outline a number of promising directions for future 
investigation.

\subsection{Discussion and expectations}

\subsubsection{Relation to the Drinfel'd centre}

Weis proved in \cite{weis2022drinfeld} that when $G$ is compact connected, the Drinfel'd centre of the string 2-group $\String_k(G)$ is the invertible part of the category of positive energy representations of $LG$ at level $k$. Our construction of a 2-representation $F_{k, V}: B\huaG_{\bullet} \to 2\Vect $ from a positive energy representation $V$ (Section   \ref{ex_str2grp_rep}) can be seen as a geometric incarnation of this algebraic correspondence. Moreover,  the associated 2-vector bundle construction in Theorem \ref{asso2vectbdl} is a categorification of this correspondence.

\subsubsection{Relation to recent work on 2-vector bundles}

Kristel, Ludewig, and Waldorf  \cite{KLW2Vect22} developed a theory of 2-vector bundles based on the bicategory of algebras, bimodules, and intertwiners, showing that bundle gerbes and ordinary algebra bundles appear as  sub-bicategories. Their model differs from ours (algebras/bimodules vs. $k$-prelinear categories), but both frameworks provide a unified approach to higher geometric structures. A comparison of these two models in the context of string 2-group representations would be an interesting direction for future research.

\subsubsection{Representations of the string 2-group}

Kristel, Ludewig, and Waldorf  constructed in \cite{KLW2Rep22} a representation of the string 2-group on a 2-vector space modeled on von Neumann algebras, aiming to establish it as a categorification of the spinor representation. Their representation lives in the analytic setting (hyperfinite type $\text{III}_1$ factor), while ours lives in the algebraic \& geometric setting. Understanding the relationship between these two representations and whether they are equivalent under a suitable 2-version of Morita equivalence is a promising avenue for future work.

\subsubsection{Flat string structures}

Berwick-Evans, Cliff, Murray, Nakade, and Phillips \cite{berwick2025flat} constructed a bicategory of flat 2-group bundles over differentiable stacks and characterized flat string structures in terms of ordinary $G$-bundles together with a trivialization of a certain 2-gerbe. 
Our framework suggests a natural higher analogue of the classical correspondence between flat bundles and fundamental groupoid representations. Recall that the classifying map of a 2-vector bundle $E_\rho$ constructed via the associated 2-vector bundle construction (Theorem \ref{asso2vectbdl}) is given by the composition
\[
X \rightarrow B\huaG_{\bullet} \xrightarrow{\ \rho\ } 2\Vect ,
\]
where the first map gives  a principal $\huaG_{\bullet}$-2-bundle $\huaP$ and $\rho$ is a $\huaG_{\bullet}$-2-representation. By analogy with the classical case, where a flat bundle corresponds to a factorization through the fundamental groupoid $\Pi_1(X)$, we propose that the 2-vector bundle $E_\rho$ should be called \emph{flat} if its classifying pseudofunctor factors through the fundamental 2-groupoid $\Pi_2(X)$. This global flatness condition can be ensured by requiring both that the underlying principal 2-bundle $\huaP$ is flat (as characterized in \cite{berwick2025flat}) and that the 2-representation $\rho$ is suitably ``flat'' in the sense that it factors through $\Pi_2(B\huaG_{\bullet})$.

\subsection{Future directions} The construction presented in this section suggests several avenues for further work:

\begin{enumerate}
    \item \textbf{Fusion product}. One of the main motivations in \cite{murray-robert-wockel} was to define a fusion product on 2-representations of the string 2-group corresponding to the fusion product of positive energy representations. While the principal 2-bundle and associated 2-vector bundle construction in this paper provide the natural geometric setting, the construction of such a fusion product remains open. A possible approach is to adapt the convolution product on the Atiyah Lie groupoid of a principal 2-bundle, with the higher coherence conditions bringing new technical obstacles, which require further investigation.

\item  \textbf{Higher Dixmier--Douady theory}. In \cite{evans2022equivariant}  Evans and Pennig  studied equivariant higher twisted $K$-groups of $SU(n)$ for twists arising from exponential functors, extending the classical Dixmier--Douady classification of twisted $K$-theory. The framework of principal 2-bundles and associated 2-vector bundles developed in this paper provides a categorified geometric model for higher twists classified by non-abelian gerbes. It is a natural question whether these two approaches to higher twisted $K$-theories can be unified under a common categorical framework, which is left for future investigation.

\item \textbf{Physical applications}.  

 In string theory, the string 2-group plays a central role in the Green--Schwarz anomaly cancellation mechanism: the cancellation condition is equivalently the statement that the structure group of the tangent bundle lifts to the string 2-group. The framework of principal 2-bundles developed in this paper provides a natural higher geometric setting for  such anomaly cancellation. Furthermore, in extended topological field theories, surface defects are known to be classified by 2-vector spaces \cite{KapustinWitten2006, Lurie08_TFT}. The associated 2-vector bundles constructed via  2-representations (Theorem \ref{asso2vectbdl}) offer a geometric model for the state spaces of such surface defects. Finally, the extension of these constructions to the equivariant setting (Section \ref{2eq_2vect}) provides the geometric data for studying such systems in the presence of background higher gauge symmetries, potentially linking to equivariant elliptic cohomology and higher twisted K-theory.

\item \textbf{Elliptic cohomology}. 
Positive energy representations of the loop group are known to give rise to modular forms \cite{brylinski1990representations} and to equivariant elliptic cohomology \cite{kitchloo2019quantization}. The interpretation of each  positive energy  representation as a 2-representations of the string 2-group (Section \ref{ex_str2grp_rep}) suggests that this elliptic information may also be captured within the 2-representation framework. A precise geometric link remains to be developed.

On the geometric side, the relationship between 2-vector bundles and elliptic cohomology has been explored from a homotopy-theoretic perspective in \cite{bdr}, and a twisted quasi-elliptic cohomology was constructed from the $K$-theory of the loop space in \cite{huan2025}. The associated  2-representation   constructed in Corollary \ref{lem:CV} offers a new categorical perspective on this circle of ideas, and we plan to develop the corresponding $2K$-theory in future work to establish a direct geometric link,  thereby contributing to a higher categorical realization of the Stolz--Teichner program \cite{st} for the string 2-group.

\end{enumerate}

In summary, this paper has developed a systematic framework for 2-representations 
and 2-vector bundles associated with coherent Lie 2-groups, with the string 2-group 
as the central example. The functorial dictionary between 2-representations and 
associated 2-vector bundles, together with its equivariant extension, provides a 
geometric setting for studying higher symmetries in gauge theory and representation 
theory. The directions outlined in Section \ref{app_add} suggest several avenues for 
further investigation; we intend to explore these in future work.

\appendix

\section{Proof of Theorem \ref{biequiv2rep_grp}}
\label{proof3.2}

This appendix contains the full proof of Theorem \ref{biequiv2rep_grp}. The proof is a direct categorification of the classical argument that equivalent groups have equivalent representation categories: given a biequivalence of 2-groups, we construct a biequivalence of their 2-representation bicategories by precomposition and post-whiskering. The verification is lengthy but entirely formal, using only the bicategorical structure of $2\Vect$ and the fact that pseudofunctors, strong transformations, and modifications form a bicategory.
\begin{proof}[Proof of Theorem \ref{biequiv2rep_grp}]
    Let $$\Bigl(  B\huaG_{\bullet} \xrightarrow{\rho}B\huaG_{\bullet}', \quad B\huaG_{\bullet}'\xrightarrow{\iota} B\huaG_{\bullet}, \quad \id_{B\huaG_{\bullet}} \buildrel{\eta}\over\Rightarrow \iota \circ \rho,  \quad 
     \rho\circ \iota \buildrel{\epsilon}\over\Rightarrow \id_{B\huaG_{\bullet}'}\Bigr) $$ denote a biequivalence from $B\huaG_{\bullet}$ to $B\huaG_{\bullet}'$, where $\rho$ and $\iota$ are both pseudofunctors, 
     $\eta$ and $\epsilon$ are both strong transformations. Moreover, there exist strong transformations $ \iota \circ \rho
      \buildrel{\eta'}\over\Rightarrow \id_{B\huaG_{\bullet}} $ and   $\id_{B\huaG_{\bullet}'} \buildrel{\epsilon'}\over\Rightarrow \rho\circ \iota $ together with invertible modifications
       \[ \Gamma: 1_{\id_{B\huaG_{\bullet}}} \xrightarrow{\cong} \eta'\horicirc \eta, \quad \Gamma': \eta\horicirc \eta' \xrightarrow{\cong} 1_{\iota \circ \rho}\]
     and \[ \Chi: 1_{\rho\circ \iota}\xrightarrow{\cong} \epsilon' \horicirc \epsilon, \quad \Chi': \epsilon\horicirc \epsilon'\xrightarrow{\cong}
     1_{\id_{B\huaG_{\bullet}'}}.\]

     \bigskip
     
     From this,  we can construct a biequivalence from $2\Rep( \huaG_{\bullet}' )$ to $2\Rep( \huaG_{\bullet} )$. We will now outline each component step by step in detail.

     \textbf{Part I: the construction of the pseudofunctors $\rho^*$ and $\iota^*$}

     By precomposing with $\rho$, we can define a strict functor \[ \rho^*: 2\Rep (\huaG_{\bullet}' )\rightarrow 2\Rep (\huaG_{\bullet}). \]
     \begin{itemize}
         \item Object level: a pseudofunctor $B\huaG_{\bullet}'\xrightarrow{F} 2\Vect$ is mapped to the composition $B\huaG_{\bullet} \xrightarrow{\rho} B\huaG_{\bullet}' \xrightarrow{F}2\Vect$;
         \begin{itemize}
             \item $F\circ \rho$ maps the only object $\ast$ in $B\huaG_{\bullet}$ to $F(\ast)$; 
             \item $F\circ \rho$ maps each $1$-morphism $g_0$ in $B\huaG_{\bullet}$ to $F(\rho(g_0))$, and each $2$-morphism $g_1$ in $B\huaG_{\bullet}$ to $F(\rho(g_1))$;
             \item The laxity constraints are defined by: for any $1$-morphisms $h$ and $g$ in $B\huaG_{\bullet}$,
             \begin{align*}
                 (F\circ \rho)^2_{h, g}  & = F(\rho^2_{h, g}) \circ F^2_{\rho(h), \rho(g)}; \\
                 (F\circ \rho)^0_{\ast} & = F(\rho^0_{\ast}) \circ F^0_{\ast}.
             \end{align*} As shown in \cite[Definition 4.1.26, Lemma 4.1.29]{JY:Bicat}, these data above give a well-defined lax functor.  
             
         \end{itemize}
         \item $1$-Morphism level:  a natural isomorphism $\xymatrix@1@C+2em{
  B\huaG_{\bullet}'  \ar@/^1pc/[r]^{F}_{}="4"
  \ar@/_1pc/[r]_{F'}^{}="5"
  \ar@{=>}"4";"5"^{\alpha}
  & 2\Vect
}$ is mapped to the pre-whiskering of $\alpha$ with $\rho$, $$\xymatrix@1@C+2em{B\huaG_{\bullet}\ar[r]^{\rho}
 & B\huaG_{\bullet}'  \ar@/^1pc/[r]^{F}_{}="4"
  \ar@/_1pc/[r]_{F'}^{}="5"
  \ar@{=>}"4";"5"^{\alpha}
  & 2\Vect,
}$$ Explicitly, $\rho^*\alpha: F\circ \rho \Rightarrow F' \circ \rho $ is the strong transformation defined as below. 
\begin{itemize}
    \item The $1$-morphism $\rho^*\alpha_{\ast}:= \alpha_{\rho(\ast)} = \alpha_{\ast}$.
    \item For each $1$-morphism $g$ in $B\huaG_{\bullet}$, $\rho^*\alpha_g:  F'(\rho(g))\circ \alpha_{\ast} \Rightarrow \alpha_{\ast} \circ F(\rho(g))$ is defined to be $\alpha_{\rho(g)}$.
    \item $\bigl(\alpha_{\ast},   \alpha_{\rho(-)}\bigr) $ indeed define a natural isomorphism, as given in \cite[Def.11.1.1(1)]{JY:Bicat}.
    We check the coherence conditions below.
    \begin{enumerate}
        \item \textbf{Lax Unity} 

        \begin{equation}
            \xymatrix{ \id_{F'\circ \rho(\ast)} \circ \alpha_{\rho(\ast)} \ar[r]^>>>>>{l_{\alpha_{\rho(\ast)}}}  \ar[d]_-{F_{\rho(\ast)}^{'0} \horiprod 1_{\alpha_{\rho(\ast)}}} & \alpha_{\rho(\ast)} \ar@{}[d]|{\circlearrowright} & \alpha_{\rho(\ast)}\circ \id_{F\circ \rho(\ast) } \ar[d]^-{1_{\alpha_{\rho(\ast)}} \horiprod F^0_{\rho(\ast)} } \ar[l]_-{r_{\alpha_{\rho(\ast)}} } \\
            (F'\circ \id_{\rho(\ast)})\circ \alpha_{\rho(\ast)} \ar[rr]|{\alpha_{\id_{\rho(\ast)}} }   \ar[d]_-{F'(\rho^0_{\ast})\horiprod 1_{\alpha_{\rho(\ast)}}} 
            && \alpha_{\rho(\ast)} \circ F(\id_{\rho(\ast)})  \ar[d]^{1_{\alpha_{\rho(\ast)}} \horiprod F(\rho^0_{\ast})} \\
            (F'\circ \rho)(\id_{\ast}) \circ \alpha_{\rho(\ast)} \ar[rr]_-{\alpha_{\rho(\id_{\ast})}}  \ar@{}[rru]|{\circlearrowright}  && \alpha_{\rho(\ast)} \circ (F\circ \rho)(\id_{\ast})
            } 
        \end{equation} where the upper diagram commutes because of the lax unity of $\alpha$, and the lower square commutes because of the naturality of $\alpha$. 

        \item \textbf{Lax Naturality}

        For any $1$-morphisms $\ast\xrightarrow{f} \ast \xrightarrow{g} \ast $ in $B\huaG_{\bullet}$, 
        \begin{equation}
         \mbox{\tiny    \xymatrix{ F'(\rho(g)) \circ (\alpha_{\rho(\ast)} \circ F(\rho(f)) )  & (F'(\rho(g)) \circ \alpha_{\rho(\ast)} )\circ F(\rho(f))  \ar[l]_{a}  \ar[r]^-{\alpha_{\rho(g)}\horiprod 1_{F(\rho(f))} } & (\alpha_{\rho(\ast)} \circ F(\rho(g))  ) \circ F(\rho(f)) \ar[d]^-{a} \\
            F'(\rho(g))\circ  \bigl(F'(\rho(f)) \circ \alpha_{\rho(\ast)} \bigr) \ar[u]^{1_{F'} \horiprod \alpha_{\rho(f)} } &\circlearrowright& \alpha_{\rho(\ast)} \circ \bigl( F(\rho(g)) \circ F(\rho(f))  \bigr) \ar[d]^-{1_{\alpha_{\rho(\ast)}} \horiprod F^2_{\rho(g), \rho(f)} } \\
             \bigl( F'(\rho(g))\circ  F'(\rho(f)) \bigr)  \circ \alpha_{\rho(\ast)} \ar[u]^-{a} \ar[r]^-{ F^{'2}_{\rho(g), \rho(f)} \horiprod 1_{\alpha_{\rho(\ast)}} }  & F'(\rho(g)\circ\rho(f)) \circ \alpha_{\rho(\ast)} \ar[r]^-{\alpha_{\rho(g)\circ\rho(f)}}
            \ar[d]_-{F'(\rho^2_{g, f}) \horiprod 1_{\alpha_{\rho(\ast)}} } & \alpha_{\rho(\ast)}\circ F(\rho(g)\circ \rho(f)) 
            \ar[d]^-{1_{\alpha_{\rho(\ast)}} \horiprod F(\rho^2_{g,f}) } \\
            & F'(\rho(g\horicirc f)) \circ \alpha_{\rho(\ast)} \ar[r]_-{\alpha_{\rho(g\horicirc f)} } \ar@{}[ru]|{\circlearrowright}& \alpha_{\rho(\ast)} \circ F(\rho(g\horicirc f))
             } }
        \end{equation} where the upper diagram commutes because of the lax naturality of $\alpha$, and lower square commutes because of the naturality of $\alpha$.  
    \end{enumerate}
\end{itemize} In addition, since $\alpha$ is a strong transformation, so is $\rho^*\alpha$.

\item the laxity constraints: \begin{itemize}
    \item 
Note  $\rho^{*}\alpha \horicirc \rho^{*}\alpha' = \rho^*(\alpha\horicirc \alpha') $. Thus, $\rho^{* 2}_{\alpha, \alpha'}$ is defined to be the identity. 


\item 
Note that $\rho^*(\id_F)$ is exactly the identity transformation $\id_{F\circ \rho}$, which is explicitly defined in \cite[Proposition 4.2.12]{JY:Bicat}.
Thus, \[\rho^{* 0}_{F}: \id_{\rho^*(F)} \Rightarrow  \rho^*(\id_F) \] is defined to be the identity.

     \end{itemize} 
     
     \item 
    The lax naturality, lax left and right unity of $\rho^*$ defined above follow immediately.  \end{itemize}

     Similarly, we can define the strict functors \begin{align*}
         \iota^*: \quad  &2\Rep ( \huaG_{\bullet} ) \longrightarrow 2\Rep ( \huaG_{\bullet}' ), \\ (\iota\circ \rho) ^*: \quad & 2\Rep (\huaG_{\bullet} )\longrightarrow 2\Rep ( \huaG_{\bullet} ), \\
         (\rho\circ \iota)^*: \quad &2\Rep (\huaG_{\bullet}') \longrightarrow 2\Rep (\huaG_{\bullet}' ).
     \end{align*} By the strict associativity and strict unitality of pseudofunctors \cite[Thm.4.1.30]{JY:Bicat}, we have the equalities 
     \[ (\iota\circ \rho) ^* = \rho^* \circ \iota^*; \quad (\rho\circ \iota)^* = \iota^* \circ \rho^*;\] and for any $\huaG_{\bullet}$, 
     \[ (\id_{B\huaG_{\bullet}})^* = \id_{2\Rep(\huaG_{\bullet})}:2 \Rep (\huaG_{\bullet}) \longrightarrow 2\Rep ( \huaG_{\bullet} ).  \]

\bigskip

\textbf{Part II: the construction of the strong transformation $\tilde{\eta}$ and $\tilde{\eta}'$}

     Next we 
     define a strong transformation $$\tilde{\eta}: \id_{2\Rep (\huaG_{\bullet} )} \Rightarrow (\iota\circ \rho)^*.$$ 
     For any object $F$ in $2\Rep(\huaG_{\bullet})$, we construct a strong transformation $\tilde{\eta}_F: F\Rightarrow F\circ \iota\circ \rho$ associated to it, which is the post-whiskering of $\eta$ with $F$, as defined in \cite[Def. 11.1.1 (2)]{JY:Bicat}. 
     \begin{itemize}
         \item $(\tilde{\eta}_F)_{\ast} : F(\ast) \rightarrow F\circ\iota\circ \rho(\ast)$ is defined to be $F(\eta_{\ast})$. 
         \item For each $g$ in $\huaG_0$, $(\tilde{\eta}_F)_g$ is defined to be the composition
         \[ (F\circ \iota\circ \rho) (g) \circ (\tilde{\eta}_F)_{\ast} \xrightarrow{F^2_{(\iota\circ \rho) (g),  \eta_{\ast} } }
         F( (\iota\circ \rho) (g) \horicirc \eta_{\ast})\xrightarrow{F(\eta_g)} F(\eta_{\ast}\horicirc g)
         \xrightarrow{(F^2_{\eta_{\ast}, g})^{-1}} (\tilde{\eta}_F)_{\ast} \circ F(g).\]
     \end{itemize} Since $\eta$ is a strong transformation, so is $\tilde{\eta}_F$.

     In addition, for any strong transformation $\xymatrix@1@C+2em{
  B\huaG_{\bullet}  \ar@/^1pc/[r]^{F}_{}="4"
  \ar@/_1pc/[r]_{F'}^{}="5"
  \ar@{=>}"4";"5"^{\alpha}
  & 2\Vect
}$, we define an invertible modification $\tilde{\eta}_{\alpha}:  \bigl(  (\iota\circ \rho)^*\alpha \bigr) \horicirc \tilde{\eta}_{F}   \longrightarrow  \tilde{\eta}_{F'} \horicirc \id^*(\alpha)$  by the invertible $2$-morphism 
\[ (\tilde{\eta}_{\alpha})_{\ast}: = \alpha_{\eta_{\ast}}. \]

We check the modification axiom for $(\tilde{\eta}_{\alpha})_{\ast}$ explicitly below.

\begin{align*}
     & \xymatrix{ F(\ast) \ar@/_4.0pc/[rdd] \ar[r]^{F(g)} \ar[d]^{\alpha_{\ast}} & F(\ast)  \ar@/^6.0pc/[rd]\ar[d]|{\alpha_{\ast}} \ar[r]^{F(\eta_{\ast})} & F(\ast)  \ar[d]_{\alpha_{\ast} }   \\
     F'(\ast)  \ar@{}@<-3.5ex>[rr]|{\Uparrow{(\tilde{\eta}_{F'})_g} } \ar[r]_{F'(g)} \ar@/_1pc/[rd]|{F'(\eta_{\ast})} \ar@{}[ru]|{\buildrel{\alpha_g}\over\Rightarrow }& F'(\ast)  \ar[r]_{F'(\eta_{\ast})} \ar@{}[ru]|{\buildrel{\alpha_{\eta_{\ast}}}\over\Rightarrow }& F'(\ast)  \\
     & F'(\ast)  \ar@/_1pc/[ru]|{(F'\circ (\iota\circ \rho) )(g)  }& }  
     &\buildrel{(1)}\over{=}
    &   \xymatrix{ F(\ast)  \ar@/_4.0pc/[rdddddd]  \ar[r]^{F(g)} \ar[d]^{\alpha_{\ast}} & F(\ast)  \ar@/^6.0pc/[rd]   \ar[d]|{\alpha_{\ast}} \ar[r]^{F(\eta_{\ast})} & F(\ast)  \ar[d]_{\alpha_{\ast}}   \\
     F'(\ast)  \ar@/_2pc/[rddddd]|{F'(\eta_{\ast})}
     \ar@/_4pc/[rr]|{F'(\eta_{\ast}\horicirc g)} \ar@{}@<-5.0ex>[rr]|{ \Downarrow{F^{'2}_{\eta_{\ast}, g}} }
     \ar@/_8pc/[rr]|{F'( (\iota\circ \rho)(g)\horicirc  \eta_{\ast})} \ar@{}@<-14.0ex>[rr]|{ \Uparrow{F'(\eta_g)} }
     \ar[r]_{F'(g)}  \ar@{}[ru]|{\buildrel{\alpha_g}\over\Rightarrow }  \ar@{}@<-22.5ex>[rr]|{\Uparrow F^{'2}_{(\iota\circ \rho)(g),  \eta_{\ast}} }
     & F'(\ast)  \ar[r]_{F'(\eta_{\ast})} \ar@{}[ru]|{\buildrel{\alpha_{\eta_{\ast}}}\over\Rightarrow }
     & F'(\ast)  
      \\
     &&\\  &&\\  &&\\  &&\\ 
     & F'(\ast)  \ar@/_2pc/[ruuuuu]|{(F'\circ (\iota\circ \rho) )(g)  } & }  \\
     \buildrel{(2)}\over{=}    &   
      \xymatrix{ & F(\ast)   \ar@/^4.0pc/[drd]  \ar@/^1pc/[rd]|{F(\eta_{\ast})} & \\
      F(\ast)  \ar@/_7.0pc/[rddd] \ar@/^1pc/[ru]|{F(g)} \ar@{}@<4.0ex>[rr]|{\Downarrow{F^2_{\eta_{\ast}, g}}}
      \ar[rr]|{F(\eta_{\ast}\horicirc g)} \ar[d]^{\alpha_{\ast}} 
         && F(\ast)  \ar[d]_{\alpha_{\ast}}   \\
     F'(\ast)  \ar@/_3pc/[rdd]|{F'(\eta_{\ast})}
     \ar@/_3pc/[rr]|{F'( (\iota\circ \rho)(g)\horicirc  \eta_{\ast})} \ar@{}@<-4.0ex>[rr]|{ \Uparrow{F'(\eta_g)} }
     \ar[rr]|{F'(\eta_{\ast}\horicirc g)}  \ar@{}[rru]|{\buildrel{\alpha_{\eta_{\ast}\horicirc g}}\over\Rightarrow }  
     \ar@{}@<-10.0ex>[rr]|{\Uparrow F^{'2}_{(\iota\circ \rho)(g),  \eta_{\ast}} }
     && F'(\ast) 
      \\
     &&\\ 
     & F'(\ast)  \ar@/_3pc/[ruu]|{(F'\circ (\iota\circ \rho) )(g)  } & } 
      &\buildrel{(3)}\over{=}  &  
       \xymatrix{ & F(\ast)   \ar@/^4.0pc/[rddd]  \ar@/^1pc/[rdd]|{F(\eta_{\ast})} & \\ && \\ 
      F(\ast)  \ar@/_5.0pc/[rdd]   \ar@/^1pc/[ruu]|{F(g)} \ar@{}@<+10.0ex>[rr]|{\Downarrow{F^2_{\eta_{\ast}, g}}}
      \ar@/^3pc/[rr]|{F(\eta_{\ast}\horicirc g)} \ar@{}@<4.0ex>[rr]|{\Uparrow{F(\eta_g)}}
      \ar[rr]|{ F((\iota\circ \rho)(g)\horicirc \eta_{\ast}) } \ar[d]^{\alpha_{\ast}} 
         && F(\ast)  \ar[d]_{\alpha_{\ast}}   \\
     F'(\ast)  \ar@/_1.5pc/[rd]|{F'(\eta_{\ast})}
     \ar[rr]|{F'( (\iota\circ \rho)(g)\horicirc  \eta_{\ast})}  
     \ar@{}[rru]|{\buildrel{\alpha_{(\iota\circ \rho)(g)\horicirc \eta_{\ast}}}\over\Rightarrow }  
     \ar@{}@<-4.0ex>[rr]|{\Uparrow F^{'2}_{(\iota\circ \rho)(g),  \eta_{\ast}} }
     && F'(\ast)  
      \\
     & F'(\ast)  \ar@/_1.5pc/[ru]|{(F'\circ (\iota\circ \rho) )(g)  } & } \\ 
       \buildrel{(4)}\over{=}  &   
  \xymatrix{ 
     &&F(\ast)   \ar@/^6.0pc/[rrddddd]  \ar@/^2pc/[rrdddd]|{(F (\eta_{\ast}) } &&\\  & &&\\  &&&\\  &&&\\ 
  F(\ast)  \ar@/_6.0pc/[rrd]  
  \ar@/^2pc/[rruuuu]|{F(g)}
     \ar@/^4pc/[rrrr]|{F'(\eta_{\ast}\horicirc g)} \ar@{}@<+22.0ex>[rrrr]|{ \Downarrow{F^{2}_{\eta_{\ast}, g}} }
     \ar@/^8pc/[rrrr]|{F'( (\iota\circ \rho)(g)\horicirc  \eta_{\ast})} \ar@{}@<+14.0ex>[rrrr]|{ \Uparrow{F(\eta_g)} }
      \ar@{}@<+5.0ex>[rrrr]|{ \Uparrow{F^2_{(\iota\circ \rho)(g), \eta_{\ast}}} }
     \ar[rr]^{F(\eta_{\ast})} \ar[d]^{\alpha_{\ast}} && F(\ast)  \ar[d]|{\alpha_{\ast}} \ar[rr]^{F((\iota\circ \rho)(g))} & &F(\ast)  \ar[d]_{\alpha_{\ast}}   \\
     F'(\ast)  
     \ar[rr]_{F'(\eta_{\ast})}  \ar@{}[rru]|{\buildrel{\alpha_{\eta_{\ast}}}\over\Rightarrow } 
     & & F'(\ast)  \ar[rr]_{F'((\iota\circ \rho)(g) )} \ar@{}[rur]|{\buildrel{\alpha_{(\iota\circ \rho)(g)}}\over\Rightarrow }
     & & F'(\ast)  
      \\ }  
       &\buildrel{(5)}\over{=} &   \xymatrix{ 
     && F(\ast) \ar@/^5.0pc/[rrdd]  \ar@/^1pc/[rrd]|{(F (\eta_{\ast}) } &&\\ 
  F(\ast)  \ar@/_7.0pc/[rrd] 
  \ar@/^1pc/[rru]|{F(g)}
    \ar@{}@<+4.0ex>[rrrr]|{ \Uparrow{ (\tilde{\eta}_F)_g} }
     \ar[rr]^{F(\eta_{\ast})} \ar[d]^{\alpha_{\ast}} && F(\ast)  \ar[d]|{\alpha_{\ast}} \ar[rr]^{F((\iota\circ \rho)(g))} & & F(\ast)  \ar[d]_{\alpha_{\ast}}   \\
     F'(\ast)  
     \ar[rr]_{F'(\eta_{\ast})}  \ar@{}[rru]|{\buildrel{\alpha_{\eta_{\ast}}}\over\Rightarrow } 
     & & F'(\ast)  \ar[rr]_{F'((\iota\circ \rho)(g) )} \ar@{}[rur]|{\buildrel{\alpha_{(\iota\circ \rho)(g)}}\over\Rightarrow }
     & & F'(\ast)  
      \\ }  
\end{align*}
     
     In the equations of pasting diagrams above, the equality $(1)$ and $(5)$ are defined from the definition of $(\tilde{\eta}_{F'})_g$ and $(\tilde{\eta}_{F})_g$ respectively; the equality $(2)$ and $(4)$ are defined from the lax naturality  of $\alpha$; the equality $(3)$ is defined from the naturality of $\alpha$.

Next we check the naturality, lax unity, lax naturality of $\tilde{\eta}$, so that the data above indeed define a strong transformation.

\begin{enumerate}
    \item \textbf{Naturality}

    For any $\huaG_{\bullet}$-representation $F$ and $F'$ and any strong transformation $\alpha, \beta: F\Rightarrow F'$, let $\Gamma$ denote a modification from $\alpha$ to $\beta$. The modification axiom of $\Gamma$ guarantees the equality \eqref{eta_nat_1}.
    \begin{equation}\label{eta_nat_1}
    \xymatrix{F(\ast)   \ar@/^0.8pc/[rr]^{\alpha_{\ast}}  \ar@/_0.8pc/[rr]_{\beta_{\ast}}
    \ar@{}[rr]|{\Uparrow{\Gamma_{\ast}}}    \ar[dd]_{F(\eta_{\ast})} &&F'(\ast)    \ar[dd]^{F'(\eta_{\ast}) \hspace{0.6cm} = } \\
    && \\
    F\circ \iota\circ \rho(\ast) \ar@/_0.8pc/[rr]_{\beta_{\iota\circ \rho(\ast)} } \ar@{}[rruu]|{\buildrel{\beta_{\eta_{\ast}}}\over\Rightarrow }
    && F'\circ\iota\circ\rho(\ast)
    } \xymatrix{F(\ast) \ar@/^0.8pc/[rr]^{\alpha_{\ast}} \ar[dd]_{F(\eta_{\ast})}  && F'(\ast)  \ar[dd]^{F'(\eta_{\ast})} \\ && \\
    F\circ\iota\circ \rho (\ast)   \ar@/^0.8pc/[rr]^{\alpha_{\iota\circ \rho (\ast)}}  \ar@/_0.8pc/[rr]_{\beta_{\iota\circ \rho (\ast)}}
    \ar@{}[rr]|{\Uparrow{\Gamma_{\iota\circ \rho (\ast)}}}  \ar@{}[rruu]|{\buildrel{\alpha_{\eta_{\ast}}}\over\Rightarrow }    &&F'\circ\iota\circ \rho (\ast)     }
    \end{equation}
    which gives the naturality of $\tilde{\eta}$, i.e.  \eqref{eta_nat_2}.
    \begin{equation} \label{eta_nat_2}
           \xymatrix{F    \ar@/^0.8pc/[rr]^{ \id^*\alpha}  \ar@/_0.8pc/[rr]_{\id^*\beta}
    \ar@{}[rr]|{\Uparrow{\id^*\Gamma}}    \ar[dd]_{\tilde{\eta}_{F}} &&F'   \ar[dd]^{\tilde{\eta}_{F'} \hspace{0.6cm} = } \\
    && \\
    F\circ \iota\circ \rho  \ar@/_0.8pc/[rr]_{(\iota\circ \rho)^*{\beta} } \ar@{}[rruu]|{\buildrel{\tilde{\eta}_{\beta}}\over\Rightarrow }
    && F'\circ\iota\circ\rho 
    } \xymatrix{F  \ar@/^0.8pc/[rr]^{\id^{*}\alpha} \ar[dd]_{\tilde{\eta}_{F}}  && F'  \ar[dd]^{\tilde{\eta}_{F'}} \\ && \\
    F\circ\iota\circ \rho     \ar@/^0.8pc/[rr]^{(\iota\circ \rho)^*\alpha}  \ar@/_0.8pc/[rr]_{(\iota\circ \rho)^*\beta}
    \ar@{}[rr]|{\Uparrow{(\iota\circ \rho)^*\Gamma}}  \ar@{}[rruu]|{\buildrel{\tilde{\eta}_{\alpha}}\over\Rightarrow }    &&F'\circ\iota\circ \rho     }
    \end{equation}

    \item \textbf{Lax Unity}

   Since    $(\iota\circ\rho)^{* 0}_{\ast}$ and $\id^{* 0}_{\ast} $ are both identity,   $(\tilde{\eta}_{1_F})_{\ast} = (1_F)_{\eta_{\ast}} = l^{-1}_{F(\eta_{\ast})}\circ r_{F(\eta_{\ast}) }$,  and $(\tilde{\eta}_F)_{\ast} = F(\eta_{\ast})$,   we have the equation
    \begin{equation}\label{eta_unit_1}
        \xymatrix{ F \ar@/^0.8pc/[rr]^{\id^{*}(1_F)} \ar[dd]_{\tilde{\eta}_F} && F \ar[dd]^{\tilde{\eta}_F \hspace{0.6cm} = } \\ && \\
        F\circ \iota\circ \rho  \ar@{}[rruu]|{\buildrel{\tilde{\eta}_{1_F}}\over\Leftarrow }
        \ar@/^0.8pc/[rr]^{(\iota\circ\rho)^* 1_F}   \ar@/_0.8pc/[rr]_{1_{F\circ \iota\circ \rho}}
        \ar@{}[rr]|{\Uparrow {(\iota\circ\rho)^{* 0}_{\ast}} } && F\circ\iota\circ\rho }
             \xymatrix{ F \ar@/^0.8pc/[rr]^{\id^{*}(1_F)}  \ar@/_0.8pc/[rr]_{1_{F\circ \id }}
             \ar[dd]_{\tilde{\eta}_F} \ar@{}@<+5.0ex>[dd]|{ \buildrel{l}\over\Rightarrow} 
               \ar@{}[rr]|{\Uparrow {\id^{* 0}_{\ast}} }     \ar@/_0.6pc/[rrdd] && F \ar[dd]^{\tilde{\eta}_F} 
             \ar@{}@<-5.0ex>[dd]|{ \buildrel{r^{-1}}\over\Rightarrow} 
               \\ && \\
        F\circ \iota\circ \rho  
       \ar@/_0.8pc/[rr]_{1_{F\circ \iota\circ \rho}}
       && F\circ\iota\circ\rho }
    \end{equation} 
 
\item \textbf{Lax Naturality}

Let $F$, $F'$ and $F''$ be $\huaG_{\bullet}$-representations and $F\buildrel{\alpha}\over\Rightarrow F' \buildrel{\alpha'}\over\Rightarrow F''$ be strong transformations. Since $\tilde{\eta}_{ \alpha} = \alpha_{\eta_{\ast}}$, the lax naturality of $\alpha$ guarantees the equality \eqref{eta_lax_nat_1}.

\begin{equation} \label{eta_lax_nat_1}
\xymatrix{ 
 F \ar@/^0.8pc/[rr]^{ \id^*(  \alpha'\horicirc \alpha) } \ar[dd]_{\tilde{\eta}_F} && F'' \ar[dd]^{\tilde{\eta}_{F''} \hspace{0.6cm} = } \\ && \\
   F\circ \iota\circ \rho  \ar@{}[rruu]|{\buildrel{\tilde{\eta}_{ \alpha'\horicirc \alpha}}\over\Rightarrow }
        \ar@/^0.8pc/[rr]|{(\iota\circ\rho)^*( \alpha'\horicirc \alpha) }   \ar@/_0.8pc/[rd]_{(\iota\circ\rho)^*\alpha}
      && F''\circ\iota\circ\rho \\
        & F' \circ \iota\circ \rho  \ar@/_0.8pc/[ru]_{(\iota\circ\rho)^*\alpha'}   \ar@{}[u]|{\Uparrow {(\iota\circ\rho)^{* 2}_{\alpha', \alpha}} } & 
}
\xymatrix{ 
 F \ar@/^0.8pc/[rr]^{ \id^*(  \alpha'\horicirc \alpha) }\ar@/_0.8pc/[rd]|{ \id^*(  \alpha) }  \ar[dd]_{\tilde{\eta}_F} & & F'' \ar[dd]^{\tilde{\eta}_{F''} } \\ & F' \ar@/_0.8pc/[ru]|{ \id^*(  \alpha') }  \ar[dd]|{\tilde{\eta}_{F'} }  \ar@{}[u]|{ \Uparrow \id^{* 2}_{\alpha, \alpha'} } & \\ 
   F\circ \iota\circ \rho  \ar@{}[ru]|{\buildrel{\tilde{\eta}_{ \alpha}}\over\Rightarrow }
      \ar@/_0.8pc/[rd]_{(\iota\circ\rho)^*\alpha}
     && F''\circ\iota\circ\rho \\   &   F\circ \iota\circ \rho    \ar@{}[ruuu]|{\buildrel{\tilde{\eta}_{ \alpha'}}\over\Rightarrow } \ar@/_0.8pc/[ru]_{(\iota\circ\rho)^*\alpha' } &
}
\end{equation}

\end{enumerate}

Thus, $\tilde{\eta}$ is a strong transformation.

\bigskip

Similarly, we can define a strong transformation \[ \tilde{\eta}': (\iota\circ \rho)^* \Rightarrow \id_{2\Rep ( \huaG_{\bullet} )}.\]
For any object $F$ in $2\Rep (\huaG_{\bullet})$,   we construct a strong transformation $\tilde{\eta}'_F: F \circ \iota\circ \rho\Rightarrow F  $  associated to it, which is the post-whiskering of $\eta'$ with $F$, as defined in \cite[Def. 11.1.1 (2)]{JY:Bicat}. 
     \begin{itemize}
         \item $(\tilde{\eta}'_F)_{\ast} : F\circ\iota\circ \rho(\ast) \rightarrow F(\ast)$ is defined to be $F(\eta'_{\ast})$. 
         \item For each $g$ in $\huaG_0$, $(\tilde{\eta}_F)_g$ is defined to be the composition
         \[ F (g) \circ (\tilde{\eta}'_F)_{\ast} \xrightarrow{F^2_{ g,  \eta'_{\ast} } }
         F(  g \horicirc \eta'_{\ast})\xrightarrow{F(\eta'_g)} F(\eta'_{\ast}\horicirc (\iota\circ \rho)g)
         \xrightarrow{(F^2_{\eta'_{\ast}, \iota\circ \rho(g)})^{-1}} (\tilde{\eta}'_F)_{\ast} \circ \bigl( F\circ\iota\circ \rho\bigr) (g).\]
     \end{itemize} Since $\eta'$ is a strong transformation, so is $\tilde{\eta}'_F$.

     In addition, for any strong transformation $\xymatrix@1@C+2em{
  B\huaG_{\bullet}  \ar@/^1pc/[r]^{F}_{}="4"
  \ar@/_1pc/[r]_{F'}^{}="5"
  \ar@{=>}"4";"5"^{\alpha}
  & 2\Vect
}$, we define an invertible modification $\tilde{\eta}'_{\alpha}:  \bigl(  (\iota\circ \rho)^*\alpha \bigr) \horicirc \tilde{\eta}'_{F}   \longrightarrow  \tilde{\eta}'_{F'} \horicirc \id^*(\alpha)$  by the $2$-isomorphism 
\[ (\tilde{\eta}'_{\alpha})_{\ast}: = \alpha_{\eta'_{\ast}}. \]

     In addition, for any strong transformation $\xymatrix@1@C+2em{
  B\huaG_{\bullet}  \ar@/^1pc/[r]^{F}_{}="4"
  \ar@/_1pc/[r]_{F'}^{}="5"
  \ar@{=>}"4";"5"^{\alpha}
  & 2\Vect
}$, we define an invertible modification $\tilde{\eta}'_{\alpha}:  \id^*(\alpha)  \horicirc \tilde{\eta}'_{F}   \longrightarrow  \tilde{\eta}'_{F'} \horicirc \bigl( (\iota\circ \rho)^*(\alpha) \bigr)$  by the  $2$-isomorphism
\[ (\tilde{\eta}'_{\alpha})_{\ast}: = \alpha_{\eta'_{\ast}}. \]
The modification axiom for $(\tilde{\eta}'_{\alpha})_{\ast}$ can be checked similarly to that for $(\tilde{\eta}_{\alpha})_{\ast}$.
In addition, we can prove the naturality, lax unity, lax naturality of $\tilde{\eta}'$ as the proof for those of $\tilde{\eta}$. Thus, $\tilde{\eta}'$ is a strong transformation.

\bigskip

\textbf{Part III: the construction of the invertible modification $\tilde{\Gamma}$ and $\tilde{\Gamma}'$}

Next, we show there is an invertible modification $\tilde{\Gamma}: 1_{\id_{2\Rep(\huaG_{\bullet})}} \longrightarrow \tilde{\eta}'\horicirc \tilde{\eta}$. For each $\huaG_{\bullet}$-representation $F$, we need to define a $2$-isomorphism \[ \tilde{\Gamma}_F: 1_F\rightarrow \tilde{\eta}'_F\horicirc \tilde{\eta} _F, \] which is a modification itself. It should be given by a $2$-isomorphism in $2\Vect$ \[ ( \tilde{\Gamma}_F)_{\ast}: 1_{F(\ast)} \longrightarrow  (\tilde{\eta}'_F)_{\ast}\circ (  \tilde{\eta}_F)_{\ast} = F(\eta'_{\ast}) \circ F(\eta_{\ast}). \]  We define it
by the composition \[ 1_{F(\ast)}\xrightarrow{F^0_{\ast}} F1_{\ast}\xrightarrow{F(\Gamma_{\ast})} F(\eta'_{\ast}\circ \eta_{\ast}) \xrightarrow{(F^2_{\eta'_{\ast}, \eta_{\ast}})^{-1}} F(\eta'_{\ast})
\circ F(\eta_{\ast}), \] which is invertible. We check that $\tilde{\Gamma}_F$ satisfy the modification axiom below, i.e. 
for any $\huaG_{\bullet}$-representation $F$ and $F'$ and strong transformation $f: F\Rightarrow F'$, we need to show the equation \eqref{modif_axiom_modif} holds. 
\begin{equation} \label{modif_axiom_modif}
  \xymatrix{ F    \ar@/^0.8pc/[rr]^{1_F}  \ar@/_0.8pc/[rr]_-{ \tilde{\eta}'_{F} \horicirc \tilde{\eta}_{F}} 
    \ar@{}[rr]|{\Downarrow{\tilde{\Gamma}_F} }    \ar[dd]_{f }  && F    \ar[dd]^{ f  \hspace{1cm} = \hspace{0.5cm}} \\
    && \\
    F'   \ar@/_0.8pc/[rr]_{ \tilde{\eta}'_{F'} \horicirc \tilde{\eta}_{F'} } \ar@{}[rruu]|{\buildrel{ (\tilde{\eta}'  \horicirc \tilde{\eta} )_f}\over\Rightarrow }
    &&  F' 
    } 
            \xymatrix{  F  
     \ar@/^1pc/[rr]^{1_{F }}
      \ar[d]_{f }   \ar@{}[rrd]|{\buildrel{ (1_{\id_{2\Rep(\huaG_{\bullet}) }})_f }\over\Leftarrow }  && F(\ast)  \ar[d]^{f}     \\
     F'     \ar[rr]|{1_{F' }}        
     \ar@{}@<-3.5ex>[rr]|{ \Downarrow{\tilde{\Gamma}_{F'}}}
     \ar@/_3pc/[rr]|{ \tilde{\eta}'_{F'} \horicirc  \tilde{\eta}_{F'}} 
     & & F'    } 
\end{equation} where $ (1_{\id_{2\Rep (\huaG_{\bullet}) }})_f = l^{-1}_{f } \circ r_{f} $.  
The equation \eqref{modif_axiom_modif} is true because we have the equation of pasting diagrams below.
\begin{align*}
     &          \xymatrix{ F (\ast)   \ar@/^0.8pc/[rr]^{ \id_{F(\ast)}}  \ar@/_0.8pc/[rr]_-{F(\eta'_{\ast}) \circ F(\eta_{\ast})} 
    \ar@{}[rr]|{\Downarrow{\tilde{(\Gamma}_F)_{\ast}}}    \ar[dd]_{f_{\ast} }  &&F(\ast)   \ar[dd]^{f_{\ast}  } \\
    && \\
    F'(\ast)  \ar@/_0.8pc/[rr]_{F'(\eta'_{\ast}) \circ F'(\eta_{\ast}) } \ar@{}[rruu]|{\buildrel{f_{\eta'_{\ast}} \horicirc f_{\eta_{\ast}}}\over\Rightarrow }
    &&  F'(\ast)
    } 
     &\buildrel{(1)}\over{=}
    &    \xymatrix{ 
 F(\ast)  
     \ar@/^11pc/[rrrr]|{\id_{F(\ast)}} \ar@{}@<+22.0ex>[rrrr]|{ \Downarrow{F^{0}_{\ast}} }
     \ar@/^8pc/[rrrr]|{ F(\id_{\ast})} \ar@{}@<+14.0ex>[rrrr]|{ \Downarrow{F(\Gamma_{\ast})} }
      \ar@/^4pc/[rrrr]|{ F(\eta'_{\ast}\circ \eta_{\ast})} 
      \ar@{}@<+5.0ex>[rrrr]|{ \Uparrow{F^2_{\eta'_{\ast}, \eta_{\ast}}} }
     \ar[rr]^{F(\eta_{\ast})} \ar[d]^{f_{\ast}} && F(\ast)  \ar[d]|{(\iota\circ \rho)^*f_{\ast}} \ar[rr]^{ F(\eta'_{\ast}) } & &F(\ast)  \ar[d]_{f_{\ast}}   \\
     F'(\ast)   \ar@/_2pc/[rrrr]
     \ar[rr]_{F'(\eta_{\ast})}  \ar@{}[rru]|{\buildrel{f_{\eta_{\ast}}}\over\Rightarrow } 
     & & F'\circ \iota\circ \rho (\ast)  \ar[rr]_{F'(\eta'_{\ast})} \ar@{}[rur]|{\buildrel{f_{\eta'_{\ast}}}\over\Rightarrow }
     & & F'(\ast)  
      \\ }   \\
      \buildrel{(2)}\over{=} &
      \xymatrix{  F(\ast)  
     \ar@/^5.5pc/[rr]|{\id_{F(\ast)}} \ar@{}@<+10.0ex>[rr]|{ \Downarrow{F^{0}_{\ast}} }
     \ar@/^3pc/[rr]|{ F(\id_{\ast})} \ar@{}@<+4.0ex>[rr]|{ \Downarrow{F(\Gamma_{\ast})} }
      \ar[rr]|{ F(\eta'_{\ast}\circ \eta_{\ast})} 
      \ar[d]^{f_{\ast}}   && F(\ast)  \ar[d]_{f_{\ast}}   \\
     F'(\ast)    
     \ar[rr]|{F'(\eta'_{\ast}\circ \eta_{\ast})} \ar@{}@<-4.0ex>[rr]|{\Uparrow{F^{'2}_{\eta'_{\ast}, \eta_{\ast}}}} \ar@{}[rru]|{\buildrel{f_{\eta'_{\ast}\circ \eta_{\ast}}}\over\Rightarrow }  \ar@/_4pc/[rr]|{ F'(\eta'_{\ast})\circ F'( \eta_{\ast})} 
     & & F'(\ast)   } 
      & \buildrel{(3)}\over{=} 
        &    \xymatrix{  F(\ast)  
     \ar@/^3pc/[rr]|{\id_{F(\ast)}} \ar@{}@<+4.0ex>[rr]|{ \Downarrow{F^{0}_{\ast}} }
     \ar[rr]|{ F(\id_{\ast})} 
      \ar[d]^{f_{\ast}}   && F(\ast)  \ar[d]_{f_{\ast}}   \\
     F'(\ast)    
     \ar[rr]|{F'(\id_{\ast})} \ar@{}@<-4.0ex>[rr]|{\Downarrow{F'(\Gamma_{\ast})}} \ar@{}[rru]|{\buildrel{f_{\id_{\ast}}}\over\Rightarrow }  \ar@/_3pc/[rr]|{ F'(\eta'_{\ast}\circ  \eta_{\ast})}      \ar@{}@<-10.0ex>[rr]|{\Uparrow{F^{'2}_{\eta'_{\ast}, \eta_{\ast}}}}
      \ar@/_6pc/[rr]|{ F'(\eta'_{\ast})\circ  F'(\eta_{\ast})} 
     & & F'(\ast)   } 
     \buildrel{(4)}\over{=} 
     \xymatrix{  F(\ast)  
     \ar@/^1pc/[rr]^{\id_{F(\ast)}}
      \ar[d]_{f_{\ast}} \ar@/^0.2pc/[rrd]  \ar@{}@<+4.0ex>[rrd]|{\buildrel{r^{-1}_{f_{\ast}}}\over\Rightarrow }  && F(\ast)  \ar[d]^{f_{\ast}}     \\
     F'(\ast)    \ar[rr]|{\id_{F'(\ast)}}      \ar@{}[ru]|{\buildrel{l_{f_{\ast}}}\over\Rightarrow }  
     \ar@{}@<-3.5ex>[rr]|{ \Downarrow{F^{'0}_{\ast}}}
     \ar@/_3pc/[rr]|{F'(\id_{\ast})}  \ar@{}@<-10.0ex>[rr]|{\Downarrow{F'(\Gamma_{\ast})}} 
 \ar@/_6pc/[rr]|{ F'(\eta'_{\ast}\circ  \eta_{\ast})}      \ar@{}@<-18.0ex>[rr]|{\Uparrow{F^{'2}_{\eta'_{\ast}, \eta_{\ast}}}}
      \ar@/_10pc/[rr]|{ F'(\eta'_{\ast})\circ  F'(\eta_{\ast})} 
     & & F'(\ast)   }  \\
      \buildrel{(5)}\over{=} &
         \xymatrix{  F(\ast)  
     \ar@/^1pc/[rr]^{\id_{F(\ast)}}
      \ar[d]_{f_{\ast}}   \ar@{}[rrd]|{\buildrel{  l^{-1}_{f_{\ast}} \circ r_{f_{\ast}} }\over\Leftarrow }  && F(\ast)  \ar[d]^{f_{\ast}}     \\
     F'(\ast)    \ar[rr]|{\id_{F'(\ast)}}        
     \ar@{}@<-3.5ex>[rr]|{ \Downarrow{(\tilde{\Gamma}_F)_{\ast}}}
     \ar@/_3pc/[rr]|{ F'(\eta'_{\ast})\circ  F'(\eta_{\ast})} 
     & & F'(\ast)   } 
      &&
\end{align*}

   In the equations of pasting diagrams above, the equality $(1)$ and $(5)$ are defined from the definition of $\tilde{(\Gamma}_F)_{\ast}$ and $\tilde{(\Gamma}_{F'})_{\ast}$ respectively; the equality $(2)$ is defined from the lax naturality of $f$; the equality $(3)$ is defined from the naturality of $\alpha$;   and $(4)$ are defined from the lax unity  of $f$. 

   Thus, the construction above indeed defines an invertible modification $\tilde{\Gamma}: 1_{\id_{2\Rep(\huaG_{\bullet}) }} \longrightarrow \tilde{\eta}'\horicirc \tilde{\eta}$.

Moreover, we can construct an invertible modification \[\tilde{\Gamma}': \tilde{\eta}\horicirc \tilde{\eta}' \longrightarrow  1_{(\iota \circ \rho)^*}. \]

\textbf{Part IV: the construction of the invertible modification $\tilde{\Chi}$ and $\tilde{\Chi}'$}

\begin{enumerate}
    \item 
Analogous to the construction of $\tilde{\eta}'$, we can construct a strong transformation \[ \tilde{\epsilon}:   (\rho\circ \iota)^*  \Rightarrow \id_{2\Rep (\huaG_{\bullet}')};  \] 
\item 
analogous to the construction of $\tilde{\eta}$, we can construct a strong transformation 
 \[ \tilde{\epsilon}':     \id_{2\Rep (\huaG_{\bullet}')}  \Rightarrow (\rho\circ \iota)^*.\] 
\item analogous to the construction of $\tilde{\Gamma}$ and $\tilde{\Gamma}'$, we can construct the invertible modifications
      \[\tilde{\Chi}: 1_{(\rho\circ \iota)^*}\longrightarrow  \tilde{\epsilon}' \horicirc \tilde{\epsilon}, \quad\text{and}\quad  \tilde{\Chi}': \tilde{\epsilon}\circ \tilde{\epsilon}'\longrightarrow 
     1_{\id_{2\Rep(\huaG_{\bullet}')}} . \]
     \end{enumerate}


     \bigskip 
     

     Thus, we have constructed a biequivalence between $2\Rep( \huaG_{\bullet}) $ and $2\Rep (\huaG_{\bullet}')$, which is given by $$(  2\Rep(\huaG_{\bullet})\xrightarrow{ \iota^*} 2\Rep (\huaG_{\bullet}'), \quad  2\Rep (\huaG_{\bullet}')\xrightarrow{ \rho^*} 2\Rep (\huaG_{\bullet}),  \quad  \id_{2\Rep (\huaG_{\bullet}) }\buildrel{\tilde{\eta}}\over\Rightarrow (\iota\circ\rho)^*, \quad  (\rho\circ \iota)^* \buildrel{\tilde{\epsilon}}\over \Rightarrow \id_{2\Rep(\huaG_{\bullet}')} ).$$

\end{proof}

\section*{Acknowledgment}

The  author is supported by the Young Scientists Fund of the National Natural Science Foundation of China (Grant No. 11901591) for the project ``Quasi-elliptic cohomology and its application in geometry and topology",  the General Program of the National Natural Science Foundation of China (Grant No. 12371068) for the project ``Quasi-elliptic cohomology and its application in topology and mathematical physics",   and the research funding from Huazhong University of Science and Technology.  

This work was initiated while the author was visiting the University of G\"{o}ttingen, hosted by Professor Chenchang Zhu. The author thanks the DAAD-K.C.Wong Postdoctoral Fellowship for hospitality and support, and Professor Zhu for her kind and enlightening advice. She introduced the author to coherent 2-groups, string 2-groups, and higher groupoids, and her insights significantly influenced the exposition of Sections 2, 3.1, and 4.1. The author also thanks Professor Fei Han, who shared the original idea of using loop group positive energy representations with the author.

The author acknowledges the use of AI-assisted tools for language polishing and organizational suggestions. All mathematical content, verification, and critical evaluation are solely the responsibility of the author.


\bibliographystyle{alpha}
\bibliography{bibz}
\end{document}